\documentclass{extarticle} 

\usepackage[utf8]{inputenc}	
\usepackage[english]{babel}
\usepackage{amsmath,amsfonts,amssymb,amsthm,mathtools} 

\usepackage{geometry}
\newgeometry{vmargin={3cm}, hmargin={30mm, 30mm}}

\usepackage{graphicx}	
\usepackage{subcaption}
\usepackage{wrapfig}  
\usepackage{float} 

\usepackage{makeidx} 

\usepackage[colorinlistoftodos]{todonotes}
\usepackage[colorlinks=true, allcolors=blue]{hyperref}

\definecolor{skyblue1}{rgb}{0.447,0.624,0.812}
\definecolor{scarletred1}{rgb}{0.937,0.561,0.561}

\usepackage{listings} 

\definecolor{codegreen}{rgb}{0,0.6,0}
\definecolor{codegray}{rgb}{0.5,0.5,0.5}
\definecolor{codepurple}{rgb}{0.58,0,0.82}
\definecolor{backcolour}{rgb}{0.95,0.95,0.92}
 
\lstdefinestyle{mystyle}{
    backgroundcolor=\color{backcolour},   
    commentstyle=\color{codegreen},
    keywordstyle=\color{magenta},
    numberstyle=\tiny\color{codegray},
    stringstyle=\color{codepurple},
    basicstyle=\footnotesize,
    breakatwhitespace=false,         
    breaklines=true,                 
    captionpos=b,                    
    keepspaces=true,                 
    numbers=left,                    
    numbersep=5pt,                  
    showspaces=false,                
    showstringspaces=false,
    showtabs=false,                  
    tabsize=2
}

\usepackage{color}
\usepackage{hyperref} 
\hypersetup{colorlinks,
citecolor = black,
linkcolor = black,
urlcolor = black}

\usepackage{siunitx}  

\usepackage{pdfpages}	
\usepackage{lscape} 

\usepackage{footnote}
\usepackage{sidecap}

\usepackage{parskip}

\usepackage{tabularx}
\makesavenoteenv{tabular}
\makesavenoteenv{table}
\usepackage{array}
\usepackage{booktabs}

\usepackage{colortbl}

\usepackage{tikz}
\usetikzlibrary{shapes.geometric,calc,positioning,3d,intersections,arrows}
\usetikzlibrary{decorations.markings}
\usepackage{pgfplots}
\usepackage{circuitikz}  

\usepackage{enumerate}

\usepackage{setspace}

\usepackage[skins,theorems]{tcolorbox}
\tcbset{highlight math style={enhanced,
  colframe=red,colback=white,arc=0pt,boxrule=1pt}}

\usepackage{float}
\restylefloat{table}

\usepackage{esvect}
\usepackage{bm}

\newcommand{\rom}[1]{\uppercase\expandafter{\romannumeral #1\relax}}

\newcommand{\romLower}[1]{\romannumeral #1\relax}

\newcommand{\hsone}{\hspace{1cm}}
\newcommand{\hstwo}{\hspace{2cm}}
\newcommand{\vsone}{\vspace{1cm}}

\newcommand{\tx}[1]{\text{#1}}

\newcommand{\parD}[2]{\frac{\partial #1}{\partial #2}}

\newcounter{region}
\newenvironment{region}{\refstepcounter{region}\equation}{\tag{\rom{\theregion}}\endequation}

\newcounter{LowerCaseRoman}
\newenvironment{equationTaggedLowerCaseRoman}{\refstepcounter{LowerCaseRoman}\equation}{\tag{\romLower{\theLowerCaseRoman}}\endequation}

\usepackage[style=numeric, sorting=nyt]{biblatex}
\numberwithin{equation}{section}

\usepackage[page,toc,titletoc,title]{appendix}

\pgfplotsset{compat=1.18} 

\makeatletter
\renewcommand*\l@section[2]
  {%
    \ifnum \c@tocdepth >\z@
      \addpenalty \@secpenalty
      \addvspace {1.0em \@plus \p@ }%
      \setlength \@tempdima {1.5em}%
      \begingroup
        \parindent \z@
        \rightskip \@pnumwidth
        \parfillskip -\@pnumwidth
        \leavevmode \bfseries
        \advance \leftskip \@tempdima
        \hskip -\leftskip
        #1\nobreak
        \hfil
        \nobreak
        \hb@xt@ \@pnumwidth {\hss #2\kern -\p@ \kern \p@ }%
        \par
      \endgroup
    \fi
  }
\makeatother

\usepackage{afterpage}

\title{\vspace{-1.5cm} Effects of Symmetry in a Diffusive Energy Balance Model\\
    \vspace{0.5cm}
   \normalsize }

\usepackage{authblk}

\author[1]{Aksel Samuelsberg}
\author[1]{Per Kristen Jakobsen}
\affil[1]{Department of Mathematics and Statistics, UiT - The Arctic University of Norway}

\date{\textbf{Correspondence:} Aksel Samuelsberg, aksel.samuelsberg@uit.no}

\begin{document}

\maketitle   

\section*{Abstract}

\textbf{
In this paper, we solve a North-type Energy Balance Model (EBM) using an analytical method, the Boundary Integral Method. This approach is discussed in light of existing analytical techniques for this type of equation. We use the method to demonstrate that the placement of a zonally symmetric continent, with an altered ice-albedo feedback dynamic, introduces new equilibrium states. A finite difference algorithm is implemented to solve the time-dependent equation and assess the stability of the equilibrium states, along with a numerical perturbation scheme. Bifurcation diagrams are drawn and we show that the bifurcation curve is extremely sensitive to the placement of a continent. The continent is initially configured with meridional symmetry, and we investigate how the system dynamics respond to a gradual reduction of the system's symmetry properties. We find that meridional symmetry increases the number of fold bifurcations and equilibria.
Additionally, we discuss how the emerging bifurcation structures may provide insights into the complex dynamics involved as one ascends the climate model hierarchy.
}


\tikzset{
  >=latex,
  inner sep=0pt,
  outer sep=2pt,
  mark coordinate/.style={inner sep=0pt,outer sep=0pt,minimum size=3pt, fill=black,circle}
}

\newcommand\pgfmathsinandcos[3]{%
  \pgfmathsetmacro#1{sin(#3)}%
  \pgfmathsetmacro#2{cos(#3)}%
}
 
\newcommand\LongitudePlane[3][current plane]{%
  \pgfmathsinandcos\sinEl\cosEl{#2} 
  \pgfmathsinandcos\sint\cost{#3} 
  \tikzset{#1/.estyle={cm={\cost,\sint*\sinEl,0,\cosEl,(0,0)}}}
}
 
\newcommand\LatitudePlane[3][current plane]{%
  \pgfmathsinandcos\sinEl\cosEl{#2} 
  \pgfmathsinandcos\sint\cost{#3} 
  \pgfmathsetmacro\yshift{\cosEl*\sint}
  \tikzset{#1/.estyle={cm={\cost,0,0,\cost*\sinEl,(0,\yshift)}}} %
}
 
\newcommand\DrawLongitudeCircle[2][1]{
  \LongitudePlane{\angEl}{#2}
  \tikzset{current plane/.prefix style={scale=#1}}
  \pgfmathsetmacro\angVis{atan(sin(#2)*cos(\angEl)/sin(\angEl))} %
  \draw[current plane,thin,black] (\angVis:1) arc (\angVis:\angVis+180:1);
  \draw[current plane,thin,dashed] (\angVis-180:1) arc (\angVis-180:\angVis:1);
}

\newcommand\DrawLongitudeCirclered[2][1]{
  \LongitudePlane{\angEl}{#2}
  \tikzset{current plane/.prefix style={scale=#1}}
  \pgfmathsetmacro\angVis{atan(sin(#2)*cos(\angEl)/sin(\angEl))} %
  \draw[current plane,red,thick] (150:1) arc (150:180:1);
}

\newcommand\DLongredd[2][1]{
  \LongitudePlane{\angEl}{#2}
  \tikzset{current plane/.prefix style={scale=#1}}
  \pgfmathsetmacro\angVis{atan(sin(#2)*cos(\angEl)/sin(\angEl))} %
  \draw[current plane,black,dashed, ultra thick] (150:1) arc (150:180:1);
}

\newcommand\DLatred[2][1]{
  \LatitudePlane{\angEl}{#2}
  \tikzset{current plane/.prefix style={scale=#1}}
  \pgfmathsetmacro\sinVis{sin(#2)/cos(#2)*sin(\angEl)/cos(\angEl)}
  \pgfmathsetmacro\angVis{asin(min(1,max(\sinVis,-1)))}
  \draw[current plane,dashed,black,ultra thick] (-50:1) arc (-50:-35:1);
}

\newcommand\fillred[2][1]{
  \LongitudePlane{\angEl}{#2}
  \tikzset{current plane/.prefix style={scale=#1}}
  \pgfmathsetmacro\angVis{atan(sin(#2)*cos(\angEl)/sin(\angEl))} %
  \draw[current plane,red,thin] (\angVis:1) arc (\angVis:\angVis+180:1);
}
 
\newcommand\DrawLatitudeCircle[2][1]{
  \LatitudePlane{\angEl}{#2}
  \tikzset{current plane/.prefix style={scale=#1}}
  \pgfmathsetmacro\sinVis{sin(#2)/cos(#2)*sin(\angEl)/cos(\angEl)}
  \pgfmathsetmacro\angVis{asin(min(1,max(\sinVis,-1)))}
  \draw[current plane,thin,black] (\angVis:1) arc (\angVis:-\angVis-180:1);
  \draw[current plane,thin,dashed] (180-\angVis:1) arc (180-\angVis:\angVis:1);
}

\newcommand\DrawLatitudeCircleThick[2][1]{
  \LatitudePlane{\angEl}{#2}
  \tikzset{current plane/.prefix style={scale=#1}}
  \pgfmathsetmacro\sinVis{sin(#2)/cos(#2)*sin(\angEl)/cos(\angEl)}
  \pgfmathsetmacro\angVis{asin(min(1,max(\sinVis,-1)))}
  \draw[current plane, very thick,dashed] (\angVis:1) arc (\angVis:-\angVis-180:1);
  \draw[current plane, very thick,black] (180-\angVis:1) arc (180-\angVis:\angVis:1);
}

\newcommand\DrawLatitudeCircleRed[2][1]{
  \LatitudePlane{\angEl}{#2}
  \tikzset{current plane/.prefix style={scale=#1}}
  \pgfmathsetmacro\sinVis{sin(#2)/cos(#2)*sin(\angEl)/cos(\angEl)}
  \pgfmathsetmacro\angVis{asin(min(1,max(\sinVis,-1)))}
  \draw[current plane,very thick,red] (\angVis:1) arc (\angVis:-\angVis-180:1);
  \draw[current plane,very thick,red, dashed] (180-\angVis:1) arc (180-\angVis:\angVis:1);
} 

\newcommand\DrawLatitudeCircleRedoncont[2][1]{
  \LatitudePlane{\angEl}{#2}
  \tikzset{current plane/.prefix style={scale=#1}}
  \pgfmathsetmacro\sinVis{sin(#2)/cos(#2)*sin(\angEl)/cos(\angEl)}
  \pgfmathsetmacro\angVis{asin(min(1,max(\sinVis,-1)))}
  \draw[current plane,very thick,red, dashed] (\angVis:1) arc (\angVis:-\angVis-180:1);
  \draw[current plane,very thick,red] (180-\angVis:1) arc (180-\angVis:\angVis:1);
}

\newcommand\DrawThetaLine[2][1]{
  \LongitudePlane{\angEl}{#2}
  \tikzset{current plane/.prefix style={scale=#1}}
  \pgfmathsetmacro\angVis{atan(sin(#2)*cos(\angEl)/sin(\angEl))} %
  \draw[current plane,thick,black] (\angVis:1) arc (\angVis:\angVis+180:1);
}

\newcommand\DrawLatitudeCircleThickNoDash[2][1]{
  \LatitudePlane{\angEl}{#2}
  \tikzset{current plane/.prefix style={scale=#1}}
  \pgfmathsetmacro\sinVis{sin(#2)/cos(#2)*sin(\angEl)/cos(\angEl)}
  \pgfmathsetmacro\angVis{asin(min(1,max(\sinVis,-1)))}
  \draw[current plane, very thick,black] (180-\angVis:1) arc (180-\angVis:\angVis:1);
}

\newcommand\DrawLatitudeCircleRedNoDash[2][1]{
  \LatitudePlane{\angEl}{#2}
  \tikzset{current plane/.prefix style={scale=#1}}
  \pgfmathsetmacro\sinVis{sin(#2)/cos(#2)*sin(\angEl)/cos(\angEl)}
  \pgfmathsetmacro\angVis{asin(min(1,max(\sinVis,-1)))}
  \draw[current plane,very thick,red] (\angVis:1) arc (\angVis:-\angVis-180:1);
} 

\newcommand\DrawLatitudeCircleRedoncontNoDash[2][1]{
  \LatitudePlane{\angEl}{#2}
  \tikzset{current plane/.prefix style={scale=#1}}
  \pgfmathsetmacro\sinVis{sin(#2)/cos(#2)*sin(\angEl)/cos(\angEl)}
  \pgfmathsetmacro\angVis{asin(min(1,max(\sinVis,-1)))}
  \draw[current plane,very thick,red] (180-\angVis:1) arc (180-\angVis:\angVis:1);
}


\clearpage
\setcounter{tocdepth}{2}
\tableofcontents
\clearpage

\section{Introduction}
\label{section: introduction}

\subsection{Energy Balance Models}
Energy balance models (EBMs) are among the simplest climate models. These models were almost simultaneously popularized by the appearance of two pioneering papers by Byduko \cite{budyko} and Sellers \cite{sellers1969global} in the late 60s, and have continued to serve as valuable tools in understanding the Earth's climate. Despite recent the advancements in the computational ability and the appearance of ever more sophisticated ESMs, conceptual models have not disappeared. In fact, the complexity of realistic models has highlighted the need for a hierarchical model structure where conceptual models provide a solid theoretical foundation as model complexity increase \cite{claussen2002earth}. The simplicity of EBMs allows for studies using standard numerical methods, as well as analytical methods in some cases. Consequently, a unique insight into the dynamics of these models is obtainable, a notable advantage over more complex models. In particular, EBMs have proved valuable in understanding the history of the Earth's climate. EMBs have a rich multiple solution structure \cite{north1990multiple}, which has demonstrated the multistability of the Earth's climate due the destabilizing effect of the ice–albedo feedback. In essence, this feedback mechanism operates on the basis that ice effectively reflects incoming solar radiation, resulting in a decrease in surface temperature. This, in turn, triggers further expansion of the ice cover on the surface. The bistable nature of the climate system was later confirmed by geological and paleomagnetic evidence \cite{hoffman1998neoproterozoic} and the same dynamic has been found in fully-coupled climate models \cite{voigt2010transition}.

\subsubsection{Zero-dimensional Models}
The so-called zero-dimensional EBMs are the simplest of all EBMs. They include no spatial variable and instead describe the Earth's global average temperature by relating the incoming solar radiation, $E_{\tx{in}}$, to the outgoing infrared radiation of the Earth, $E_{\tx{out}}$. The time rate of global average surface temperature change is described in terms of the global radiative imbalance,
\begin{equation*}
    C\frac{dT}{dt} = E_{\tx{in}} - E_{\tx{out}}.
\end{equation*}
Here $C$ is the average heat capacity of the Earth. The energy the Earth receives from the sun in the form of electromagnetic radiation may be modelled as
\begin{equation*}
    E_{\tx{in}} = Q(1-a).
\end{equation*}
Here $Q$ denotes the total solar irradiance, also called solar constant,  divided by 4, as the disk silhouette capturing radiation is $4$ times less than the total area of the Earth. The albedo, denoted $a$, is included in the model. The albedo is a measure of the diffuse reflection of solar radiation and varies with surface and atmospheric conditions. Zero-dimensional EBMs may include an average planetary albedo ($a_p \approx 0.3$), but more interesting is to include an elementary form of the ice-albedo feedback mechanism by allowing for a temperature dependent albedo,
\begin{equation*}
    a = a(T).
\end{equation*}
In polar regions, the presence of snow cover and a more concentrated cloud cover typically results in a reduction of albedo. Utilizing a temperature-dependent albedo offers a acceptable approximation of this dynamic on a global scale and stays within the framework of such a simplified model. The Earth emits longwave radiation into outer space at the top of the atmosphere, a phenomenon that can be approximated by Budyko's empirical formula \cite{budyko} under conditions akin to the current climate state of the Earth.  By assuming a constant lapse rate, there will be a linear relation between the surface temperature and the top-of-the-atmosphere temperature, and consequently the outgoing energy. This outgoing energy is modelled as
\begin{equation*}
    E_{\tx{out}} = A + BT,
\end{equation*}
where $A$ and $B$ are constants. This form for the outgoing radiation ensures that the model exhibits an important stabilizing mechanism: Increased (surface) temperature ensures an increased emission of energy into outer space. With the prescribed outgoing and incoming energy (and parameters in Appendix \ref{section: Parameter values}) it is possible to obtain a reasonable approximation of the Earth's global average temperature.

\subsubsection{Budyko-type Models}
Building upon the notion of a global EBM, the pursuit of an improved temperature estimate naturally leads to the inclusion an additional degree of freedom to the model, extending the concept an energy budget to encompass an infinitesimal segment of the Earth's surface. By incorporating a spatial variable into the model, the objective is to encompass the most important global disparities, making a latitude dependence a reasonable first axis of variability. Invoking a spherical coordinate system where the polar angle is latitude, designated by $\theta$, the latitude ranges from $\theta = 0$, the North pole, to $\theta=\pi$, the South pole. Assuming rotational symmetry, averaging over long times scales the diurnal cycle may be neglected and the model varies only in time and with latitude. Consequently, the temperature is constant along zonal strips around the sphere for a given time $t$,
\begin{equation*}
    T=T(\theta, t).
\end{equation*}
Therefore, the model describes the zonal mean surface temperature, or in other words the surface temperature averaged over latitude. Satellite observations indicate a considerable disparity in the net absorption of solar radiation between tropical latitudes and polar regions \cite{north2017energy2}. To account for this, the model includes an average annual latitudinal radiation distribution function, $s(\theta)$. The current obliquity of the Earth and its resulting seasonal effects are built into $s(\theta)$. The averaging effect of $s(\theta)$ naturally constrains the valid variability of the model to within this timescale, restricting it to describing the annual zonal mean surface temperature. To sustain a steady state, it becomes imperative to invoke a heat transport mechanism facilitating the transfer of heat to colder regions. Budyko \cite{budyko} and Sellers \cite{sellers1969global} independently proposed a model where meridonal heat transport in the atmosphere and hydrosphere is described through a linear term on the form $\beta(T - T_0)$, where $T = T(\theta)$ is the temperature at a given latitude, and $T_0$ is the planetary mean temperature. These models have subsequently been dubbed, \textit{Budyko-Seller-type} EBMs, or simply Budyko-type models, and may be expressed as 
\begin{equation*}
    C \parD{T}{t}  - \beta(T - T_0) = Qs(\theta)(1-a(T)) - A + BT.
\end{equation*}
The Budyko-Sellers-type model incorporates the ice-albedo feedback mechanism, as well as the radiative feedback response, and exhibits interesting ice line dynamics \cite{dynamics2011widiasih} and multistability \cite{north1990multiple}. The simplicity of the governing equation ensure that the system is easily studied using analytical methods \cite{held1974simple}. The addition of a transport term ensures that heat energy fluxes are carried across latitude circles from warm to cool areas. However, the geophysical fluids transport heat through both a mean flow and numerous uncorrelated eddy processes \cite{north1981energy}. Thus, the empirical linear transport term may not be representative for conditions very far from those of the current climate, a central query for EBMs. Below we will discuss a meridional transport term that more closely resembles the physical process of heat dispersion, namely a diffusive heat transport.

\subsubsection{North-type Models}
\label{section: Model}
One-dimensional EBMs with a diffusive heat transport term, hereafter called \textit{North-type models}, have been highly influential and are well-studied. Gerald North has conducted pioneering research, with a focus on analytical investigations, into these models \cite{north1975analytical, north1975theory, north1981energy}. In these EBMs, various physical prosesses are lumped into a the thermal diffusion term $-D\nabla^2 T$, ensuring that thermal energy is being redistributed across latitudes with a flux density proportional to the negative local temperature gradient and scaled by the thermal diffusion coefficient, $D$. With the appropriate expression of the Laplace operator in spherical coordinates the time-dependent energy balance equation takes the form
\begin{equation}
     C \parD{T}{t}  -D\frac{1}{\sin\theta}\parD{}{\theta }\big( \sin\theta \parD{T}{\theta} \big)  + BT= Qs(\theta)(1-a) - A.
    \label{eqn: EBM on symmetric water sphere}
\end{equation}
In this paper, we study this EBM. North showed that by subjecting the model to some basic physical constraints, it is possible to obtain an analytical solution using Fourier-Legendre series for a step-function albedo \cite{north1975analytical}. These constrains are retained in what follows. A restatement of these constraints is: Firstly, we will allow for ice/snow cover on the surface if the local temperature is sufficiently low. We model this threshold as
\begin{equation*}
    T_s(\theta, t) = -T_s,
\end{equation*}
where $-T_s$ refers to the critical temperature for the ice-water phase transition. Consequently, there will be a \textit{critical latitude}, $\theta_{c_i}$, at which the ice-cover ends and begins. The temperature at the critical latitude must necessarily be
\begin{equation*}
    T(\theta_{c_i}) = -T_s.
\end{equation*}
Secondly, the solution must be continuous and continuously differentiable across the critical latitude. Finally, the gradient must vanish at the boundaries as we allow for no heat transport at the poles, leading to the boundary conditions
\begin{align}
    \lim_{\theta \rightarrow 0} \sin \theta \parD{T}{\theta}(\theta)  &=0
    \label{eqn: BC theta=0}
    \intertext{and}
    \lim_{\theta \rightarrow \pi} \sin \theta \parD{T}{\theta}(\theta) &=0.
    \label{eqn: BC theta=pi}
\end{align}
 
When solving (\ref{eqn: EBM on symmetric water sphere}), a major challenge arises in determining the location of the critical latitude under the given constrains. We show that the \textit{Boundary Integral Method} provides a framework for addressing this problem, even for states with multiple critical latitudes. Before proceeding with this analysis, a non-dimensional form of (\ref{eqn: EBM on symmetric water sphere}) will be derived.

\vsone

\subsubsection{Scaling the Model}
\label{section: Nondimensionalize the model}
Following standard non-dimensionalization techniques, (\ref{eqn: EBM on symmetric water sphere}) is  non-dimensionalized using $T_s$ as a scale for temperature, $T=T_sT'$. A timescale, $t_0$, is included in the model, $t = t_0 t'$. Here prime denotes non-dimensional variables. With the prescribed scale, the time-dependent energy balance equation is written on dimensionless form;
\begin{equation}
     \gamma \partial_t T -\frac{1}{\sin\theta}\parD{}{\theta }\big( \sin\theta \parD{T}{\theta} \big) + \beta T= \eta s(\theta)(1-a(T)) - \alpha,
     \label{eqn: non-dimensional time dependent EBM}
\end{equation}
where
\begin{equation}
    \begin{aligned}
        \gamma &= \frac{C}{t_0 D},\\
    \beta &= \frac{B}{D},\\
    \alpha &= \frac{A}{T_s D},\\
    \eta &= \frac{Q}{T_s D}.
    \end{aligned}
    \label{eqn: dimensionless parameters}
\end{equation}
The numerical values for the parameters used in the presented work can be found in Appendix \ref{section: Parameter values}. In equation (\ref{eqn: non-dimensional time dependent EBM}) the primes denoting non-dimensional temperature and time have been omitted as the subsequent focus will be exclusively directed towards the non-dimensional form of the energy balance equation.

\clearpage
\section{Boundary Integral Method}
\label{section: Boundary Integral Method}
In this section, the Boundary Integral Method (BIM) is applied to solve the stationary energy balance equation,
\begin{equation}
    -\frac{1}{\sin\theta}\parD{}{\theta} \big( \sin\theta \parD{T}{\theta}  \big)     + \beta T= \eta s(\theta)(1-a(T)) - \alpha.
    \label{eqn: dimensionless Energy balance equation on symmetric water sphere}
\end{equation}
In the following analysis, a step-function albedo is applied,
\begin{equation}
    a(T) = \begin{cases} 
    a_1 , & T>-1\\
    a_2 , & T< -1.
    \end{cases}
    \label{eqn: step-function albedo}
\end{equation}
Note that this albedo is a function of non-dimensional temperature. The investigation starts with an idealized aquaplanet, establishing a framework for deriving stationary solutions based on the state of the system. Following this, a rotationally symmetric continent is introduced, and the appearance of new stationary states is discussed. Initially, the continent is put in a configuration preserving meridional symmetry. Subsequently, the continent is gradually relocated disrupting the meridional symmetry and stationary solutions associated with these modified continent configurations are discussed. This systematic reduction of meridional symmetry is intended to reveal patterns relating key system dynamics, such as the emergence of multiple equilibria, to symmetry in subsequent analyses.



\vsone

\subsection{About the Method}
Before proceeding to solve (\ref{eqn: dimensionless Energy balance equation on symmetric water sphere}) using the BIM, a recap of the general procedure of the method is presented. The aim of the BIM is to describe the solution to a differential equation in the interior of a given domain, through a relation including certain boundary points of the domain. Consider a general equation on the form
\begin{equation}
    \mathcal{G}y(\bm{x}) = h(\bm{x}), \hsone \bm{x} \: \in \: D,
    \label{eqn: general DE BIE}
\end{equation}
where $\mathcal{G}$ is some differential operator. Solving (\ref{eqn: general DE BIE}) using the BIM necessitates a \textit{Green's function} for the operator $\mathcal{G}$. A Green's function for the operator $\mathcal{G}$, is any solution to the equation
\begin{equation}
    \mathcal{G}K(\bm{x}, \bm{\xi}) = \delta_{\bm{\xi}}(\bm{x}),
    \label{eqn: general greens function}
\end{equation}
where $K$ is the Green's function and $\delta_{\bm{\xi}}$ is an expression of the the Dirac-delta function on the domain $D$. The Green's function can be applied, in conjunction with an integral identity associated with the operator $\mathcal{G}$, to derive integral equations for the solution to (\ref{eqn: general DE BIE}). These \textit{boundary integral equations} relate values of the solution in the interior of the domain, $D$, to values on the boundary of the domain, $\partial D$. To illustrate this, consider the following example: Suppose an integral identity associated with the operator $\mathcal{G}$ has the form
\begin{equation}
    \int_D dV \: \{u\mathcal{G} v - v\mathcal{G}u  \} = \int_{\partial D} dS \: \{u\mathcal{H}v - v\mathcal{H}u\},
    \label{eqn: general fund integral identity}
\end{equation}
where $u$ and $v$ are some functions defined on $D$ and $\mathcal{H}$ is some operator that ensures that (\ref{eqn: general fund integral identity}) holds. Let $K$ be a Green's function that solves (\ref{eqn: general greens function}). Let $u=y$ and $v=K$ in (\ref{eqn: general fund integral identity}) such that
\begin{align}
    \int_D dV \: \{y\mathcal{G} K - K\mathcal{G}y  \} &= \int_{\partial D} dS \: \{y\mathcal{H}K - K\mathcal{H}y\}
    \nonumber\\
    \int_D dV \: \{y\delta_{\bm{\xi}} - Kh  \} &= \int_{\partial D} dS \: \{y\mathcal{H}K - K\mathcal{H}y\}
    \nonumber\\
    y(\bm{\xi}) &= \int_D dV \: Kh + \int_{\partial D} dS \: \{y\mathcal{H}K - K\mathcal{H}y\} .
    \label{eqn: general integral identity for solution y}
\end{align}
The boundary integral equation (\ref{eqn: general integral identity for solution y}) relates values of the solution to (\ref{eqn: general DE BIE}) inside the domain $D$ to values on the boundary $\partial D$. The unknown boundary values in (\ref{eqn: general integral identity for solution y}) are typically addressed by employing a combination of techniques, including the elimination of unknowns through a specific choice of Green's function, invoking boundary conditions, and the derivation of additional boundary integral equations, resulting in a system of equations involving these unknown boundary values.

\vsone

\subsection{Aquaplanet}
\label{section: BIM for a aquaplanet}
Returning to the stationary energy balance equation (\ref{eqn: dimensionless Energy balance equation on symmetric water sphere}), we initially consider an idealized aquaplanet. The planet's surface is considered as one uniform medium and model parameters take on globally-averaged values. We will solve (\ref{eqn: dimensionless Energy balance equation on symmetric water sphere}) for an idealized aquaplanet using the BIM. As established, the BIM aims to describe the solution to a differential equation within a  domain though boundary integral equations (BIEs). Therefore, we seek a set of BIEs for the solution to (\ref{eqn: dimensionless Energy balance equation on symmetric water sphere}). The differential operator of interest is
\begin{equation}
    \mathcal{L}(\cdot) = -\frac{1}{\sin\theta}\parD{}{\theta} \big( \sin\theta \parD{}{\theta} (\cdot) \big)     + \beta (\cdot).
    \label{eqn: operator with beta not appendix}
\end{equation}
In Appendix \ref{Appendix: Fundamental integral identity and a Green's function} an integral identity associated with the operator (\ref{eqn: operator with beta not appendix}) is derived. For two functions $v$ and $u$ on the longitude line $[\theta_1, \theta_2]$ on the surface of a sphere, we have
\begin{equation}
    \int_{\theta_1}^{\theta_2} d\theta\: \sin \theta \big\{v\mathcal{L}u - u \mathcal{L}v \big\} = \left\{ u \sin\theta \parD{v}{\theta} -v \sin\theta \parD{u}{\theta}\right\}\bigg|_{\theta_1}^{\theta_2}.
    \label{eqn: integral identity not appendix}
\end{equation}
A Green's function for the operator (\ref{eqn: operator with beta not appendix}) is required. That is, any solution to the equation
\begin{equation}
    \mathcal{L}K(\theta, \xi) = \delta_\xi (\theta),
    \label{eqn: equation for greens function not appendix}
\end{equation}
where $\delta_\xi (\theta)$ is the Dirac-delta function. In Appendix \ref{Appendix: Fundamental integral identity and a Green's function} a suitable Green's function is found,
\begin{equation}
           K(\theta, \xi) = \begin{cases}
    \frac{P_\lambda (\cos \xi)(\pi \cot (\pi \lambda)P_\lambda(\cos \theta) - 2 Q_\lambda(\cos \theta) )}{ 2(1+\lambda)(P_\lambda(\cos\xi)Q_{\lambda+1}(\cos\xi) - P_{\lambda+1}(\cos\xi)Q_\lambda(\cos\xi))}        , & \theta > \xi\\[10pt]
    \frac{P_\lambda (\cos \theta)(\pi \cot (\pi \lambda)P_\lambda(\cos \xi) - 2 Q_\lambda(\cos \xi) )}{ 2(1+\lambda)(P_\lambda(\cos\xi)Q_{\lambda+1}(\cos\xi) - P_{\lambda+1}(\cos\xi)Q_\lambda(\cos\xi))}        , & \theta < \xi
    \end{cases}.
    \label{eqn: Green's function}
\end{equation} Here $P_\lambda$ and $Q_\lambda$ are Legendre functions of order $\lambda$, where
\begin{equation*}
    \lambda = \frac{1}{2}\left(\sqrt{1-4\beta}-1 \right).
\end{equation*}
This Green's function is continuous and bounded on the domain $\theta \in [0, \pi]$, and its derivative is also bounded at the boundaries $\theta = 0$ and $\theta = \pi$, for a given $\xi \in [0, \pi]$.

\vsone

\subsubsection{No Partial Ice Cover}
\label{section: The case of no ice edge, aquaplanet}
Let's start by obtaining solutions when the surface has no partial ice cover. That is the linear problem where the ice albedo feedback is inactive due to extreme temperatures. Defining
\begin{equation*}
    h(T, \theta) = \eta s(\theta)(1-a(T))-\alpha,
\end{equation*}
we may write the stationary energy balance equation (\ref{eqn: dimensionless Energy balance equation on symmetric water sphere}) on the compact form
\begin{equation}
    \mathcal{L}T = h.
    \label{eqn: compact EBM}
\end{equation}
We let $v=T$ and $u=K$ in the integral identity (\ref{eqn: integral identity not appendix}),
\begin{equation*}
            \int_{\theta_1}^{\theta_2} d\theta\: \sin \theta \left\{T\mathcal{L}K - K \mathcal{L}T \right\} = \left\{ K \sin\theta \parD{T}{\theta} - T \sin\theta \parD{K}{\theta}\right\}\bigg|_{\theta_1}^{\theta_2} ,
\end{equation*}
and apply (\ref{eqn: compact EBM}) and (\ref{eqn: equation for greens function not appendix}) to express the solution to (\ref{eqn: compact EBM}) through the relation
\begin{equation}
    T(\xi) = \int_{\theta_1}^{\theta_2} d\theta\: \sin \theta \: K(\theta, \xi)h(T, \theta) + \left\{ K \sin\theta \parD{T}{\theta} - T \sin\theta \parD{K}{\theta}\right\}\bigg|_{\theta_1}^{\theta_2}.
    \label{eqn: general solution on integral form to stationary EBM }
\end{equation}
The step-function albedo (\ref{eqn: step-function albedo}) evidently leads to 
\begin{equation*}
    h(T, \theta) = \begin{cases} 
    \eta s(\theta)(1-a_1) - \alpha, & T>-1\\
    \eta s(\theta)(1-a_2) - \alpha, & T< -1.
    \end{cases}
\end{equation*}
For a surface either entirely covered by ice or entirely devoid of ice, the source term, $h$, will be identical throughout the domain. Consequently, we may apply (\ref{eqn: general solution on integral form to stationary EBM }) on the entire domain, letting $\theta_1 \rightarrow 0^+$ and $\theta_2 \rightarrow \pi^-$,
\begin{equation}
    \begin{aligned}
        T(\xi) &= \int_0^\pi d\theta \: \sin\theta \:K(\theta, \xi) h(T, \theta) + \lim_{\theta_2 \rightarrow \pi^-} K(\theta_2, \xi) \sin \theta_2 \parD{T}{\theta}(\theta_2)\\
                &- \lim_{\theta_2 \rightarrow \pi^-} T(\theta_2) \sin \theta_2 \parD{K}{\theta}(\theta_2, \xi)
                - \lim_{\theta_1 \rightarrow 0^+} K(\theta_1, \xi) \sin \theta_1 \parD{T}{\theta}(\theta_1)\\
                &+ \lim_{\theta_1 \rightarrow 0^+} T(\theta_1) \sin \theta_1 \parD{K}{\theta}(\theta_1, \xi).
    \end{aligned}
    \label{eqn: integral identity applied for no ice edges}
\end{equation}
The fact that the Green's function, $K$, and its derivative is bounded at the boundary, in combination with boundary conditions (\ref{eqn: BC theta=0}) and (\ref{eqn: BC theta=pi}), ensures that equation (\ref{eqn: integral identity applied for no ice edges}) takes on the simpler form
\begin{align}
    T(\xi) &= \int_0^\pi d\theta \: \sin\theta\: K(\theta, \xi) h(T, \theta).
    \label{eqn: BIE with h(T, theta)}
\end{align}
Note that since our Green's function is piecewise, we must necessarily split any relevant integral in order to integrate over the appropriate bounds,
\begin{equation*}
    T(\xi) = \int_0^\xi d\theta \: \sin\theta \:K_-(\theta, \xi) h(T, \theta) + \int_\xi^\pi d\theta\: \sin\theta \:K_+(\theta, \xi) h(T, \theta),
\end{equation*}
where 
\begin{equation}
        K(\theta, \xi) = \begin{cases}
    K_+ (\theta, \xi), & \theta > \xi\\[4pt]
    K_- (\theta, \xi), & \theta < \xi
    \end{cases}.
    \label{eqn: notation of K_+ and K_-}
\end{equation}
We define 
\begin{equation}
    h_1(\theta) = \eta s(\theta)(1-a_1) - \alpha,
    \label{eqn: h1}
\end{equation}
wherever there is no ice cover. Likewise, we define
\begin{equation}
    h_2(\theta) = \eta s(\theta)(1-a_2) - \alpha,
    \label{eqn: h2}
\end{equation}
wherever the is ice cover. The solution to (\ref{eqn: dimensionless Energy balance equation on symmetric water sphere}) for the case where the surface is entirely cover by water is
\begin{align}
    T(\xi) &= \int_0^\pi d\theta \: \sin\theta\: K(\theta, \xi) h_1(\theta), \hsone  \forall \: \xi \in [0, \pi].
    \label{eqn: BIE no ice cover, no ice edge}
    \intertext{Similarly, the solution to (\ref{eqn: dimensionless Energy balance equation on symmetric water sphere}) for the case where the surface is entirely cover by ice or snow is}
    T(\xi) &= \int_0^\pi d\theta \: \sin\theta\: K(\theta, \xi) h_2(\theta), \hsone  \forall \: \xi \in [0, \pi].
    \label{eqn: BIE ice cover, no ice edge}
\end{align}

\vsone

The solutions  (\ref{eqn: BIE no ice cover, no ice edge}) and (\ref{eqn: BIE ice cover, no ice edge}) are readily obtainable from quadrature methods. The integrals cannot be evaluated analytically due to the form of the chosen Green's function. Corresponding temperature profiles are shown in dimensional form in Figure \ref{fig: plot of solution no egdes}. Model parameters are found in Appendix  \ref{section: Parameter values}.

\begin{figure}[H]
\centering
\begin{subfigure}{.5\textwidth}
  \centering
  \includegraphics[width=\linewidth]{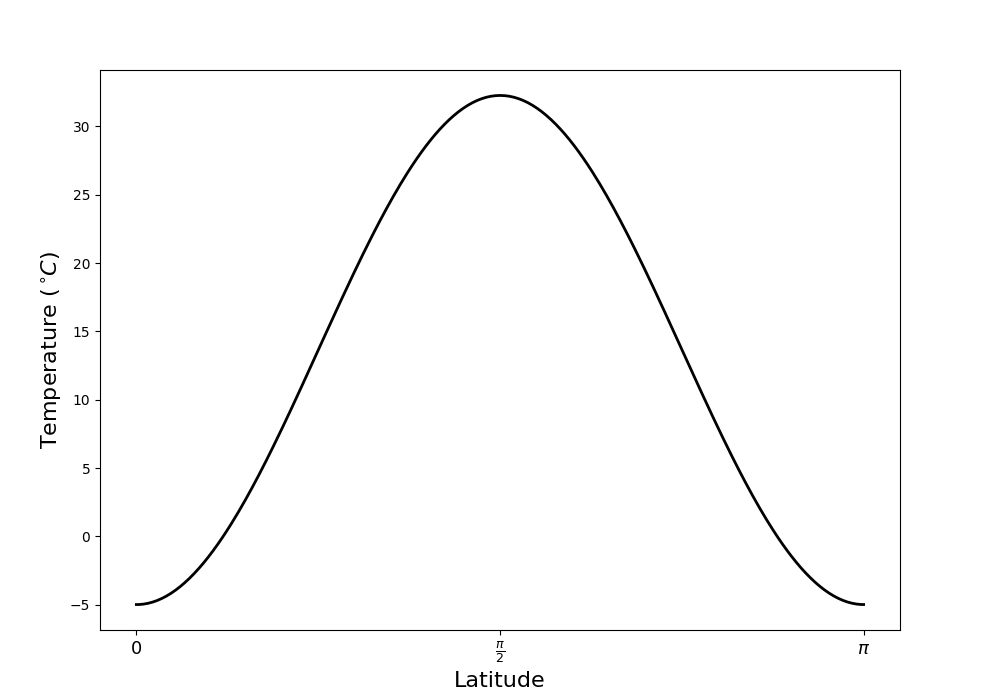}
  \caption{}
  \label{fig: plot of solution no egde, no ice}
\end{subfigure}%
\begin{subfigure}{.5\textwidth}
  \centering
  \includegraphics[width=\linewidth]{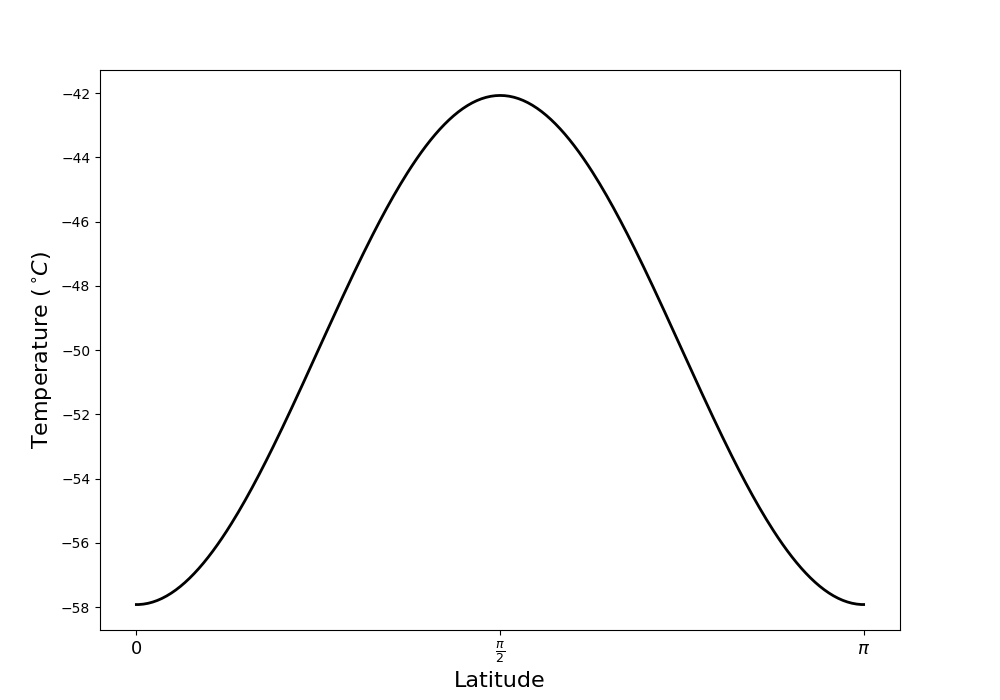}
  \caption{}
  \label{fig: plot of solution no egde, only ice}
\end{subfigure}
\caption{Temperature profile of an aquaplanet with no ice or snow present on the surface (a) and an aquaplanet with its surface entirely covered by ice and/or snow (b). }
\label{fig: plot of solution no egdes}
\end{figure}

\vsone

\subsubsection{Partial Ice Cover}
\label{section: Boundary integral equation for the case of two ice edges}

For solutions to (\ref{eqn: dimensionless Energy balance equation on symmetric water sphere}) where the zonal mean temperature profile is not strictly above or below the critical temperature, $T_s$, the surface will have a partial ice cover analogous to the Earth's current climate state. There will be certain critical latitudes at which ice cover ends and begins, denoted $\theta_{c_1}$ and $\theta_{c_2}$. In solving (\ref{eqn: dimensionless Energy balance equation on symmetric water sphere}), a key concern is to identify these critical latitudes. The BIM provides a convenient way to obtain these. Since the latitudinal energy distribution function, $s(\theta)$, is symmetrical about equator and decreasing towards the poles, the ice cover on the surface will stretch from the poles to some unknown latitudes on the northern and southern hemisphere. Our aim will be to derive a system of BIEs, which can be solved to determine the critical latitudes, $\theta_{c_1}$ and $\theta_{c_2}$.
In order to obtain a temperature profile, the values of $\theta_{c_1}$ and $\theta_{c_2}$ are essential.

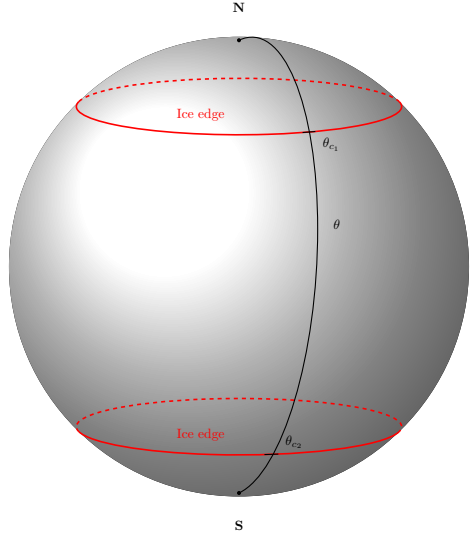
\begin{figure}[H]
\centering
\resizebox{0.4\textwidth}{!}{ 

\begin{tikzpicture}[scale=1,every node/.style={minimum size=1cm}]
    
    \def\R{6} 
    
    \def\angEl{10} 
    \def\angPhiThetaline{-110} 
    \def\angIceEdgeN{45} 
    \def\angIceEdgeS{-45} 
    
    \pgfmathsetmacro\CN{\R*cos(\angEl)*sin(\angIceEdgeN)} 
    \pgfmathsetmacro\CS{\R*cos(\angEl)*sin(\angIceEdgeS)} 
    \pgfmathsetmacro\H{\R*cos(\angEl)} 
    \LongitudePlane[Theta_line_plane]{\angEl}{\angPhiThetaline}
    \LatitudePlane[North_Ice_edge]{\angEl}{\angIceEdgeN}
    \LatitudePlane[South_Ice_edge]{\angEl}{\angIceEdgeS}
    \fill[ball color=white!10] (0,0) circle (\R); 
    \coordinate (CN) at (0,\CN);
    \coordinate (CS) at (0,\CS);
    \coordinate[mark coordinate] (N) at (0,\H);
    \coordinate[mark coordinate] (S) at (0,-\H);

    \DrawLatitudeCircleRed[\R]{\angIceEdgeN}   
    \DrawLatitudeCircleRed[\R]{\angIceEdgeS}
    \node[above=8pt] at (N) {$\mathbf{N}$};
    \node[below=8pt] at (S) {$\mathbf{S}$};
    \node[xshift=-1cm, yshift=-0.2cm] at (CN) {\textcolor{red}{Ice edge}};
    \node[xshift=-1cm, yshift=-0.2cm] at (CS) {\textcolor{red}{Ice edge}};
 
    \path[Theta_line_plane] (\angIceEdgeN:\R) coordinate (IceN);
    \path[Theta_line_plane] (\angIceEdgeS:\R) coordinate (IceS);
    \draw[Theta_line_plane, thick, black] ([shift={(90:\R)}]0,0) arc  (90:90+180:\R) node[midway, above=-13pt, xshift=0.5cm] {$\theta$};
    \draw[North_Ice_edge, thick, black] ([shift={(-67:\R)}]CN) arc  (-67:-62:\R) node[near end, xshift=0.5cm, below=-7pt] {$\theta_{c_1}$};
    \draw[South_Ice_edge, thick, black] ([shift={(-81:\R)}]CS) arc  (-81:-76:\R) node[near end, xshift=0.5cm, above=-7pt] {$\theta_{c_2}$};

\end{tikzpicture}
}
\caption{Schematic of the surface of an aquaplanet with a partial ice cover.} 
\label{fig: Schematic of the surface of the water planet with a partial ice cover}
\end{figure}

We divide our domain into three sub-domains based on the ability of the surface to absorb incident radiation,
\begin{region}
    \theta \in (0, \theta_{c_1}),
    \label{eqn: region 1}
\end{region}
\begin{region}
    \theta \in (\theta_{c_1}, \theta_{c_2})
    \label{eqn: region 2}
\end{region}
and
\begin{region}
    \theta \in (\theta_{c_2}, \pi).
    \label{eqn: region 3}
\end{region}
As there is no ice cover in region (\ref{eqn: region 2}) we let (\ref{eqn: h1}) be the source term in the governing equation for region (\ref{eqn: region 2}),
\begin{equation*}
    \mathcal{L}T = h_1.
\end{equation*}
In region (\ref{eqn: region 1}) and (\ref{eqn: region 3}), where there is ice cover, we let (\ref{eqn: h2}) be the source term in the governing equation,
\begin{equation*}
    \mathcal{L}T = h_2.
\end{equation*}
The solution within three regions (\ref{eqn: region 1}), (\ref{eqn: region 2}) and (\ref{eqn: region 3}) are obtained by apply the general relation (\ref{eqn: general solution on integral form to stationary EBM }) in the respective regions, letting $\theta_1$ and $\theta_2$ approach the boundaries of the sub-domains. Starting with region (\ref{eqn: region 1}), we apply (\ref{eqn: general solution on integral form to stationary EBM }) and let $\theta_1 \rightarrow 0^+$ and $\theta_2 \rightarrow \theta_{c_1}^-$,

\begin{align}
    T(\xi) &= \int_0^{\theta_{c_1}}d\theta \: \sin\theta \: K(\theta, \xi) h_2(\theta) +
    \lim_{\substack{\theta_1  \rightarrow 0^+\\
        \theta_2 \rightarrow \theta_{c_1}^-}} \left\{ K \sin\theta \parD{T}{\theta} - T \sin\theta \parD{K}{\theta}\right\}\bigg|_{\theta_1}^{\theta_2}\nonumber\\
            &
            \begin{aligned}
                &= \int_0^\xi d\theta \: \sin\theta \:K_-(\theta, \xi) h_2(\theta) + \int_\xi^{\theta_{c_1}} d\theta\: \sin\theta \:K_+(\theta, \xi) h_2(\theta)\\
                &+ \lim_{\theta_2 \rightarrow \theta_{c_1}^-} K(\theta_{2}, \xi) \sin(\theta_{2}) \parD{T}{\theta}(\theta_{2})\\
                &- \lim_{\theta_2 \rightarrow \theta_{c_1}^-} T(\theta_{2}) \sin(\theta_{2}) \parD{K}{\theta}(\theta_{2},\xi)\\
                &-\lim_{\theta_1  \rightarrow 0^+}K(\theta_1, \xi) \sin (\theta_1) \parD{T}{\theta}(\theta_1)\\
                &+\lim_{\theta_1  \rightarrow 0^+}T(\theta_1) \sin (\theta_1) \parD{K}{\theta}(\theta_1, \xi).
            \end{aligned}
            \label{eqn: general integral relation applied region 1}
\end{align}
Assuming that the solution is bounded at $\theta = 0$, the last term in (\ref{eqn: general integral relation applied region 1}) is zero. The Green's function is also bounded at $\theta = 0$ and therefore the second-to-last term in (\ref{eqn: general integral relation applied region 1}) is also zero due to the boundary condition (\ref{eqn: BC theta=0}). Furthermore, we will assume that $T$ is continuous across $\theta = \theta_{c_1}$ and $\theta = \theta_{c_2}$, and the temperature at these critical latitudes must be $-T_s$. The value of the non-dimensional temperature at the latitudes is therefore
\begin{equation*}
    T(\theta_{c_1})=T(\theta_{c_2})=-1.
\end{equation*}
Furthermore, assuming that $\parD{T}{\theta}$ is continuous at the critical latitudes we must have
\begin{align*}
    \lim_{\theta \rightarrow \theta_{c_1}^-}\parD{T}{\theta}(\theta) &= \lim_{\theta \rightarrow \theta_{c_1}^+}\parD{T}{\theta}(\theta)=\parD{T}{\theta}(\theta_{c_1})\\
    \lim_{\theta \rightarrow \theta_{c_2}^-}\parD{T}{\theta}(\theta) &= \lim_{\theta \rightarrow \theta_{c_2}^+}\parD{T}{\theta}(\theta)=\parD{T}{\theta}(\theta_{c_2}).
\end{align*}
Applying this, as well as the fact that the Green's function is continuous at $\theta = \theta_{c_1}$, to (\ref{eqn: general integral relation applied region 1}) we get
\begin{equation}
    \begin{aligned}
    T(\xi)  &= \int_0^\xi d\theta \: \sin\theta \:K_-(\theta, \xi) h_2(\theta) + \int_\xi^{\theta_{c_1}} d\theta\: \sin\theta \:K_+(\theta, \xi) h_2(\theta)\\
            &+  K(\theta_{c_1}, \xi) \sin(\theta_{c_1}) \lim_{\theta_2 \rightarrow \theta_{c_1}}\parD{T}{\theta}(\theta_{2})
            + \sin(\theta_{c_1}) \lim_{\theta_2 \rightarrow \theta_{c_1}^-} \parD{K}{\theta}(\theta_{2},\xi).
    \end{aligned}
    \label{eqn: solution in region 1}
\end{equation}

\vspace{0.7cm}

The relation (\ref{eqn: general solution on integral form to stationary EBM }) is also applied in region (\ref{eqn: region 2}), letting $\theta_1 \rightarrow \theta_{c_1}^+$ and $\theta_2 \rightarrow \theta_{c_2}^-$,
\begin{equation}
    \begin{aligned}
        T(\xi) &= \int_{\theta_{c_1}}^\xi d\theta \: \sin\theta \:K_-(\theta, \xi) h_1(\theta) + \int_\xi^{\theta_{c_2}} d\theta\: \sin\theta \:K_+(\theta, \xi) h_1(\theta)\\
                &+ \lim_{\theta_2 \rightarrow \theta_{c_2}^-} K(\theta_{2}, \xi) \sin(\theta_{2}) \parD{T}{\theta}(\theta_{2})\\
                &- \lim_{\theta_2 \rightarrow \theta_{c_2}^-} T(\theta_{2}) \sin(\theta_{2}) \parD{K}{\theta}(\theta_{2},\xi)\\
                &-\lim_{\theta_1  \rightarrow \theta_{c_1}^+}K(\theta_1, \xi) \sin (\theta_1) \parD{T}{\theta}(\theta_1)\\
                &+\lim_{\theta_1  \rightarrow \theta_{c_1}^+}T(\theta_1) \sin (\theta_1) \parD{K}{\theta}(\theta_1, \xi).
    \end{aligned}
    \label{eqn: general integral relation applied region 2}
\end{equation}
Under the assumption that $T$, $\parD{T}{\theta}$ and $K$ are continuous at $\theta = \theta_{c_1}$ and $\theta = \theta_{c_2}$, (\ref{eqn: general integral relation applied region 2}) becomes
\begin{equation}
    \begin{aligned}
        T(\xi) &= \int_{\theta_{c_1}}^\xi d\theta \: \sin\theta \:K_-(\theta, \xi) h_1(\theta) + \int_\xi^{\theta_{c_2}} d\theta\: \sin\theta \:K_+(\theta, \xi) h_1(\theta)\\
                &+  K(\theta_{c_2}, \xi) \sin(\theta_{c_2}) \lim_{\theta_2 \rightarrow \theta_{c_2}}\parD{T}{\theta}(\theta_{2})\\
                &+ \sin(\theta_{c_2}) \lim_{\theta_2 \rightarrow \theta_{c_2}^-}\parD{K}{\theta}(\theta_{2},\xi)\\
                &-K(\theta_{c_1}, \xi) \sin (\theta_{c_1}) \lim_{\theta_1  \rightarrow \theta_{c_1}}\parD{T}{\theta}(\theta_1)\\
                &- \sin (\theta_{c_1}) \lim_{\theta_1  \rightarrow \theta_{c_1}^+} \parD{K}{\theta}(\theta_1, \xi).
    \end{aligned}
    \label{eqn: solution in region 2}
\end{equation}

\vspace{0.7cm}

Finally, applying (\ref{eqn: general solution on integral form to stationary EBM }) in region (\ref{eqn: region 3}) letting $\theta_1 \rightarrow \theta_{c_2}^+$ and $\theta_2 \rightarrow \pi^-$ we get
\begin{equation}
    \begin{aligned}
        T(\xi) &= \int_{\theta_{c_2}}^\xi d\theta \: \sin\theta \:K_-(\theta, \xi) h_2(\theta) + \int_\xi^{\pi} d\theta\: \sin\theta \:K_+(\theta, \xi) h_2(\theta)\\
                &+ \lim_{\theta_2 \rightarrow \pi^-} K(\theta_{2}, \xi) \sin(\theta_{2}) \parD{T}{\theta}(\theta_{2})\\
                &- \lim_{\theta_2 \rightarrow \pi^-} T(\theta_{2}) \sin(\theta_{2}) \parD{K}{\theta}(\theta_{2},\xi)\\
                &-\lim_{\theta_1  \rightarrow \theta_{c_2}^+}K(\theta_1, \xi) \sin (\theta_1) \parD{T}{\theta}(\theta_1)\\
                &+\lim_{\theta_1  \rightarrow \theta_{c_2}^+}T(\theta_1) \sin (\theta_1) \parD{K}{\theta}(\theta_1, \xi).
    \end{aligned}
    \label{eqn: general integral relation applied region 3}
\end{equation}
Assuming that $T$ bounded at $\theta = \pi$ and that $T$ and $\parD{T}{\theta}$ are continuous at $\theta=\theta_{c_2}$, in combination with the boundary condition (\ref{eqn: BC theta=pi}) and the critical latitude condition, (\ref{eqn: general integral relation applied region 3}) becomes
\begin{equation}
    \begin{aligned}
        T(\xi) &= \int_{\theta_{c_2}}^\xi d\theta \: \sin\theta \:K_-(\theta, \xi) h_2(\theta) + \int_\xi^{\pi} d\theta\: \sin\theta \:K_+(\theta, \xi) h_2(\theta)\\
                &-K(\theta_{c_2}, \xi) \sin (\theta_{c_2}) \lim_{\theta_1  \rightarrow \theta_{c_2}}\parD{T}{\theta}(\theta_1)
                - \sin (\theta_{c_2}) \lim_{\theta_1  \rightarrow \theta_{c_2}}\parD{K}{\theta}(\theta_1, \xi).
    \end{aligned}
    \label{eqn: solution in region 3}
\end{equation}

\vsone

Equation (\ref{eqn: solution in region 1}), (\ref{eqn: solution in region 2}) and (\ref{eqn: solution in region 3}) gives the solution in region (\ref{eqn: region 1}), (\ref{eqn: region 2}) and (\ref{eqn: region 3}), respectively, expressed through an integral and the unknown values $\theta_{c_1}$, $\theta_{c_2}$, $\parD{T}{\theta}(\theta_{c_1})$ and $\parD{T}{\theta}(\theta_{c_2})$. These unknown values may be determined through developing a set of BIEs. The approach will be to let $\xi$ approach the boundaries in each of the sub-domains (\ref{eqn: region 1}), (\ref{eqn: region 2}) and (\ref{eqn: region 3}).

\vsone

Staring with region (\ref{eqn: region 1}), letting $\xi \rightarrow 0^+$ in equation (\ref{eqn: solution in region 1}), we get
\begin{equationTaggedLowerCaseRoman}
\begin{aligned}
    T(0) &=\int_0^{\theta_{c_1}}d\theta \: \sin\theta \: K_+(\theta, 0) h_2(\theta)\\
        &+ K_+(\theta_{c_1}, 0) \sin(\theta_{c_1}) \parD{T}{\theta}(\theta_{c_1})
                + \sin(\theta_{c_1}) \lim_{\xi \rightarrow 0^+}\lim_{\theta \rightarrow \theta_{c_1}^-}\parD{K}{\theta}(\theta, \xi).
                \label{eqn: boundary formulation 1}
\end{aligned}
\end{equationTaggedLowerCaseRoman}
Letting $\xi \rightarrow \theta_{c_1}^-$ in (\ref{eqn: solution in region 1}) we get
\begin{equationTaggedLowerCaseRoman}
    \begin{aligned}
        -1 &=\int_0^{\theta_{c_1}}d\theta \: \sin\theta \: K_-(\theta, \theta_{c_1}) h_2(\theta)\\
        &+ K_+(\theta_{c_1}, \theta_{c_1}) \sin(\theta_{c_1}) \parD{T}{\theta}(\theta_{c_1})
            + \sin(\theta_{c_1}) \lim_{\xi \rightarrow \theta_{c_1}^-}\lim_{\theta \rightarrow \theta_{c_1}^-}\parD{K}{\theta}(\theta,\xi).
    \end{aligned}
    \label{eqn: boundary formulation 2}
\end{equationTaggedLowerCaseRoman}

Similarly, for region (\ref{eqn: region 2}) we let $\xi \rightarrow \theta_{c_1}^+$ in (\ref{eqn: solution in region 2}) and get
\begin{equationTaggedLowerCaseRoman}
        \begin{aligned}
        -1 &=\int_{\theta_{c_1}}^{\theta_{c_2}}d\theta\:\sin\theta \:K_+(\theta, \theta_{c_1}) h_1(\theta)\\
        &+ K_+(\theta_{c_2}, \theta_{c_1})\sin(\theta_{c_2})\parD{T}{\theta}(\theta_{c_2}) + \sin(\theta_{c_2})\lim_{\xi \rightarrow \theta_{c_1}^+}\lim_{\theta \rightarrow \theta_{c_2}^-}\parD{K}{\theta}(\theta,\xi)\\
        &-K_-(\theta_{c_1}, \theta_{c_1})\sin(\theta_{c_1})\parD{T}{\theta}(\theta_{c_1})
        -\sin(\theta_{c_1})\lim_{\xi \rightarrow \theta_{c_1}^+}\lim_{\theta \rightarrow \theta_{c_1}^+}\parD{K}{\theta}(\theta, \xi).
    \end{aligned}
    \label{eqn: boundary formulation 3}
\end{equationTaggedLowerCaseRoman}
Letting $\xi \rightarrow \theta_{c_2}^-$ in (\ref{eqn: solution in region 2}) we get
\begin{equationTaggedLowerCaseRoman}
        \begin{aligned}
        -1 &=\int_{\theta_{c_1}}^{\theta_{c_2}}d\theta\:\sin\theta \:K_-(\theta, \theta_{c_2}) h_1(\theta)\\
        &+ K_+(\theta_{c_2}, \theta_{c_2})\sin(\theta_{c_2})\parD{T}{\theta}(\theta_{c_2}) + \sin(\theta_{c_2})\lim_{\xi \rightarrow \theta_{c_2}^-}\lim_{\theta \rightarrow \theta_{c_2}^-}\parD{K}{\theta}(\theta,\xi)\\
        &-K_-(\theta_{c_1}, \theta_{c_2})\sin(\theta_{c_1})\parD{T}{\theta}(\theta_{c_1})
        -\sin(\theta_{c_1})\lim_{\xi \rightarrow \theta_{c_2}^-}\lim_{\theta \rightarrow \theta_{c_1}^+}\parD{K}{\theta}(\theta, \xi).
    \end{aligned}
    \label{eqn: boundary formulation 4}
\end{equationTaggedLowerCaseRoman}

Finally, for region (\ref{eqn: region 3}) we let $\xi \rightarrow \theta_{c_2}^+$ in (\ref{eqn: solution in region 3}), yielding
\begin{equationTaggedLowerCaseRoman}
\begin{aligned}
        -1 &= \int_{\theta_{c_2}}^\pi d\theta \: \sin\theta \: K_+(\theta, \theta_{c_2}) h_2(\theta) \\
        &-K_-(\theta_{c_2}, \theta_{c_2}) \sin (\theta_{c_2}) \parD{T}{\theta}(\theta_{c_2})
        - \sin (\theta_{c_2})\lim_{\xi \rightarrow \theta_{c_2}^+}\lim_{\theta \rightarrow \theta_{c_2}^+}\parD{K}{\theta}(\theta, \xi).
    \label{eqn: boundary formulation 5}
\end{aligned}
\end{equationTaggedLowerCaseRoman}
Letting $\xi \rightarrow \pi^-$ in (\ref{eqn: solution in region 3}) we get
\begin{equationTaggedLowerCaseRoman}
\begin{aligned}
        T(\pi) &= \int_{\theta_{c_2}}^\pi d\theta \: \sin\theta \: K_-(\theta, \pi) h_2(\theta) \\
        &-K_-(\theta_{c_2}, \pi) \sin (\theta_{c_2}) \parD{T}{\theta}(\theta_{c_2})
        - \sin (\theta_{c_2})\lim_{\xi \rightarrow \pi^-}\lim_{\theta \rightarrow \theta_{c_2}^+}\parD{K}{\theta}(\theta, \xi).
    \label{eqn: boundary formulation 6}
\end{aligned}
\end{equationTaggedLowerCaseRoman}

\begin{figure}[H]
\centering
\begin{subfigure}{.5\textwidth}
  \centering
  \includegraphics[width=\linewidth]{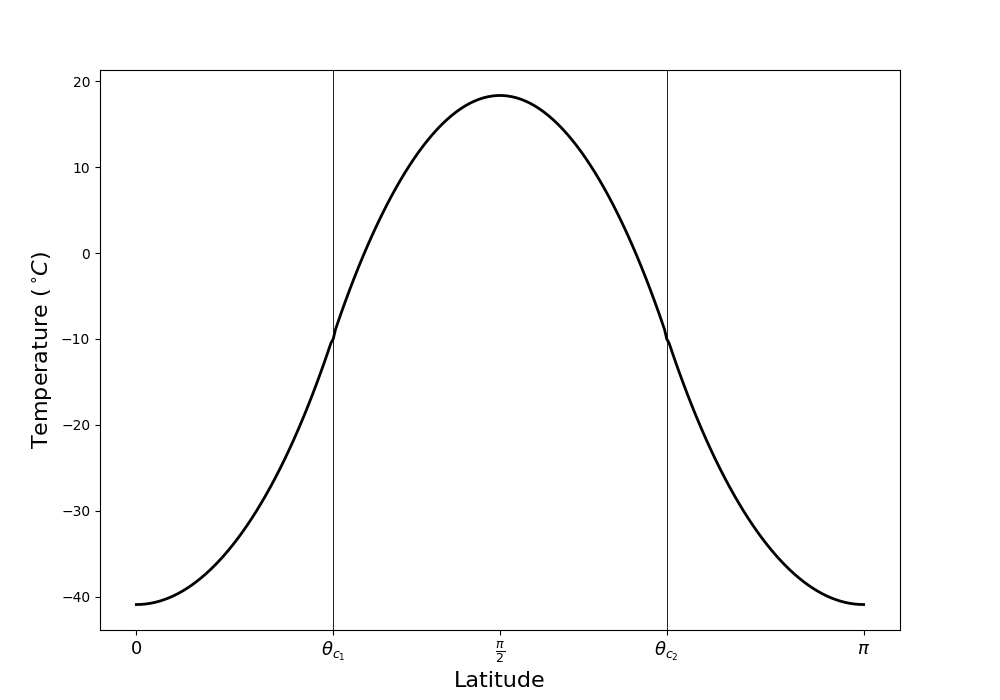}
  \caption{}
  \label{fig: plot of solution w ice egde 1}
\end{subfigure}%
\begin{subfigure}{.5\textwidth}
  \centering
  \includegraphics[width=\linewidth]{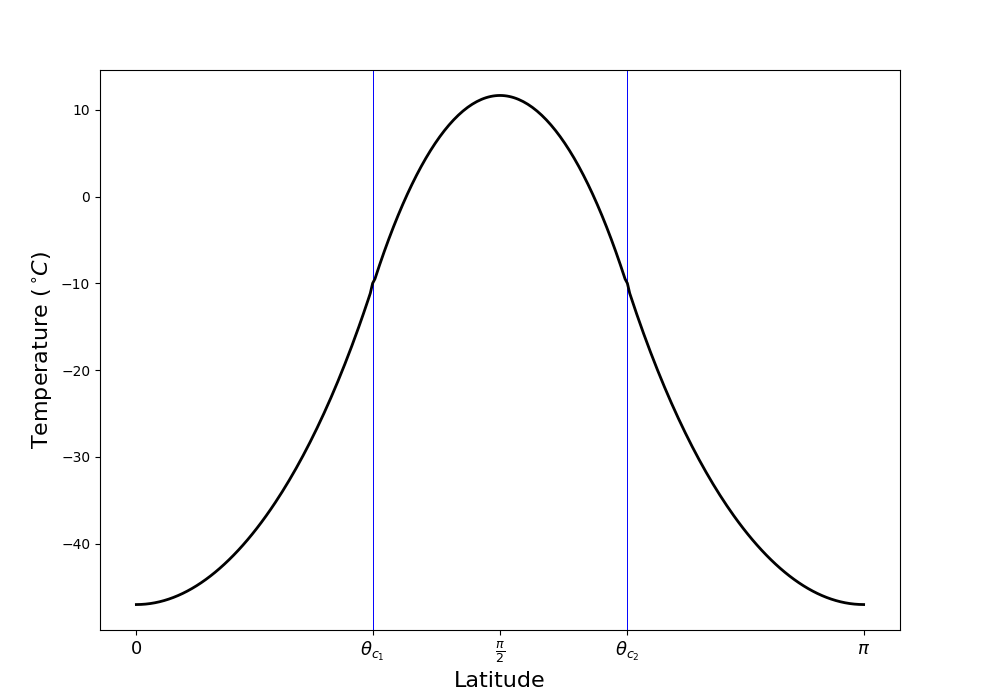}
  \caption{}
  \label{fig: plot of solution w ice egde 2}
  \end{subfigure}\\[1ex]
        \begin{subfigure}{\textwidth} 
            \centering 
            \includegraphics[width=.5\linewidth]{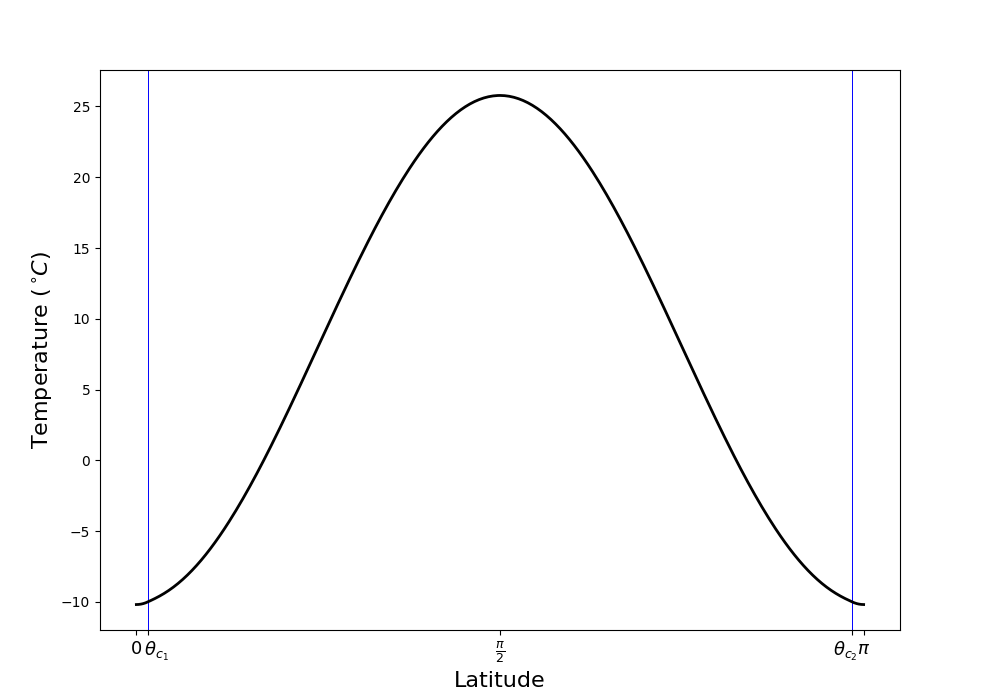}
            \caption{}
            \label{fig: plot of solution w ice egde 3} 
        \end{subfigure} %
\caption{Temperature profile of an aquaplanet with two ice edges. In the relevant parameter regime, two ice edges may be realized in three possible ways, shown in (a), (b) and (c).}
\label{fig: plot of solution w ice egdes}
\end{figure}

The system of equations (\ref{eqn: boundary formulation 1})-(\ref{eqn: boundary formulation 6}) is a system of 6 equations for 6 unknown boundary values, $\theta_{c_1}$, $\theta_{c_2}$, $\parD{T}{\theta}(\theta_{c_1})$, $\parD{T}{\theta}(\theta_{c_2})$, $T(0)$ and $T(\pi)$. Solving this system of equations ultimately require a combination of both analytical and numerical methods as the integrals cannot be evaluated analytically. We start by solving equation (\ref{eqn: boundary formulation 3}) and (\ref{eqn: boundary formulation 4}) for $\parD{T}{\theta}(\theta_{c_1})$ and $\parD{T}{\theta}(\theta_{c_2})$. We substitute the expression for $\parD{T}{\theta}(\theta_{c_1})$ and $\parD{T}{\theta}(\theta_{c_2})$ into equation (\ref{eqn: boundary formulation 2}) and (\ref{eqn: boundary formulation 5}) to create two functions of $\theta_{c_1}$ and $\theta_{c_2}$,
\begin{align*}
    f_1 &= f_1(\theta_{c_1}, \theta_{c_2}),\\
    f_2 &= f_2(\theta_{c_1}, \theta_{c_2}).
\end{align*}
A Newton's iteration is applied to the functions $f_1$ and $f_2$ to obtain $\theta_{c_1}$ and $\theta_{c_2}$. Given the roots $\theta_{c_1}$ and $\theta_{c_2}$, the remaining values, $\parD{T}{\theta}(\theta_{c_1})$, $\parD{T}{\theta}(\theta_{c_2})$, $T(0)$ and $T(\pi)$, are easily found. The solution in the interior of the three sub-domains, (\ref{eqn: region 1}), (\ref{eqn: region 2}) and (\ref{eqn: region 3}), are given by (\ref{eqn: solution in region 1}), (\ref{eqn: solution in region 2}) and (\ref{eqn: solution in region 3}), respectively. The solution at the boundaries of the sub-domains is naturally $T(0)$ and $T(\pi)$ at the North and South pole, respectively, and at the critical latitudes $T(\theta_{c_1})=T(\theta_{c_2})=-1$. 
Using the model parameters in Appendix \ref{section: Parameter values} and $Q=247$, we obtain the solutions shown in dimensional form in Figure \ref{fig: plot of solution w ice egdes}. The blue vertical lines indicate the critical latitudes. For $Q=247$ the system of equations (\ref{eqn: boundary formulation 1})-(\ref{eqn: boundary formulation 6}) has three solutions, meaning that there are three possible ways the planet can realize a temperature distribution consistent with two critical latitudes. Note that the limits of the Green's function's derivative must be handled with meticulous attention, hence the awkward notation.

\subsection{Introducing a Continent}
\label{section: Introducing a continent}
We are going to extend the EBM to include a continent. In order to persist with a one-dimensional analysis, the continent must be introduced in a manner that preserves the rotational symmetry. Ensuring rotational symmetrical demands that a continent is a zonal strip stretching from the latitude $\theta_{l_1}$ to latitude $\theta_{l_2}$. The following analysis will only include one such continent. Parameter values on the continent may be altered to distinguish land from ocean and better capture the thermal response of the lithosphere. Here the heat capacity, $C$, the critical temperature for ice formation, $T_s$ and the albedo, $a(T)$ was altered on the part of the domain corresponding to the continent. This has effect of changing the ice dynamics, and subsequently the ice-albedo feedback, on the continent. Following the non-dimensionalization approach outlined in section \ref{section: Nondimensionalize the model}, yields the following energy balance equation
\begin{equation}
    \begin{cases}
    \gamma_{\text{water}} \partial_t T +\mathcal{L}T + \beta T = \eta s(\theta)(1 - a_{\text{water}}(T)) - \alpha, &\hsone \theta\in [0, \theta_{l_1}) \cup (\theta_{l_2}, \pi]\\[10pt]
    \gamma_{\text{land}} \partial_t T +\mathcal{L}T + \beta T = \eta s(\theta)(1 - a_{\text{land}}(T)) - \alpha, &\hsone  \theta\in [\theta_{l_1}, \theta_{l_2}],
    \end{cases}
    \label{eqn: time dependent EBM w continent}
\end{equation}
where the non-dimensional parameters are given in (\ref{eqn: dimensionless parameters}) and 
\begin{equation}
    \begin{aligned}
        \gamma_{\text{water}}&= \frac{C_{\tx{water}}}{t_0 D},\\
        \gamma_{\text{land}}&= \frac{C_{\tx{land}}}{t_0 D}.
    \end{aligned}
    \label{eqn: gamma}
\end{equation}
Like above, a step function albedo is applied on the water and on the continent. However, the critical temperature for ice formation is set to $-T_{s,\tx{land}}$ within the continent. Consequently, the non-dimensional temperature at critical latitudes within the continent must be
\begin{equation}
    T(\theta_{c_i})=- \frac{T_{s_2}}{T_s} = - T_c.
    \label{eqn: critical latitude condition on continent}
\end{equation}
Equation (\ref{eqn: critical latitude condition on continent}) is the defining condition for a critical latitude on the continent. Owing to a difference in the critical temperature for ice formation between water and land, an ice cover may initiate at the edges of the continent and extend to a critical latitude within the continent's interior. However, the continental edge is not considered a critical latitude unless (\ref{eqn: critical latitude condition on continent}) is satisfied. Consequently, characterizing a critical latitude as an "ice line" or "ice edge" is misleading in this context. The rationale behind the definition (\ref{eqn: critical latitude condition on continent}) is that only critical latitudes need to be identified through the BIM.

\vsone

We are going to apply the BIM to solve the stationary form of 
(\ref{eqn: time dependent EBM w continent}). This will largely follow the procedure outlined in section \ref{section: BIM for a aquaplanet}, with the added complication of a continent. The domain must be partitioned into subdomains as in section \ref{section: Boundary integral equation for the case of two ice edges}. As the continent and the water will have diverging parameter values, the part of the domain corresponding to water and land must be separated into subdomains. This ensure that there will generally be more boundary integral equations, but the procedure for obtaining a stationary solution will be the same as in section \ref{section: Boundary integral equation for the case of two ice edges}. Mathematical details are left out here and the reader is referred to Appendix \ref{appendix: BIEs with a continent} for a thorough description of how to apply the method in all analyzed cases. Since the continent has an altered albedo, as well as an another critical temperature, the continent will have an slightly different ice-albedo dynamic to that of water. The resulting dynamic differs vastly from that of an aquaplanet, and this results in some very exotic stationary solutions. The introduction of a continent also ensures the appearance of multiple, new stationary solutions as we vary the bifurcation parameter, $Q$.

\vsone

\subsubsection{Meridional Symmetry}
\label{section: BIM symmetrical continent}
Initially, let's consider a meridionally symmetric continent configuration. We let the extent of the continent be
\begin{equation*}
    l = \frac{\pi}{4}.
\end{equation*}
The continent extends from $\theta_{l_1}$ to $\theta_{l_2}$ where
\begin{equation}
    \begin{aligned}
    \theta_{l_1} &= \frac{\pi}{2} - \frac{l}{2} = \frac{3\pi}{8}\\
    \theta_{l_2} &= \frac{\pi}{2} + \frac{l}{2} = \frac{5\pi}{8}.
    \end{aligned}
    \label{eqn: meridionally symmetric continent configuration}
\end{equation}
As we vary the control parameter, we find that the constraints of the system may be satisfied by some fascinating land-ice distributions on the continent. Due to the differing ice-albedo dynamic for water and land, there may be a partial ice cover on the two oceans and on the continent (see Figure \ref{fig: eps=0, 4 ice edges}). Furthermore, there may even be several bands of ice on the continent (see Figure \ref{fig: eps=0, 6 ice edges}).

\begin{figure}[H]
\centering
\begin{subfigure}{.5\textwidth}
  \centering
  \includegraphics[width=\linewidth]{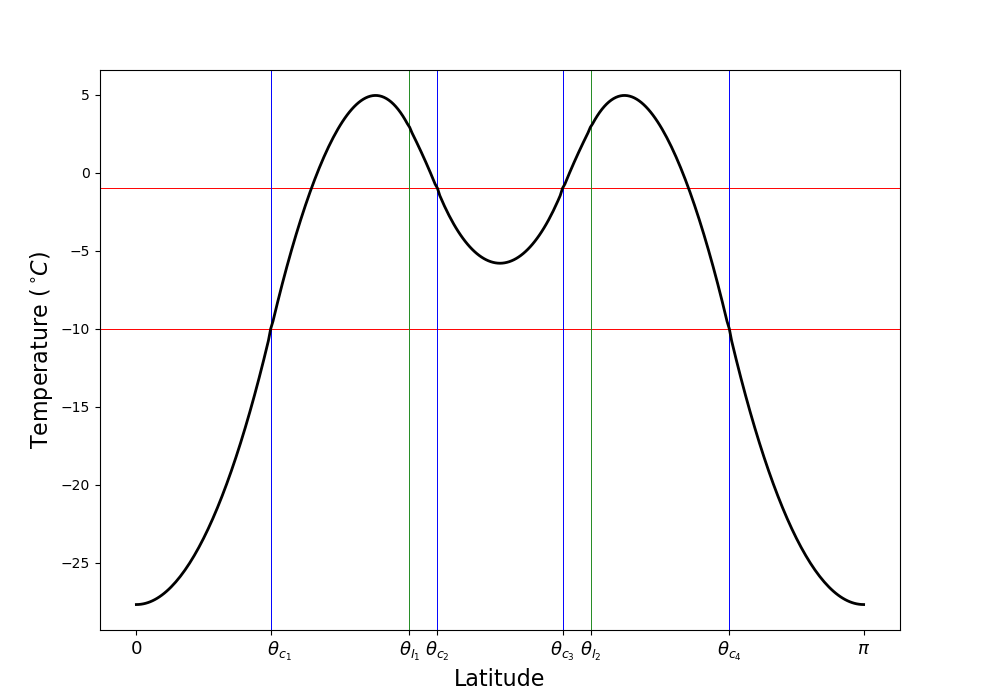}
  \caption{}
  \label{fig: plot of solution no ice edges no ice }
\end{subfigure}%
\begin{subfigure}{.5\textwidth}
  \centering
  \resizebox{0.5\textwidth}{!}{ 

\begin{tikzpicture}[scale=1,every node/.style={minimum size=1cm}]
    
    \def\R{6} 
    
    \def\angEl{10} 
    \def\angPhiThetaline{-110} 
    \def\angIceEdgeN{45} 
    \def\angIceEdgeS{-45} 
    \def\angContinentStart{157.5} 
    \def\angContinentEnd{202.5} 
    \def\angIceEdgeNOnCont{157.5 + 10} 
    \def\angIceEdgeSOnCont{202.5 - 10} 
    
    \pgfmathsetmacro\CN{\R*cos(\angEl)*sin(\angIceEdgeN)} 
    \pgfmathsetmacro\CS{\R*cos(\angEl)*sin(\angIceEdgeS)} 
    \pgfmathsetmacro\CNoncont{\R*cos(\angEl)*sin(\angIceEdgeNOnCont)} 
    \pgfmathsetmacro\CSoncont{\R*cos(\angEl)*sin(\angIceEdgeSOnCont)} 
    \pgfmathsetmacro\ContN{\R*cos(\angEl)*sin(\angContinentStart)} 
    \pgfmathsetmacro\ContS{\R*cos(\angEl)*sin(\angContinentEnd)} 
    \pgfmathsetmacro\H{\R*cos(\angEl)} 
    \LongitudePlane[Theta_line_plane]{\angEl}{\angPhiThetaline}
    \LatitudePlane[North_Ice_edge]{\angEl}{\angIceEdgeN}
    \LatitudePlane[South_Ice_edge]{\angEl}{\angIceEdgeS}
    \LatitudePlane[North_Ice_edge_on_cont]{\angEl}{\angIceEdgeNOnCont}
    \LatitudePlane[South_Ice_edge_on_cont]{\angEl}{\angIceEdgeSOnCont}
    \LatitudePlane[Continent_start]{\angEl}{\angContinentStart}
    \LatitudePlane[Continent_end]{\angEl}{\angContinentEnd}
    \fill[ball color=white!10] (0,0) circle (\R); 
    \coordinate (CN) at (0,\CN);
    \coordinate (CS) at (0,\CS);
    \coordinate (CNoncont) at (0,\CNoncont);
    \coordinate (CSoncont) at (0,\CSoncont);
    \coordinate (ContN) at (0,\ContN);
    \coordinate (ContS) at (0,\ContS);
    \coordinate[mark coordinate] (N) at (0,\H);
    \coordinate[mark coordinate] (S) at (0,-\H);

    \DrawLatitudeCircleRedNoDash[\R]{\angIceEdgeN}   
    \DrawLatitudeCircleRedNoDash[\R]{\angIceEdgeS}
    \DrawLatitudeCircleRedoncontNoDash[\R]{\angIceEdgeNOnCont}   
    \DrawLatitudeCircleRedoncontNoDash[\R]{\angIceEdgeSOnCont}
    \DrawLatitudeCircleThickNoDash[\R]{\angContinentStart}   
    \DrawLatitudeCircleThickNoDash[\R]{\angContinentEnd}
    \node[above=8pt] at (N) {$\mathbf{N}$};
    \node[below=8pt] at (S) {$\mathbf{S}$};

    \path[Theta_line_plane] (\angIceEdgeN:\R) coordinate (IceN);
    \path[Theta_line_plane] (\angIceEdgeS:\R) coordinate (IceS);
    \draw[Theta_line_plane, thick, black] ([shift={(90:\R)}]0,0) arc  (90:90+180:\R) ;
    \draw[North_Ice_edge, thick, black] ([shift={(-67:\R)}]CN) arc  (-67:-62:\R) node[near end, xshift=0.5cm, above=-7pt] {$\theta_{c_1}$};
    \draw[South_Ice_edge, thick, black] ([shift={(-81:\R)}]CS) arc  (-81:-76:\R) node[near end, xshift=0.5cm, above=-7pt] {$\theta_{c_4}$};
    
    \draw[North_Ice_edge_on_cont, thick, black] ([shift={(108:\R)}]CNoncont) arc  (108:113:\R) node[near end, xshift=0.5cm, above=-7pt] {$\theta_{c_2}$};
    \draw[South_Ice_edge_on_cont, thick, black] ([shift={(105:\R)}]CSoncont) arc  (105:110:\R) node[near end, xshift=0.5cm, above=-7pt] {$\theta_{c_3}$};
    
    \draw[Continent_start, very thick, black] ([shift={(109:\R)}]ContN) arc  (109:114:\R) node[near end, xshift=0.5cm, above=-7pt] {$\theta_{l_1}$};
    \draw[Continent_end, very thick, black] ([shift={(103:\R)}]ContS) arc  (103:108:\R) node[near end, xshift=0.5cm, above=-7pt] {$\theta_{l_2}$};

\end{tikzpicture}
}
  \caption{}
  \label{fig: plot of solution no ice edges snowball}
\end{subfigure}
\caption{A temperature profile (a) and a schematic of the planet (b) with four critical latitudes with the continent configuration (\ref{eqn: meridionally symmetric continent configuration}). Model parameters are those in Appendix \ref{section: Parameter values} and $Q=294 \: \text{W m}^{-2}$. Green vertical lines mark the continent borders and blue vertical lines mark critical latitudes. The red horizontal lines mark the critical temperatures for ice formation.}
\label{fig: eps=0, 4 ice edges}
\end{figure}
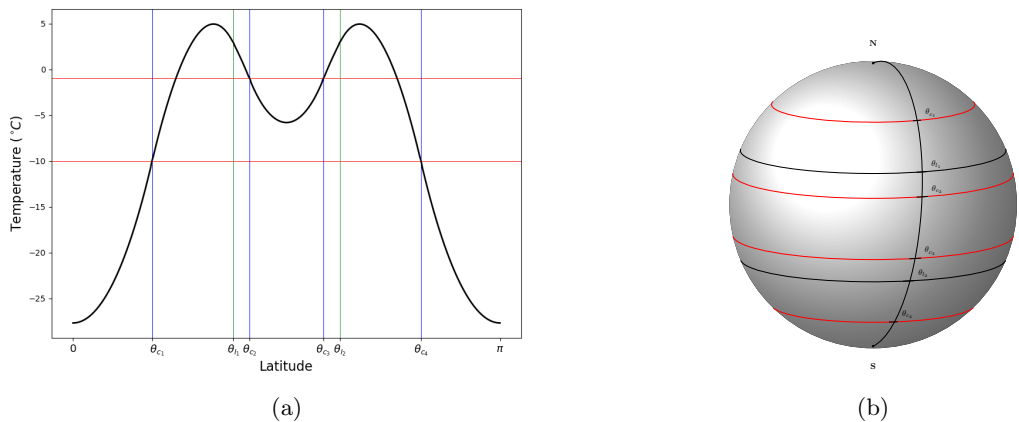

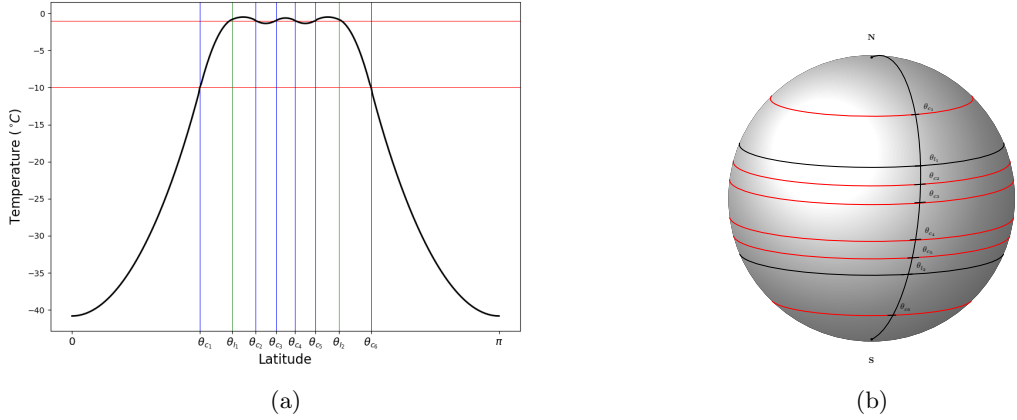
\begin{figure}[H]
\centering
\begin{subfigure}{.5\textwidth}
  \centering
  \includegraphics[width=\linewidth]{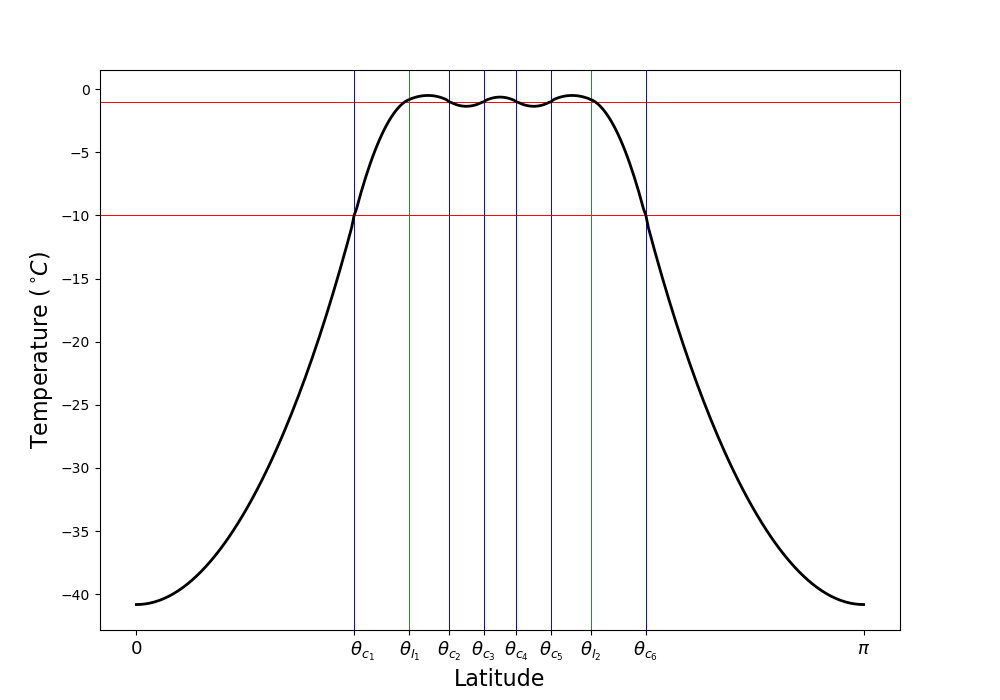}
  \caption{}
  \label{fig: plot of solution no ice edges no ice }
\end{subfigure}%
\begin{subfigure}{.5\textwidth}
  \centering
  \resizebox{0.5\textwidth}{!}{ 

\begin{tikzpicture}[scale=1,every node/.style={minimum size=1cm}]
    
    \def\R{6} 
    
    \def\angEl{10} 
    \def\angPhiThetaline{-110} 
    \def\angIceEdgeN{45} 
    \def\angIceEdgeS{-45} 
    \def\angContinentStart{157.5} 
    \def\angContinentEnd{202.5} 
    \def\angIceEdgeNOnCont{157.5 + 7.5} 
    \def\angIceEdgeSOnCont{202.5 - 7.5} 
    \def\angIceEdgeNOnContII{157.5 + 15} 
    \def\angIceEdgeSOnContII{202.5 - 15} 
    
    \pgfmathsetmacro\CN{\R*cos(\angEl)*sin(\angIceEdgeN)} 
    \pgfmathsetmacro\CS{\R*cos(\angEl)*sin(\angIceEdgeS)} 
    \pgfmathsetmacro\CNoncont{\R*cos(\angEl)*sin(\angIceEdgeNOnCont)} 
    \pgfmathsetmacro\CSoncont{\R*cos(\angEl)*sin(\angIceEdgeSOnCont)} 
    \pgfmathsetmacro\CNoncontII{\R*cos(\angEl)*sin(\angIceEdgeNOnContII)} 
    \pgfmathsetmacro\CSoncontII{\R*cos(\angEl)*sin(\angIceEdgeSOnContII)} 
    \pgfmathsetmacro\ContN{\R*cos(\angEl)*sin(\angContinentStart)} 
    \pgfmathsetmacro\ContS{\R*cos(\angEl)*sin(\angContinentEnd)} 
    \pgfmathsetmacro\H{\R*cos(\angEl)} 
    \LongitudePlane[Theta_line_plane]{\angEl}{\angPhiThetaline}
    \LatitudePlane[North_Ice_edge]{\angEl}{\angIceEdgeN}
    \LatitudePlane[South_Ice_edge]{\angEl}{\angIceEdgeS}
    \LatitudePlane[North_Ice_edge_on_cont]{\angEl}{\angIceEdgeNOnCont}
    \LatitudePlane[South_Ice_edge_on_cont]{\angEl}{\angIceEdgeSOnCont}
    \LatitudePlane[North_Ice_edge_on_contII]{\angEl}{\angIceEdgeNOnCont}
    \LatitudePlane[South_Ice_edge_on_contII]{\angEl}{\angIceEdgeSOnCont}
    \LatitudePlane[Continent_start]{\angEl}{\angContinentStart}
    \LatitudePlane[Continent_end]{\angEl}{\angContinentEnd}
    \fill[ball color=white!10] (0,0) circle (\R); 
    \coordinate (CN) at (0,\CN);
    \coordinate (CS) at (0,\CS);
    \coordinate (CNoncont) at (0,\CNoncont);
    \coordinate (CSoncont) at (0,\CSoncont);
    \coordinate (CNoncontII) at (0,\CNoncontII);
    \coordinate (CSoncontII) at (0,\CSoncontII);
    \coordinate (ContN) at (0,\ContN);
    \coordinate (ContS) at (0,\ContS);
    \coordinate[mark coordinate] (N) at (0,\H);
    \coordinate[mark coordinate] (S) at (0,-\H);

    \DrawLatitudeCircleRedNoDash[\R]{\angIceEdgeN}   
    \DrawLatitudeCircleRedNoDash[\R]{\angIceEdgeS}
    \DrawLatitudeCircleRedoncontNoDash[\R]{\angIceEdgeNOnCont}   
    \DrawLatitudeCircleRedoncontNoDash[\R]{\angIceEdgeSOnCont}
    \DrawLatitudeCircleRedoncontNoDash[\R]{\angIceEdgeNOnContII}   
    \DrawLatitudeCircleRedoncontNoDash[\R]{\angIceEdgeSOnContII}
    \DrawLatitudeCircleThickNoDash[\R]{\angContinentStart}   
    \DrawLatitudeCircleThickNoDash[\R]{\angContinentEnd}
    \node[above=8pt] at (N) {$\mathbf{N}$};
    \node[below=8pt] at (S) {$\mathbf{S}$};

    \path[Theta_line_plane] (\angIceEdgeN:\R) coordinate (IceN);
    \path[Theta_line_plane] (\angIceEdgeS:\R) coordinate (IceS);
    \draw[Theta_line_plane, thick, black] ([shift={(90:\R)}]0,0) arc  (90:90+180:\R) ;
    \draw[North_Ice_edge, thick, black] ([shift={(-67:\R)}]CN) arc  (-67:-62:\R) node[near end, xshift=0.5cm, above=-7pt] {$\theta_{c_1}$};
    \draw[South_Ice_edge, thick, black] ([shift={(-81:\R)}]CS) arc  (-81:-76:\R) node[near end, xshift=0.5cm, above=-7pt] {$\theta_{c_6}$};
    
    \draw[North_Ice_edge_on_cont, thick, black] ([shift={(108:\R)}]CNoncont) arc  (108:113:\R) node[near end, xshift=0.5cm, above=-7pt] {$\theta_{c_2}$};
    \draw[South_Ice_edge_on_cont, thick, black] ([shift={(105:\R)}]CSoncont) arc  (105:110:\R) node[near end, xshift=0.5cm, above=-7pt] {$\theta_{c_5}$};
    
    \draw[North_Ice_edge_on_contII, thick, black] ([shift={(108:\R)}]CNoncontII) arc  (108:113:\R) node[near end, xshift=0.5cm, above=-7pt] {$\theta_{c_3}$};
    \draw[South_Ice_edge_on_contII, thick, black] ([shift={(106:\R)}]CSoncontII) arc  (106:111:\R) node[near end, xshift=0.5cm, above=-7pt] {$\theta_{c_4}$};
    
    \draw[Continent_start, very thick, black] ([shift={(109:\R)}]ContN) arc  (109:114:\R) node[near end, xshift=0.5cm, above=-7pt] {$\theta_{l_1}$};
    \draw[Continent_end, very thick, black] ([shift={(103:\R)}]ContS) arc  (103:108:\R) node[near end, xshift=0.5cm, above=-7pt] {$\theta_{l_2}$};

\end{tikzpicture}
}
  \caption{}
  \label{fig: plot of solution no ice edges snowball}
\end{subfigure}
\caption{A temperature profile (a) and a schematic of the planet (b) with six critical latitudes with the continent configuration (\ref{eqn: meridionally symmetric continent configuration}). Model parameters are those in Appendix \ref{section: Parameter values} and $Q=299 \: \text{W m}^{-2}$. Green vertical lines mark the continent borders and blue vertical lines mark critical latitudes. The red horizontal lines mark the critical temperatures for ice formation.}
\label{fig: eps=0, 6 ice edges}
\end{figure}

\subsubsection{Continent Breaking the Meridional Symmetry}
\label{section: Continent breaking the north-south symmetry}
\begin{figure}[H]
        \begin{subfigure}{.5\textwidth}
          \centering
          \includegraphics[width=\linewidth]{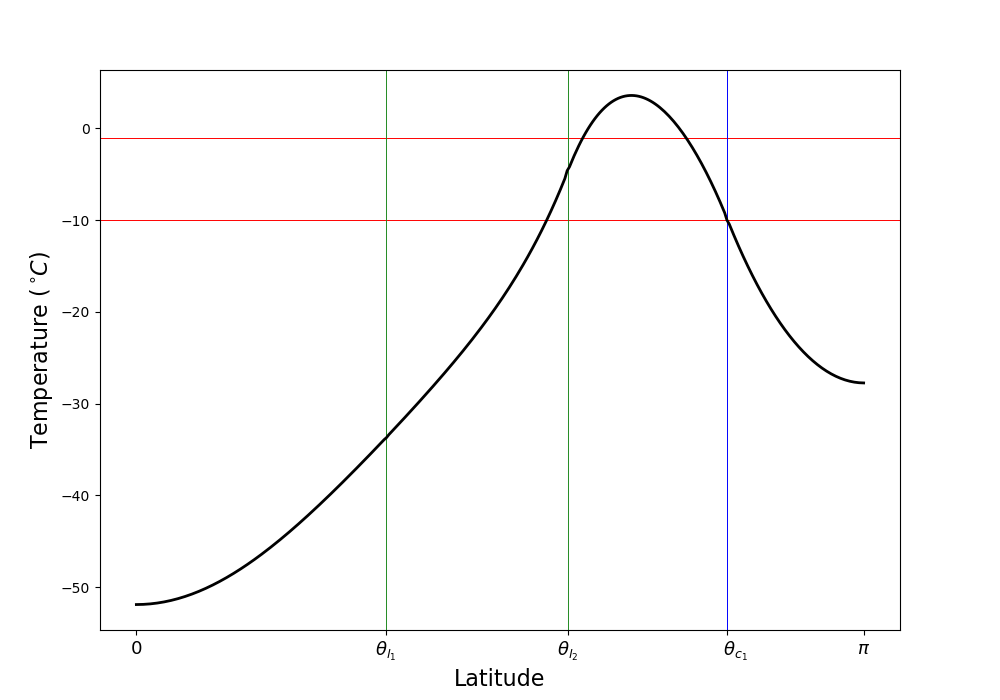}
          \caption{}
          \label{fig: plot of solution eps=0.1 one ice edge a}
        \end{subfigure}%
        \begin{subfigure}{.5\textwidth}
          \centering
          \includegraphics[width=\linewidth]{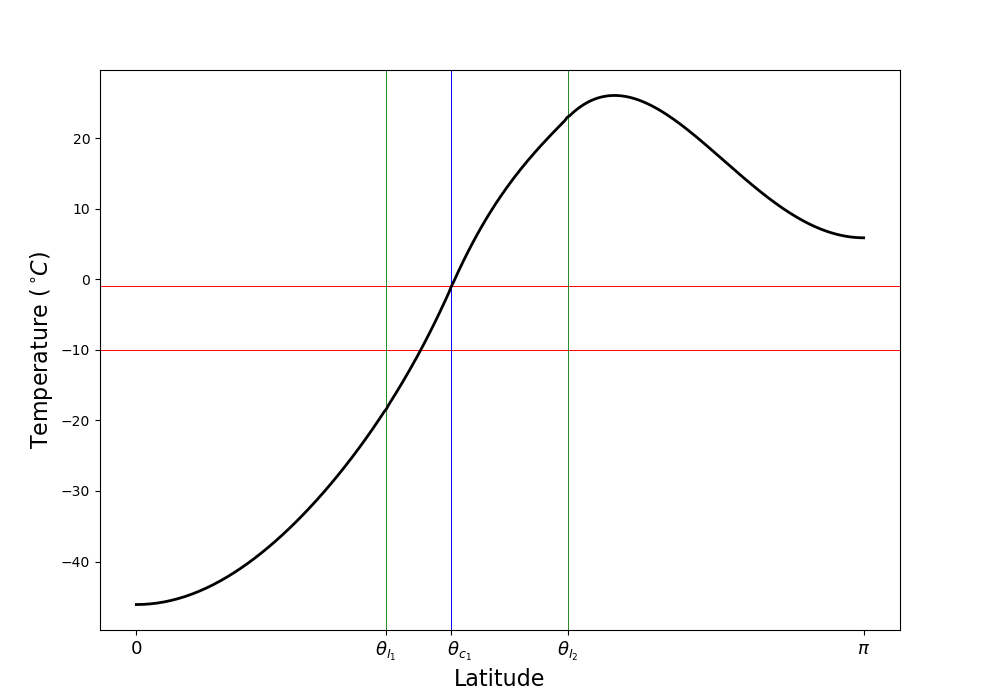}
          \caption{}
          \label{fig: plot of solution eps=0.1 one ice edge b}
        \end{subfigure}
    \caption{Temperature profiles with one critical latitude on the southern ocean (a) and one critical latitude on the continent (b) for a planet with a continent configuration of the kind (\ref{eqn: continent configuration}) and $\varepsilon = 0.1$. Model parameters are those in Appendix \ref{section: Parameter values}. In (a) $Q=300 \: \text{W m}^{-2}$ and in (b) $Q=310 \: \text{W m}^{-2}$. Green vertical lines mark the continent borders and blue vertical lines mark critical latitudes. The red horizontal lines mark the critical temperatures for ice formation. Note that in (a) both the continent and the northern ocean is completely ice covered, however, the continental edge, $\theta_{l_2}$, is not considered a critical latitude since the condition (\ref{eqn: critical latitude condition on continent}) does not hold.}
    \label{fig: plot of solution eps=0.1}
    \end{figure}

The aim of this study is to investigate how symmetry affect the model, it therefore makes sense to initially consider how the system respond to a minor symmetry violation. Asymmetry will be introduced through slightly shift the continent to the north, thus breaking the meridional symmetry. We let the extent of continent remain at $l = \frac{\pi}{4}$, but the continent will now extend from the latitude
\begin{equation}
    \begin{aligned}
        \theta_{l_1} &= \frac{\pi}{2} - \frac{l}{2} - \varepsilon\\
        \theta_{l_2} &= \frac{\pi}{2} + \frac{l}{2} - \varepsilon,
    \end{aligned}
    \label{eqn: continent configuration}
\end{equation}
where $\varepsilon = 0.1$. Applying the BIM to find stationary solutions to (\ref{eqn: time dependent EBM w continent}) with this continent configuration will be analogous to the case considered above. However, the asymmetrical continent configuration is expected to result in an asymmetrical temperature distributions and the appearance of stationary solutions with an odd number of critical latitudes. Moreover, the fine dynamic resulting in stationary solutions with numerous critical latitudes (four and six) has disappeared follow the continent shift. Figure \ref{fig: plot of solution eps=0.1} shows two examples of the asymmetrical temperature profiles resulting from the continent configuration (\ref{eqn: continent configuration}) and $\varepsilon = 0.1$.

\subsubsection{Increasing the Asymmetry}
\label{section: BIM eps=0.5}
\begin{figure}[H]
        \begin{subfigure}{.5\textwidth}
          \centering
          \includegraphics[width=\linewidth]{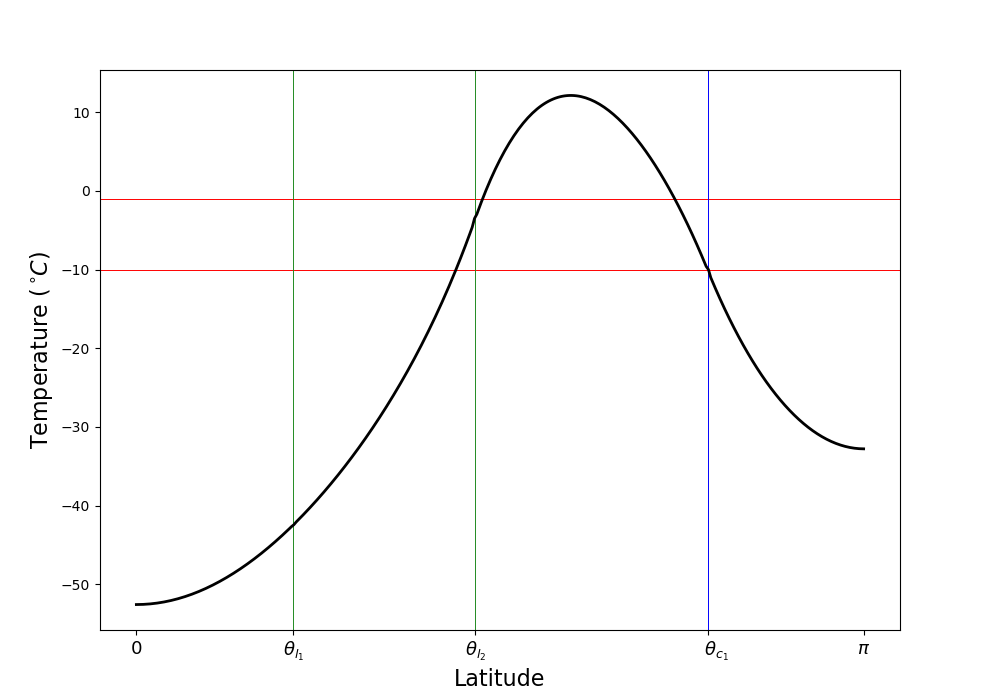}
          \caption{}
          \label{fig: plot of solution eps=0.1 one ice edge a}
        \end{subfigure}%
        \begin{subfigure}{.5\textwidth}
          \centering
          \includegraphics[width=\linewidth]{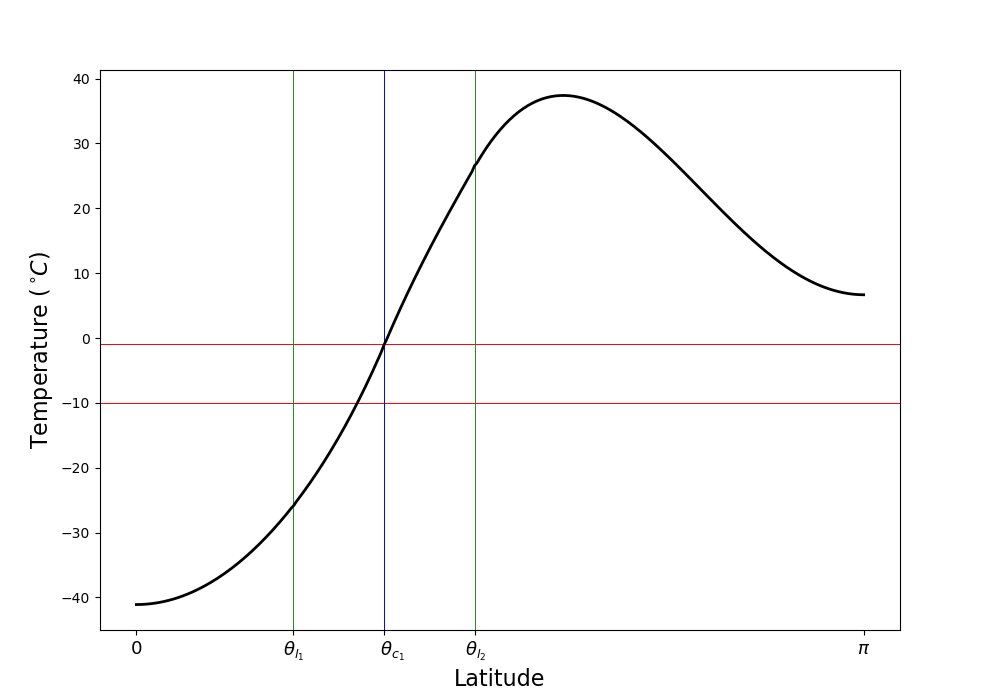}
          \caption{}
          \label{fig: plot of solution eps=0.1 one ice edge b}
        \end{subfigure}
    \caption{A temperature profile with one critical latitudes (a) and a temperature profile with one critical latitude (b) for a planet with a continent configuration of the kind (\ref{eqn: continent configuration}) and $\varepsilon = 0.5$. Model parameters are those in Appendix \ref{section: Parameter values}. In (a) $Q=267 \text{Wm}^{-2}$ and in (b) $Q=300 \: \text{Wm}^{-2}$. Green vertical lines mark the continent borders and blue vertical lines mark critical latitudes. The red horizontal lines mark the critical temperatures for ice formation.}
    \label{fig: plot of solution eps=0.5}
    \end{figure}

To further investigate how asymmetry affect the model, the continent is shifted further north. We let the continent configuration be like (\ref{eqn: continent configuration}), with $\varepsilon = 0.5$. Naturally, this continent configuration also gives rise to asymmetrical temperature profiles. Two such temperature profiles are displayed in Figure \ref{fig: plot of solution eps=0.5}.

\clearpage

\section{Finite Difference Method}
\label{section: Finite difference method}
In this section we apply a finite difference method to solve the nondimensional time-dependent energy balance equation from section \ref{section: Nondimensionalize the model},
\begin{equation}
    \gamma \partial_t T + \mathcal{L}T + \beta T= \eta s(\theta)(1-a(T)) - \alpha,
    \label{eqn: nondimensional time dependent EBM}
\end{equation}
where $\mathcal{L}$ is the operator (\ref{eqn: operator with beta not appendix}) and the dimensionless parameters are defined in (\ref{eqn: dimensionless parameters}). The implemented finite difference algorithm is then validated by an \textit{artificial source test}. The albedo, $a(T)$, may have an arbitrary temperature and/or latitude dependence in the following finite difference scheme. We will use the smooth albedo function is given in Appendix \ref{section: Parameter values}.

\subsection{Implementing a Finite Difference Algorithm}
\label{section: Implement a numerical solution to the time dependent energy balance model}
Finite difference methods relies on approximating differential operators using finite differences. Appendix \ref{section: Approximate the differential operator} includes a derivation of a centered difference formula for operator (\ref{eqn: operator with beta not appendix}) applied to a smooth function, $f$. For the discretized function evaluated on a uniform grid, $x_i$, where $f_i = f(x_i)$, a discrete approximation of $\mathcal{L}f$ will be
\begin{equation}
    \hat{\mathcal{L}}f_i = -\frac{2(f_{i-1}-2f_i + f_{i+1}) + (f_{i-1}-4f_i + 3f_{i+1})h\cot x_i}{2h^2}.
    \label{eqn: finite diff approximation of L}
\end{equation}
Applying this approximation to the discretized form of (\ref{eqn: nondimensional time dependent EBM}) turns the problem into a system of ordinary differential equations. The following algorithm will apply to an aquaplanet. Necessary modifications to the algorithm to include a continent will also be discussed.

\vsone

Partition the domain into a uniform angular-time grid,
\begin{align*}
    \theta_{i+1}&=\theta_i + h,\\
    \theta_0 &= 0,\\
    \theta_N &= \pi
    \intertext{and}
    t_{j+1}&=t_j +k.
\end{align*}
Introducing the notation 
\begin{equation*}
    T_{i,j}=T(\theta_i, t_j),
\end{equation*}
a discrete form of (\ref{eqn: nondimensional time dependent EBM}) is
\begin{equation}
     \gamma \partial_t T_{i,j} + \hat{\mathcal{L}}T_{i,j} + \beta T_{i,j}= \eta s(\theta_i)a_c(T_{i,j}) - \alpha.
     \label{eqn: discret non-domensional time dependent EBM}
\end{equation}
Assuming that $T(\theta, t)$ is smooth throughout the domain, the finite difference approximation (\ref{eqn: finite diff approximation of L}) is applied to approximate the term $\hat{\mathcal{L}}T$ in equation (\ref{eqn: discret non-domensional time dependent EBM}) such that
\begin{equation}
\begin{aligned}
         \gamma \partial_t T_{i,j} -\frac{2(T_{i-1, j}-2T_{i,j} + T_{i+1, j}) + (T_{i-1, j}-4T_{i, j} + 3T_{i+1, j})h\cot \theta_i}{2h^2}&\\
         + \: \beta T_{i,j} &= \eta s(\theta_i)a_c(T_{i,j}) - \alpha.
\end{aligned}
     \label{eqn: discret non-domensional time dependent EBM 2}
\end{equation}
Equation (\ref{eqn: discret non-domensional time dependent EBM 2}) paired with an initial condition,
\begin{equation*}
    T(\theta, t=0) = T_{IC}(\theta),
\end{equation*}
gives rise to the system of initial value problems
\begin{equation}
    \begin{aligned}
            \parD{}{t}T_{i,j}&= f(\theta_i, T_{i-1, j}, T_{i,j}, T_{i+1,j}), \hstwo i =1,2,3, ..., N-1,\\
        T_{i,0}&=T_{IC}(\theta_i),
    \end{aligned}
    \label{eqn: sys of IVPs}
\end{equation}
where
\begin{equation*}
\begin{aligned}
            f(\theta_i, T_{i-1, j}, T_{i,j}, T_{i+1,j}) = \frac{1}{\gamma}\bigg[&\frac{2(T_{i-1, j}-2T_{i,j} + T_{i+1, j}) + (T_{i-1, j}-4T_{i, j} + 3T_{i+1, j})h\cot \theta_i}{2h^2} \\
            &- \beta T_{i,j} + \eta s(\theta_i)a_c(T_{i,j}) - \alpha \bigg].
\end{aligned}
\end{equation*}
The values at the endpoints of grid, $T_{0, j}$ and $T_{N, j}$ are required to fully specify the solution. These may be obtained by invoking the boundary conditions. Recall that 
\begin{equation*}
    \lim_{\theta \rightarrow 0} \sin \theta \parD{T}{\theta}(\theta)  =0.
\end{equation*}
Assuming that $h$ is arbitrarily small, at the second grid point, $\theta_1 = h$, we must have
\begin{equation*}
    \sin \theta_1 \parD{T}{\theta}(\theta_1)=0.
\end{equation*}
Applying the centered difference formula for first order derivatives, we can conclude that
\begin{equation*}
    \sin h \frac{T_2 -T_0}{2h}=0.
\end{equation*}
For an arbitrarily small $h<<1$,
\begin{equation*}
    \sin h \approx h.
\end{equation*}
Consequently,
\begin{equation*}
    \frac{T_2 -T_0}{2}=0 \hsone \Longrightarrow \hsone T_0 = T_2.
\end{equation*}
Following a similar approach for the boundary condition at $\theta=\pi$, we get the analogous condition
\begin{equation*}
    T_N = T_{N-2}.
\end{equation*}
Thus, the boundary conditions (\ref{eqn: BC theta=0}) and  (\ref{eqn: BC theta=pi}) gives rise to the following rules in the numerical scheme;
\begin{equation}
    \begin{aligned}
        T_{0, j} &= T_{2, j}\\
        T_{N, j} &= T_{N-2, j}.
    \end{aligned}
    \label{eqn: rule at either end of grid}
\end{equation}
The system of initial value problems (\ref{eqn: sys of IVPs}), paired the rules at either end of the angular-time grid (\ref{eqn: rule at either end of grid}), is easily solved by any IVP-solver. Between every time step, the boundary conditions (\ref{eqn: rule at either end of grid}) is invoked allowing the algorithm to fill the angular-time grid. This algorithm is easily modified to include a continent: The model parameters that differ on land and water are changed on the interval of the domain corresponding to the continent. This is done by altering the inhomogeneous term in the system of IVP's (\ref{eqn: sys of IVPs}), $f(\theta_i, T_{i-1, j}, T_{i,j}, T_{i+1,j})$, for every $i$ such that $\theta_i \in [\theta_{l_1}, \theta_{l_2}]$ .

\vsone

\vsone

\subsection{Artificial Source Test}
We will now employ an artificial source test to validate the numerical solution from above. The artificial source test serves as a general method for assessing the accuracy of a numerical solution to a differential equation. The prescribed procedure for conducting an artificial source test is as follows: Consider a general differential equation in the form

\begin{equation}
\mathcal{G}f(\vv{x}, t) = h(\vv{x}, t),
\label{eqn: general DE}
\end{equation}

where $\mathcal{G}$ is a differential operator. Assuming the implementation of a numerical method to solve (\ref{eqn: general DE}), we introduce an exact solution

\begin{equation*}
f_e(\vv{x}, t),
\end{equation*}

such that the exact function

\begin{equation*}
\rho(\vv{x}, t) = \mathcal{G}f_e(\vv{x}, t)
\end{equation*}

is known. This exact solution serves as the "artificial source." Subsequently, we allow the implemented numerical method to solve the modified problem

\begin{equation}
\mathcal{G}f(\vv{x}, t) = \rho(\vv{x}, t).
\label{eqn: related general DE}
\end{equation}

The resulting numerical solution, denoted as $f(\vv{x}, t)$, should align perfectly with the assumed exact solution $f_e(\vv{x}, t).$ A successful match validates the implemented numerical method, as the exact solution to (\ref{eqn: related general DE}) is precisely $f_e(\vv{x}, t)$. The chosen assumed exact solution must also adhere to any constraints or conditions imposed on the solution of the original problem (\ref{eqn: general DE}).

\vsone

Multiple artificial sources are applied to test the finite difference method. Assuming that

\begin{equation*}
T_e = T_e(\theta, t),
\end{equation*}

is a solution to (\ref{eqn: nondimensional time dependent EBM}) that adherence to the boundary conditions (\ref{eqn: BC theta=0}) and (\ref{eqn: BC theta=pi}). The associated artificial source is derived,

\begin{equation*}
\rho(\theta, t) = \gamma \partial_t T_e + \mathcal{L}T_e + \beta T_e.
\end{equation*}

This artificial source is an exact function, and the finite difference algorithm is employed to solve the modified problem

\begin{align*}
\gamma \partial_t T + \mathcal{L}T + \beta T &= \rho(\theta, t)\\
\lim_{\theta \rightarrow 0} \sin \theta \parD{T}{\theta}(\theta) &=0\\
\lim_{\theta \rightarrow \pi} \sin \theta \parD{T}{\theta}(\theta) &=0.
\end{align*}

Various assumed solutions were tested, including once based on physical intuition and others specifically chosen to assess the robustness of the implemented code. Two representative assumed solutions, along with the corresponding results from the artificial source test, are presented below in Figures \ref{fig: artificial source test 1} and \ref{fig: artificial source test 2}. In Figure \ref{fig: artificial source test 1}, the assumed solution is given by
\begin{align}
T_e(\theta, t) &= (3 - \sin(t_0t))e^{-3(\theta - \frac{\pi}{2})^2}-1
\label{eqn: assumed solution 1}
\intertext{and in Figure \ref{fig: artificial source test 2}, the assumed solution is}
T_e(\theta, t) &= 2e^{-5(\theta - \frac{\pi}{2}\cos(0.2t_0t))^2}.
\label{eqn: assumed solution 2}
\end{align}

The chosen assumed solutions (\ref{eqn: assumed solution 1}) and (\ref{eqn: assumed solution 2}) are non-dimensional and designed to yield dimensional equivalents, $T_s T_e$, within a plausible range based on the current state of Earth's climate. These solutions also vary on the timescale $t_0$.

\begin{figure}[H]
\centering
\begin{subfigure}{.5\textwidth}
\centering
\hspace{-.75cm}
\includegraphics[width=1.1\linewidth]{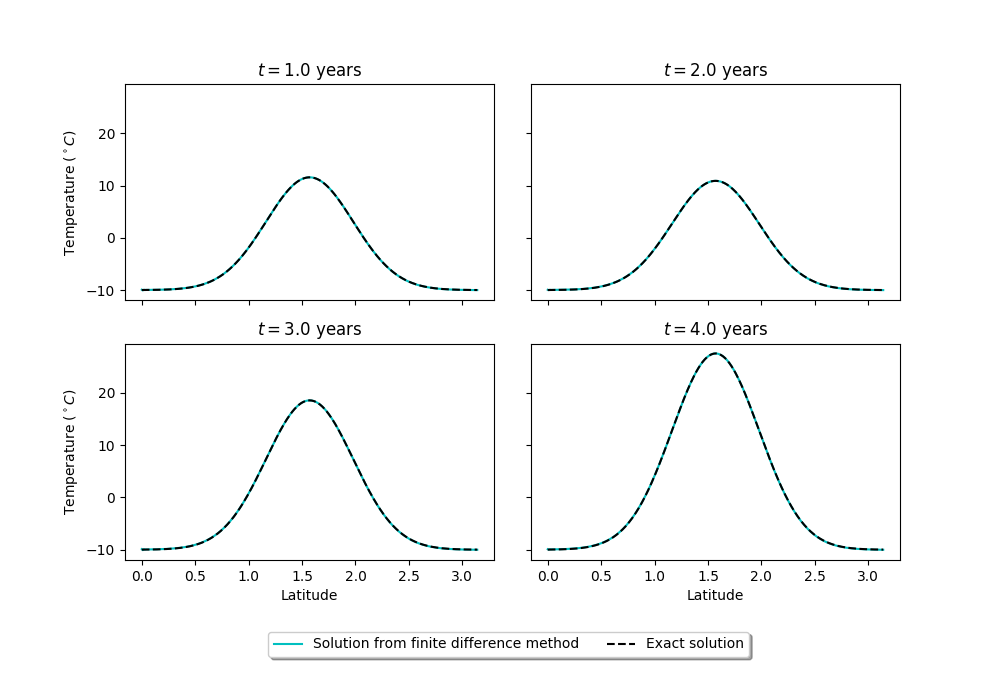}
\caption{}
\label{fig: artificial source test 1}
\end{subfigure}%
\begin{subfigure}{.5\textwidth}
\centering
\includegraphics[width=1.1\linewidth]{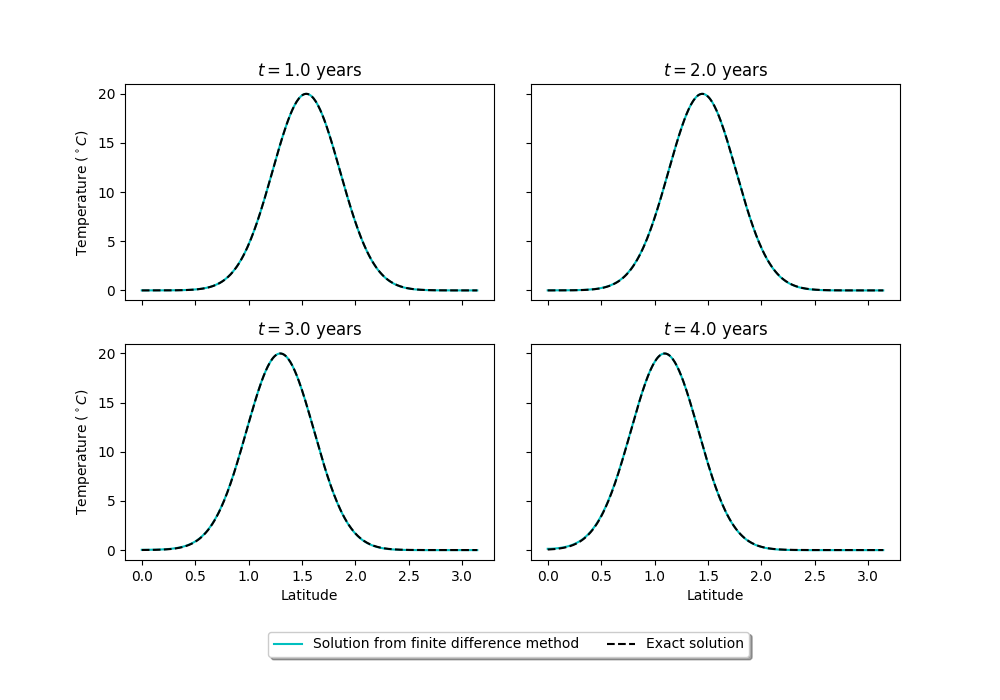}
\caption{}
\label{fig: artificial source test 2}
\end{subfigure}
\caption{Artificial source test for the assumed solution (\ref{eqn: assumed solution 1}) (a) and (\ref{eqn: assumed solution 2}) (b). }
\label{fig: plot of solution no edges}
\end{figure}

\clearpage

\section{Stability Analysis}
\label{section: stability analysis}
In this section, we assess the stability properties of the stationary solutions from section \ref{section: Boundary Integral Method}. Stability properties are inferred using up-to three techniques: A numerical perturbation scheme, a \textit{slope-stability} theorem and, if necessary, a heuristic approach. These techniques will be discussed below, starting with the perturbation scheme.

\subsection{Numerical Perturbation Scheme}
Considering the time-dependent energy balance equation
\begin{equation}
    \gamma \partial_t T -\frac{1}{\sin\theta}\parD{}{\theta }\big( \sin\theta \parD{T}{\theta} \big) + \beta T= \eta s(\theta)(1-a(T)) - \alpha.
    \label{eqn: time dep ebm, general form}
\end{equation}
Applying the notation from section \ref{section: Boundary Integral Method} we may rewrite equation (\ref{eqn: time dep ebm, general form}) as
\begin{equation}
    \gamma \partial_t T + \mathcal{L}T = h,
    \label{eqn: time dep EBM compact notation}
\end{equation}
where
\begin{align*}
    h&=h(T, \theta)
    \intertext{and}
    \mathcal{L}(\cdot) &= -\frac{1}{\sin\theta}\parD{}{\theta} \big( \sin\theta \parD{}{\theta} (\cdot) \big)     + \beta (\cdot).
\end{align*}
In section \ref{section: Boundary Integral Method} a set of stationary solutions was found,
\begin{equation*}
    T_0 = T(\theta, t=0),
\end{equation*}
such that
\begin{equation*}
    \mathcal{L}T_0 = h \hsone \implies \hsone \partial_t T_0=0.
\end{equation*}
The stationary solution, $T_0$, is given a small perturbation, 
\begin{equation*}
    T(\theta, t) = T_0(\theta) + \delta(\theta, t).
\end{equation*}
Stability properties of the stationary solution, $T_0$, are assessed by consideration of the following problem: Does the perturbation, $\delta$, grow in time or relax to zero keeping the dynamical solution in a neighborhood around the equilibrium?

\vsone

Substituting the perturbed stationary solution back into equation (\ref{eqn: time dep EBM compact notation}) we get
\begin{align}
    \gamma  \partial_t (T_0 + \delta) + \mathcal{L}(T_0 + \delta) &= h(T_0 + \delta, \theta)
    \nonumber\\
    \gamma \partial_t \delta + \mathcal{L}T_0 + \mathcal{L}\delta &= h(T_0 + \delta, \theta).
    \label{eqn: time dep EBM copact form w zeta}
\end{align}

The perturbation is initially small, therefore
\begin{equation*}
    \delta(\theta, t) \approx 0
\end{equation*}
for small $t$. Taylor expanding $h(T_0 + \delta, \theta)$ around $(T_0, \theta)$ we have
\begin{equation*}
    h(T_0 + \delta, \theta) \approx h(T_0, \theta) + h_T(T_0, \theta)\delta,
\end{equation*}
neglecting higher order terms since $\delta$ is initially small. Here
\begin{equation*}
    h_T(T_0, \theta) = \parD{h}{T}(T_0, \theta)
\end{equation*}
and can be found analytically for a given $T_0$. Substituting this $h(T_0+ \delta, \theta)$ into (\ref{eqn: time dep EBM copact form w zeta}) we get
\begin{align}
    \gamma \partial_t \delta + \mathcal{L}T_0 + \mathcal{L}\delta &= h(T_0, \theta) + h_T(T_0, \theta)\delta \nonumber\\
    \partial_t \delta &= \mathcal{H}\delta,
    \label{eqn: equation with operator H}
    \intertext{where}
    \mathcal{H}(\cdot) &= \frac{1}{\gamma} \bigg[h_T(T_0, \theta)(\cdot) - \mathcal{L}(\cdot) \bigg].
    \nonumber
\end{align}

Suppose $\delta$ has the form
\begin{equation}
    \delta (\theta, t) = e^{\lambda t}\delta_0(\theta).
    \label{eqn: form of zeta}
\end{equation}
Substitude this $\delta$ into (\ref{eqn: equation with operator H}) we get the following eigenvalue problem:
\begin{equation}
        \lambda\delta_0 = \mathcal{H}\delta_0
    \label{eqn: eigenvalue problem}
\end{equation}
If the eigenvalue $\lambda$ is real and positive, the perturbation (\ref{eqn: form of zeta}) will grow exponentially and the stationary solution $T_0$ is unstable.

\vsone
We will solve this eigenvalue problem (\ref{eqn: eigenvalue problem}) with a numerical method that relays on discretizing $\delta_0$ and $\mathcal{H}\delta_0$. We introduce a uniform angular grid for $\theta$,
\begin{align*}
    \theta_{i+1} &=\theta_i + d\theta,
    \intertext{where}
    \theta_0 &= 0\\
    \theta_N &= \pi
\end{align*}
Let
\begin{align*}
    T_0^{i}&=T_0(\theta_i)
    \intertext{and}
    \delta_0(\theta_i)&=\delta_0^{i}.
    \intertext{The discrete form of $\delta_0$, that is $\delta_0$ evaluated on the grid, is evidently a vector, }
    \bm{\delta_0} &= \begin{bmatrix} \delta_0^0 \\[1.5ex] \delta_0^1\\[1.5ex]  \vdots\\[1.5ex] \delta_0^N\end{bmatrix}.
    \end{align*}
Let
\begin{equation*}
    f^{i} = \big( \mathcal{H}\delta_0 \big)^{i} = \frac{1}{\gamma}\bigg[h_T(T_0^{i}, \theta_{i})\delta_0^{i} - \hat{\mathcal{L}}\delta_0^{i} \bigg].
\end{equation*}
be a discrete approximation to the function $\mathcal{H}\delta_0$. We can now approximate the eigenvalue problem (\ref{eqn: eigenvalue problem}) as
\begin{equation}
    \lambda \bm{\delta_0} = H \bm{\delta_0},
    \label{eqn: discrete eigenvalue problem}
\end{equation}
where the set of linear equations
\begin{equation}
    \big\{\lambda \delta_0^{i} = f^{i}\big\}_{i=0}^N
    \label{eqn: system of equations f^i}
\end{equation}
gives rise to the coefficient matrix $H$. Solving the discretized eigenvalue problem (\ref{eqn: discrete eigenvalue problem}) is a numerical exercise where we are looking for eigenvalues, $\lambda_j$, where $ \operatorname{Re}\{\lambda_j\} > 0$. If the matrix $H$, associated with a stationary solution, $T_0$, has at least one eigenvalue with a positive real part, the stationary solution $T_0$ is deemed to be unstable.

\vsone
In order to build the matrix $H$, a discrete approximation, $f^{i}$, of the function $\mathcal{H}\delta_0$ is necessary. Since $h_T(T, \theta)$ is an analytical function, only an approximation of $\mathcal{L}\delta_0$ is needed. This will be obtained using finite differences. For $i = 1, 2, 3,... , N-1$ the centered difference approximation of the operator $\mathcal{L}$ from section \ref{section: Finite difference method} may be applied,
\begin{align}
    f^{i} &= \frac{1}{\gamma}\bigg[h_T(T_0^{i}, \theta_{i})\delta_0^{i} - \hat{\mathcal{L}}\delta_0^{i} \bigg],
    \label{eqn: f^i}
    \intertext{where}
    \hat{\mathcal{L}}\delta_0^{i} &= \beta \delta_0^{i} - \frac{2(\delta_0^{i-1}-2\delta_0^{i} + \delta_0^{i+1}) + (\delta_0^{i-1}-4\delta_0^{i} + 3\delta_0^{i+1})d\theta\cot \theta_i}{2d\theta^2}.\nonumber
\end{align}
For $f^0$ a forward difference approximation is required. Likewise, a backward difference approximation is required for $f^N$. Appendix \ref{section: Approximate the differential operator} includes a derivation of both a forward and a backward difference approximation of the operator $\mathcal{L}$. With the appropriate approximations of $\mathcal{L}\delta_0$, $f^0$ and $f^N$ takes the form;
 \begin{align*}
    f^{0} &= \frac{1}{\gamma}\bigg[h_T(T_0^{0}, \theta_{0})\delta_0^{0} - \hat{\mathcal{L}}_f\delta_0^{0} \bigg],\\
    f^{N} &= \frac{1}{\gamma}\bigg[h_T(T_0^{N}, \theta_{N})\delta_0^{N} - \hat{\mathcal{L}}_b\delta_0^{N} \bigg],
    \intertext{where}
    \hat{\mathcal{L}}_f\delta_0^{0} &= \beta \delta_0^{0} + \frac{-2(\delta_0^{0}-2\delta_0^{1} + \delta_0^{2}) + (5\delta_0^{0}-8\delta_0^{1} + 3\delta_0^{2})d\theta\cot \theta_0}{2d\theta^2}
    \intertext{and}
    \hat{\mathcal{L}}_b\delta_0^{N} &= \beta \delta_0^{N} + \frac{-2(\delta_0^{N-2}-2\delta_0^{N-1} + \delta_0^{N}) + (\delta_0^{N-2}-\delta_0^{N})d\theta\cot \theta_N}{2d\theta^2}.
\end{align*}

\vsone

As $t$ grows the perturbation (\ref{eqn: form of zeta}) must necessarily follow the same boundary conditions as the solution. This is only possible if the discretized form of the function $\delta_0(\theta)$, $\delta_0^{i}$, follow the discrete boundary conditions (\ref{eqn: rule at either end of grid}) from section \ref{section: Implement a numerical solution to the time dependent energy balance model}. Therefore, $\delta_0^{i}$ must satisfy the following condition;
\begin{equation}
    \begin{aligned}
        \delta_{0}^{0} &= \delta_{0}^{2}\\
        \delta_{0}^{N} &= \delta_{0}^{N-2}.
    \end{aligned}
    \label{eqn: BC for zeta}
\end{equation}

The approximation $f^{i}$ will generally have the form
\begin{equation*}
    f^{i} = f^{i}_{1} \delta_0^{i-1} + f^{i}_{2} \delta_0^{i} + f^{i}_{3} \delta_0^{i+1},
\end{equation*}
but for every $f^{i}$ involving $\delta_{0}^{2}$ and $\delta_{0}^{N-2}$ we must impose condition (\ref{eqn: BC for zeta}), thus eliminating $\delta_{0}^{2}$ and $\delta_{0}^{N-2}$ from the system of equations (\ref{eqn: system of equations f^i}). With these expression for $f^{i}$ the matrix $H$ takes the form

\vspace{0.3cm}
\begin{equation*}
    H = \begin{bmatrix} 
    f^{0}_1 &   f^{0}_2     & 0         & 0         & 0     &0          &0         &0             &&&&&&&0\\[1ex]
    f^{1}_1 &   f^{1}_2     & 0         & 0         & 0     &0          &0         &0            &&&&&&&0\\[1ex]
    f^{2}_1 &   f^{2}_2     & 0         & f^{2}_3   & 0         &0         &0      &0    &&&&\hdots&&&0\\[1ex]
    f^{3}_1 &   0           & 0         &f^{3}_2    & f^{3}_3   &0         &0       &0              &&&&&&&0\\[1ex]
    0       &    0          &    0      &f^{4}_1    & f^{4}_2   & f^{4}_3 & 0       &0       &&&&&&&0\\[1ex]
    0       &    0          &    0      &0          &f^{5}_1    & f^{5}_2 & f^{5}_3 &0         &&&&&&&0\\[3ex]
    &&&&&&& \ddots &&&&&&&\\[3ex]
    0&&&&&&&           0&f^{N-5}_{1}   & f^{N-5}_{2}& f^{N-5}_{3}   & 0                 & 0                 & 0     & 0\\[1ex]
    0&&&&&&&            0&0             &f^{N-4}_{1} & f^{N-4}_{2}   & f^{N-4}_{3}       & 0                 & 0                  & 0\\[1ex]
    0&&&&&&&            0&0            &0             & f^{N-3}_{1}   & f^{N-3}_{2}       & 0                 & 0                 & f^{N-3}_{3} \\[1ex]
    0&&&&\hdots&&&       0&0            &0             & 0             & f^{N-2}_{1}       & 0                 & f^{N-2}_{2}       & f^{N-2}_{3}\\[1ex]
    0&&&&&&&            0&0           &0           &0             &0               & 0           & f^{N-1}_{1}       & f^{N-1}_2\\[1ex]
    0&&&&&&&              0&0          &0             &0              &0                  & 0                 & f^{N}_{1}         & f^{N}_{2}
    \end{bmatrix}.
\end{equation*}
\vspace{0.3cm}

Note that if the system includes a continent, the parameters $\gamma$ and $a(T)$ in (\ref{eqn: f^i}) must be altered for every $i$ such that $\theta_i \in [\theta_{l_1}, \theta_{l_2}]$. We build the matrix $H$ and examine the associated eigenvalues for a large number of  points in the bifurcation diagram. Stability properties are subsequently inferred from the ensemble of stationary solutions within the same branch.

\subsection{Slope-stability Theorem}
Another method for assessing the stability of stationary solutions in EBMs is the so-called "slope-stability theorem" \cite{cahalan1979stability}. The theorem states that the slope of the branch in the bifurcation diagram may be used to determine the stability of the stationary solution $T_0$. For a control parameter $Q$, the theorem may be stated as:
\begin{align*}
    \frac{d Q}{d T_0} &> 0 \hspace{0.5cm}\implies\hspace{0.5cm}T_0\tx{ is stable}\\
    \frac{d Q}{d T_0} &< 0 \hspace{0.5cm}\implies\hspace{0.5cm}T_0\tx{ is unstable}
\end{align*}
A proof of this theorem is omitted, interested readers are directed to \cite{shen1999simple} and \cite{north1990multiple}. In the bifurcation diagrams presented in \ref{section: results} this amounts to: Positive sloped branches are stable, negative sloped branches are unstable. 

\subsection{Heuristic}
Finally, a heuristic stability analysis may be applied to further solidify conclusions regarding the stability properties of a stationary solution. In the heuristic approach, the finite difference algorithm is employed to dynamically run a perturbed stationary solution. The output of the simulation should retain the form of the initial condition for a sufficiently long time in order for the stationary solution to be deemed stable. For unstable stationary solution the dynamical solution will tend to another stationary solution as $t$ grows. However, this approach is expensive. Furthermore, it tends to be inconclusive if there are close-laying points in the bifurcation diagram sharing the same $Q$-value. 

\clearpage
\section{Results}
\label{section: results}

In this section, we present the bifurcation diagram for all the systems analyzed in the paper: the aquaplanet and the planet with a zonally symmetric continent (\ref{eqn: continent configuration}) with the three configurations $\varepsilon = 0$, $\varepsilon = 0.1$ and $\varepsilon = 0.5$. The control parameter is chosen to be the scaled solar constant $Q$. Changes to the system's dynamics and bifurcation phenomena in response to the gradual disruption of meridional symmetry are of interest, and are subsequently discussed.

\begin{figure}[H]
    \centering
    \hspace*{-1cm}
    \includegraphics[scale=0.40]{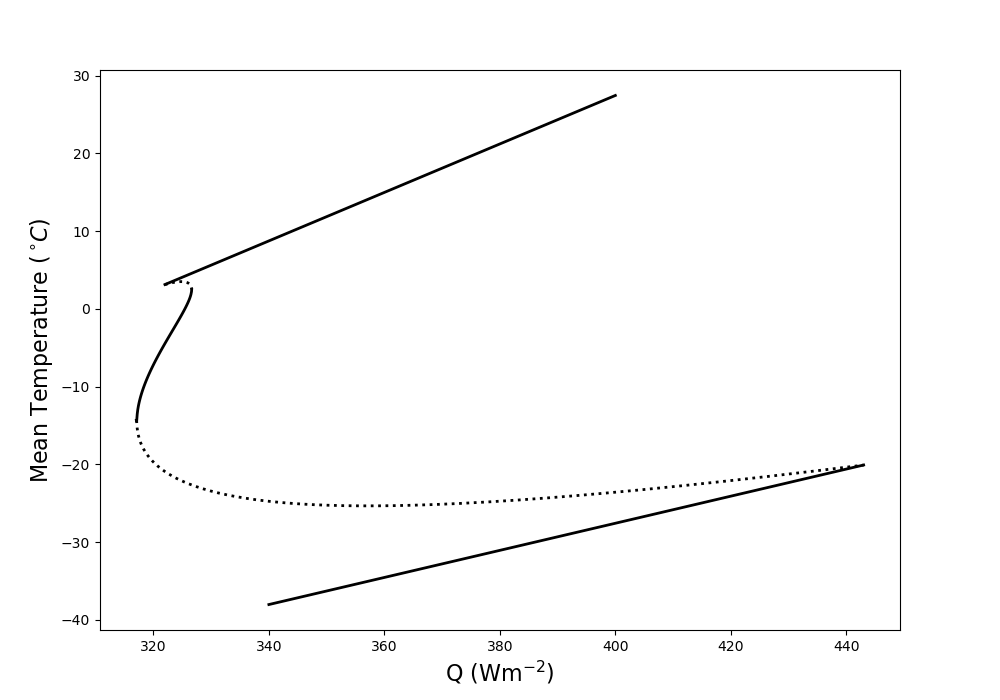}
    \caption{Annual mean equilibrium surface temperatures plotted against the control parameter $Q$ for an aquaplanet. Model parameters are those in Appendix \ref{section: Parameter values}. Solid lines indicate stable solutions and dotted lines are unstable solutions. }
    \label{fig: bifurcation diagram aquaplanet}
\end{figure}

\subsection{Bifurcation Diagrams}
\begin{figure}[H]
    \centering
    \hspace*{-1cm}
    \includegraphics[scale=0.40]{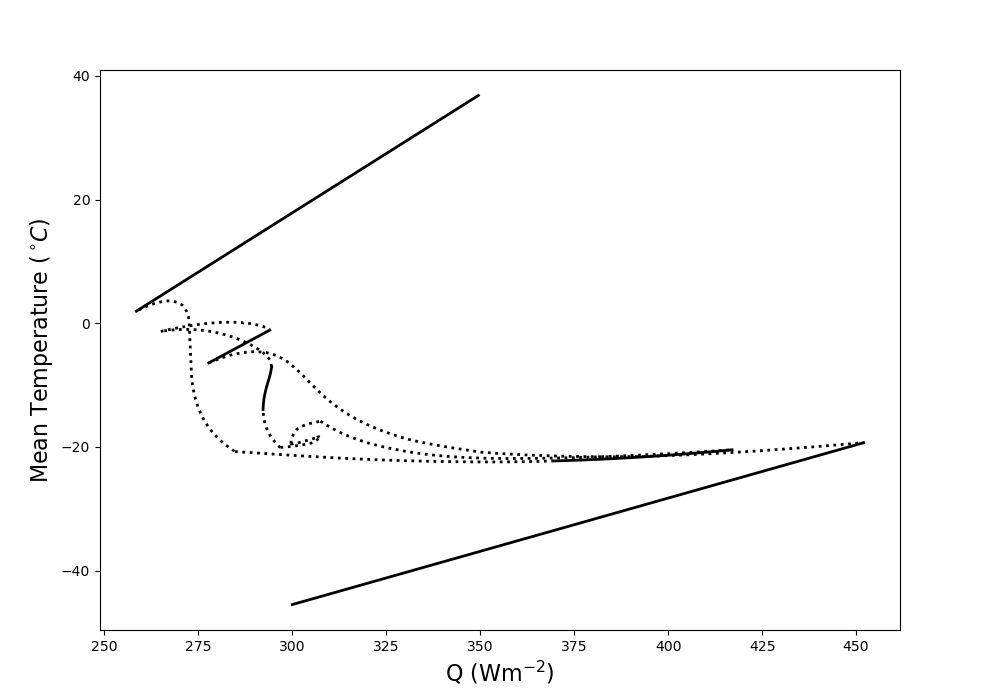}
    \caption{Annual mean equilibrium surface temperatures plotted against the control parameter $Q$ for a planet with a continent configuration for the kind (\ref{eqn: continent configuration}) where $\varepsilon = 0$. Model parameters are those in Appendix \ref{section: Parameter values}. Solid lines indicate stable solutions and dotted lines are unstable solutions. }
    \label{fig: bifurcation diagram eps=0}
\end{figure}

\begin{figure}[H]
    \centering
    \hspace*{-1cm}
    \includegraphics[scale=0.40]{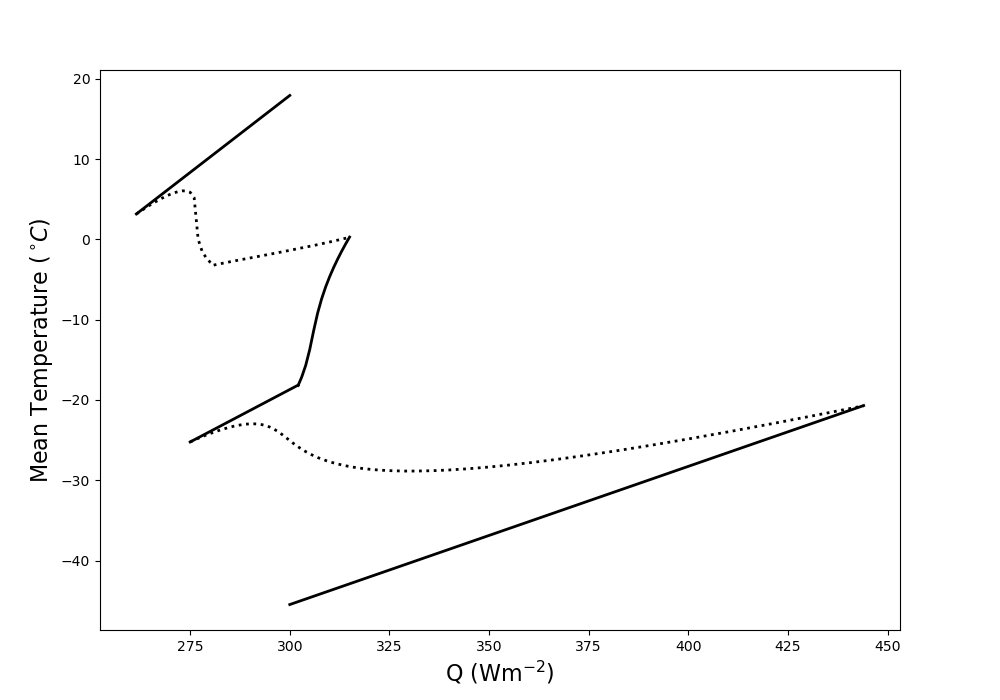}
    \caption{Annual mean equilibrium surface temperatures plotted against the control parameter $Q$ for a planet with a continent configuration for the kind (\ref{eqn: continent configuration}) where $\varepsilon = 0.1$. Model parameters in Appendix \ref{section: Parameter values}. Solid lines indicate stable solutions and dotted lines are unstable solutions. }
    \label{fig: bifurcation diagram eps=0.1}
\end{figure}

\begin{figure}[H]
    \centering
    \hspace*{-1cm}
    \includegraphics[scale=0.40]{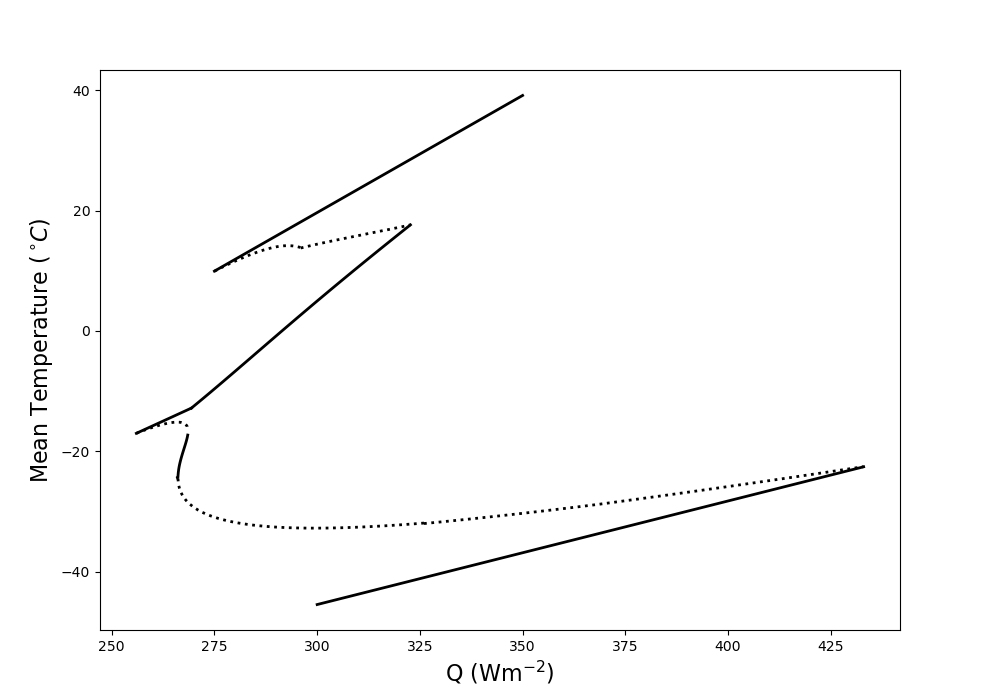}
    \caption{Annual mean equilibrium surface temperatures plotted against the control parameter $Q$ for a planet with a continent configuration for the kind (\ref{eqn: continent configuration}) where $\varepsilon = 0.5$. Model parameters in Appendix \ref{section: Parameter values}. Solid lines indicate stable solutions and dotted lines are unstable solutions. }
    \label{fig: bifurcation diagram eps=0.5}
\end{figure}

Starting with the aquaplanet, this system is well studied, but its bifurcation diagram is included for completeness. Figure \ref{fig: bifurcation diagram aquaplanet} shows the annual mean surface temperature plotted against the control parameter, the scaled solar constant. This system famously has three stable equilibria: the two extreme no partial ice cover solutions, as well as an intermediate state with a partial ice cover. Given the query of this paper, it is interesting investigation is how the bifurcation curve is affected by the introduction of a continent. Starting with the meridionally symmetric continent configuration, Figure \ref{fig: bifurcation diagram eps=0} shows the annual mean surface temperature plotted the scaled solar constant. The introduction of a continent has evidently introduced new equilibrium state. This feature is also observed in the system with the continent configuration (\ref{eqn: continent configuration}) where $\varepsilon = 0.1$ and $\varepsilon = 0.5$, shown the bifurcation diagrams in Figure \ref{fig: bifurcation diagram eps=0.1} and Figure \ref{fig: bifurcation diagram eps=0.5}, respectively. These diagrams are discussed below.

\subsection{Discussion}
The above bifurcation diagrams showcase the sensitivity of the system to the placement of a zonally symmetric continent. For the meridionally symmetric continent configuration, the system exhibits a very fine dynamic where the constraints of the system may be satisfied by a number of ice distributions. The result is an interesting bifurcation curve (see Figure \ref{fig: bifurcation diagram eps=0}) with multiple equilibria in close proximity in the phase space. Branches intersect at multiple points in the bifurcation diagram, forming several "loop" structures. EBMs are highly sensitive to model parameters \cite{soldatenko2019climate} and the extent of these loops in the phase space is affected by the model parameters \cite{samuelsberg2022effects}. Furthermore, for the model parameters used here,  the meridionally symmetric system reveals the presence of up to seven equilibria within a narrow parameter range of the control parameter. With a higher resolution in the bifurcation diagram, there is a possibility that yet more solutions may exist within this regime.

Evidently, the extremely fine dynamic discussed above disappears as the meridional symmetry is broken. Meridional symmetry introduces complexity to the bifurcation curve, giving rise to a dynamic distinct from that observed in asymmetrical cases. For the model parameters used, the parameter regime where five equilibria exist is largest for the symmetrical case. However, the range in which three stable equilibria exists is more or less the same in all the cases analyzed. The meridionally symmetric system exhibits an increased number of bifurcation points, and the stable branches are prominently disjointed. Also noteworthy is the observation that the northward shift of the continent induces a convergence of the associated bifurcation diagram towards that of an aquaplanet.

\clearpage
\section{Conclusion}
\label{section: conclusion}
In this paper, a Boundary Integral Method (BIM) was applied to find stationary solutions to a North-type EBM (\ref{eqn: EBM on symmetric water sphere}). The time dependent equation was solved using finite differences. This algorithm was subsequently applied along with a numerical perturbation scheme to infer stability properties of the stationary solutions. Bifurcation diagrams were presented for the systems analysed: the idealized aquaplanet and a planet with a zonally symmetric continent in three different configurations. These continent configurations were chosen specifically to investigate how the system is affected by symmetry. The first configuration was such that meridional symmetry was preserved. This symmetry was subsequently violated by gradually shifting the continent northward.


The BIM offers an analytical method for solving North-type EBMs. Solutions are expressed through explicit expressions, readily obtainable from quadrature methods. The presented method has some notable advantages compared to other analytical methods \cite{north1975analytical} for solving energy balance equations of this kind: It does not rely on truncating series expansions. Furthermore, the method is computationally speaking very fast and remains straightforward even for problems with partial land and sea geographies and several critical latitudes $\theta_{c_i}$. North and Mengel has discussed the application of Fourier-Legendre series in EBMs with parameter discontinuities at the continent edges \cite{mengel1988seasonal}\cite{north2017energy}: A discontinuity in the albedo and heat capacity parameter $C$ causes a solution expressed through Legendre modes to converge slowly. Moreover, North discussed the potential for changing parameter values on finite zonal strips using spectral methods \cite{north1975analytical}, but failed to address how this may result in several ice edges and subsequent ice albedo feedback dynamics. The dynamics of EBMs are highly sensitive to model parameters \cite{soldatenko2019climate}, the constraints of the model studied may be satisfied by a non-trivial ice distribution once a continent is introduced (See Figure \ref{fig: eps=0, 6 ice edges}).

Although a step function albedo is used here, the method is easily extended to any latitude dependence for the albedo. However, this method, like the spectral method, is limited to a step function-like temperature response at the ice edge. An arbitrary temperature dependent albedo on either side of the ice edge renders the energy balance equation unsolvable through the presented method. An additional, obvious drawback of the presented method is that a root search is required to find the critical latitudes, $\theta_{c_i}$. For solutions with a low number of critical latitudes, this poses no significant challenge. However, as the number critical latitudes increase, so does the complexity of the root search and the necessity for a good starting point in the iteration.

We apply the presented method to solve a North-type EBM with a zonally symmetric continent in three different configurations, showcasing the sensitivity of the system to the placement of such a continent. Meridional symmetry results in a very fine dynamic and a highly interesting bifurcation curve with numerous complicated attractors. With a higher resolution in the bifurcation diagram, there is a possibility that yet more solutions may exist for a small band of $Q$ values for the north-south symmetrical continent configuration. The complex structure of the bifurcation diagram in Figure \ref{fig: bifurcation diagram eps=0} may serve as a hint towards the complexity one might expect in a hypothetical bifurcation diagram for fully-coupled Earth System Models. A dynamical run sweeping the control parameter through such bifurcation structures is bound to produce the noisy output synonymous with realistic climate models.

\clearpage
\begin{appendices}

\section{Parameter Values}
\label{section: Parameter values}
In this study, the parameters in Figure \ref{fig: parameters} were used. Stationary solutions were found using a step function albedo,
\begin{equation*}
    a(T) = \begin{cases} a_1, & T>-T_s \\ a_2, & T<-T_s \end{cases}.
\end{equation*}
The stability analysis was performed using a smooth function replicating the behavior of the step function,
\begin{equation*}
    a(T) = a_1 + \frac{a_2 - a_1}{2}(1 + \tanh(-\sigma(T + 1))),
\end{equation*}
where the slope parameter $\sigma = 50$.

\begin{figure}[H]
    \centering
    \begin{table}[H]
        \centering
        \begin{tabular}{|l | l|}
            \hline
            \cellcolor[gray]{0.9}\textbf{Parameter} & \cellcolor[gray]{0.9}\textbf{Value} \\
            \hline
            $s(\theta)$ & $s_0 + s_1\cos^2(\theta - \frac{\pi}{2})$\\
            $s_0$ & $0.523$ \\
            $s_1$ & $0.716$ \\
            $A$ & $203 \: \text{W m}^{-2}$ \\
            $B$ & $2.09 \: \text{W m}^{-2}(^\circ\text{C})^{-1}$ \\
            $D$ & $0.208\cdot B$ \\
            $C$ & $4.7 \: B \: t_0 $ \\
            $C_{\tx{land}}$ & $0.16 \: B \: t_0 $ \\
            $t_0$ & $1 \: \text{year}$ \\
            $T_s$ & $10 \: ^\circ\text{C}$ \\
            $T_{s,\tx{land}}$ & $1 \: ^\circ\text{C}$ \\
            $a_1$ & $0.06$ \\
            $a_2$ & $0.6$ \\
            $a_{1,\tx{land}}$ & $0.3$ \\
            $a_{2,\tx{land}}$ & $0.6$ \\
            \hline
        \end{tabular}
    \end{table}
    \vspace{-0.5cm} 
    \caption{Model parameters used in the presented work. $A$, $B$, $C$, $C_{\tx{land}}$ and $T_s$ are from \cite{north1981energy}. $D$ is from \cite{MathematicsAndClimate} and $s(\theta)$ is from \cite{mcgehee2012paleoclimate}.}
    \label{fig: parameters}
\end{figure}

\clearpage

\section{Approximate the Differential Operator}
\label{section: Approximate the differential operator}
In this appendix, a discrete approximation of
\begin{align}
    &\mathcal{L}f(x),
    \nonumber\intertext{for the operator}
    \mathcal{L}(\cdot) &= -\frac{1}{\sin x} \parD{}{x}\bigg(\sin x \parD{}{x}(\cdot)\bigg) 
    \label{eqn: operator}
\end{align}
is derived using finite differences. Assuming that $f$ is smooth, it may be expressed in terms of a power series expansion around a point $x_0$ sufficiently close to $x$, 
\begin{equation}
    f(x) \approx a_1 + a_2(x-x_0) + a_3(x-x_0)^2 .
    \label{eqn: f}
\end{equation}
Applying the operator (\ref{eqn: operator}) to $f$ we get
\begin{equation}
    \mathcal{L}f(x) = - \frac{\cos x}{\sin x}a_2 + \bigg( -2 - 2(x-x_0)\frac{\cos x}{\sin x} \bigg) a_3.
    \label{eqn: Lf}
\end{equation}
Introducing a uniform grid for $x$,
\begin{align}
    x_{i+1} &= x_i + h \hsone\Longrightarrow\hsone h=x_{i+1} - x_i,
    \label{eqn: grid}
    \intertext{let}
    f_i &= f(x_i).
    \label{eqn: notation}
\end{align}
The discretized form of equation (\ref{eqn: Lf}) is
\begin{equation}
    \hat{\mathcal{L}}f_i = - \frac{\cos x_i}{\sin x_i}a_2 + \bigg( -2 - 2h\frac{\cos x_i}{\sin x_i} \bigg) a_3.
    \label{eqn: discret Lf}
\end{equation}
The coefficients $a_2$ and $a_3$ are contingent upon the chosen representation of $f$ and the designated point $x_0$. Specifically, the coefficients' values will vary depending on whether a centered, forward, or backward difference approximation is sought. This appendix includes the derivation of all these approximations, starting with the centered difference approximation. 

\vsone

\subsection{Centered Difference Approximation}
\label{appendix: Centered difference approximation}
In order to develop a centered difference approximation to $\mathcal{L}f$ we assume that $f$ can be expressed in terms of a power series expansion around a point $x_i$, 
\begin{equation*}
    f(x) \approx a_1 + a_2(x-x_i) + a_3(x-x_i)^2 .
\end{equation*}
Evaluating $f$ on the grid (\ref{eqn: grid}) and using the notation (\ref{eqn: notation}), for some arbitrarily small $h$ we must have
\begin{align}
    f_{i-1} &= a_1 + a_2(x_{i-1}-x_i)+ a_3(x_{i-1}-x_i)^2\nonumber\\
            &= a_1 -a_2 h + a_3 h^2 \label{eqn: lin sys 1}\\
    f_{i} &= a_1 + a_2(x_{i}-x_i)+ a_3(x_{i}-x_i)^2\nonumber\\
            &= a_1\\
    f_{i+1} &= a_1 + a_2(x_{i+1}-x_i)+ a_3(x_{i+1}-x_i)^2\nonumber\\
            &= a_1 +a_2 h + a_3 h^2 \label{eqn: lin sys 3}.
    \intertext{For non-zero $h$, the solution to the linear system (\ref{eqn: lin sys 1})-(\ref{eqn: lin sys 3}) is}
    a_1 &= f_i \label{eqn: a_1 centered}\\
    a_2 &= \frac{f_{i+1}-f_{i-1}}{2h}\\
    a_3 &= \frac{f_{i-1}-2f_i + f_{i+1}}{2h^2} \label{eqn: a_3 centered}.
\end{align}
Applying the coefficients (\ref{eqn: a_1 centered})-(\ref{eqn: a_3 centered}) to (\ref{eqn: discret Lf}) we get
\begin{equation}
    \hat{\mathcal{L}}f_i = -\frac{2(f_{i-1}-2f_i + f_{i+1}) + (f_{i-1}-4f_i + 3f_{i+1})h\cot x_i}{2h^2}.
    \label{eqn: approx Lf centered}
\end{equation}

The formula (\ref{eqn: approx Lf centered}) gives a centered difference approximation of $\mathcal{L}f$ for a function, $f_i$, evaluated on a uniform grid, $x_i$.

\vsone

\subsection{Forward Difference Approximation}
\label{appendix: Forward difference approximation}
In order to develop a forward difference approximation to $\mathcal{L}f$ we assume that $f$ can be expressed in terms of a power series expansion around a point $x_{i-1}$, 
\begin{equation*}
    f(x) \approx a_1 + a_2(x-x_{i-1}) + a_3(x-x_{i-1})^2 .
\end{equation*}
Evaluating $f$ on the grid (\ref{eqn: grid}) and using the notation (\ref{eqn: notation}), for some arbitrarily small $h$ we must have
\begin{align}
    f_{i-1} &= a_1 + a_2(x_{i-1}-x_{i-1})+ a_3(x_{i-1}-x_{i-1})^2\nonumber\\
            &= a_1 \label{eqn: lin sys 1 forward}\\
    f_{i} &= a_1 + a_2(x_{i}-x_{i-1})+ a_3(x_{i}-x_{i-1})^2\nonumber\\
            &= a_1 + a_2 h + a_3 h^2\\
    f_{i+1} &= a_1 + a_2(x_{i+1}-x_{i-1})+ a_3(x_{i+1}-x_{i-1})^2\nonumber\\
            &= a_1 + 2a_2 h + 4 a_3 h^2 \label{eqn: lin sys 3 forward}.
    \intertext{For non-zero $h$, the solution to the linear system (\ref{eqn: lin sys 1 forward})-(\ref{eqn: lin sys 3 forward}) is}
    a_1 &= f_{i-1} \label{eqn: a_1 forward}\\
    a_2 &= -\frac{3f_{i-1}-4f_i + f_{i+1}}{2h}\\
    a_3 &= \frac{f_{i-1}-2f_i + f_{i+1}}{2h^2} \label{eqn: a_3 forward}.
\end{align}
Applying the coefficients (\ref{eqn: a_1 forward})-(\ref{eqn: a_3 forward}) to (\ref{eqn: discret Lf}) we get
\begin{equation}
    \hat{\mathcal{L}}f_{i-1} = \frac{-2(f_{i-1}-2f_i + f_{i+1}) + (5f_{i-1}-8f_i + 3f_{i+1})h\cot x_i}{2h^2}.
    \label{eqn: approx Lf forward}
\end{equation}

The formula (\ref{eqn: approx Lf forward}) gives a forward difference approximation of $\mathcal{L}f$ for a function, $f_i$, evaluated on a uniform grid, $x_i$.

\vsone

\subsection{Backward Difference Approximation}
\label{appendix: Backward difference approximation}
In order to develop a backward difference approximation to $\mathcal{L}f$ we assume that $f$ can be expressed in terms of a power series expansion around a point $x_{i+1}$, 
\begin{equation*}
    f(x) \approx a_1 + a_2(x-x_{i+1}) + a_3(x-x_{i+1})^2 .
\end{equation*}
Evaluating $f$ on the grid (\ref{eqn: grid}) and using the notation (\ref{eqn: notation}), for some arbitrarily small $h$ we must have
\begin{align}
    f_{i-1} &= a_1 + a_2(x_{i-1}-x_{i+1})+ a_3(x_{i-1}-x_{i+1})^2\nonumber\\
            &= a_1 -2 a_2 h + 4 a_3 h^2 \label{eqn: lin sys 1 backward}\\
    f_{i} &= a_1 + a_2(x_{i}-x_{i+1})+ a_3(x_{i}-x_{i+1})^2\nonumber\\
            &= a_1 - a_2 h + a_3 h^2\\
    f_{i+1} &= a_1 + a_2(x_{i+1}-x_{i+1})+ a_3(x_{i+1}-x_{i+1})^2\nonumber\\
            &= a_1  \label{eqn: lin sys 3 backward}.
    \intertext{For non-zero $h$, the solution to the linear system (\ref{eqn: lin sys 1 backward})-(\ref{eqn: lin sys 3 backward}) is}
    a_1 &= f_{i+1} \label{eqn: a_1 backward}\\
    a_2 &= \frac{f_{i-1}-4f_i + 3f_{i+1}}{2h}\\
    a_3 &= \frac{f_{i-1}-2f_i + f_{i+1}}{2h^2} \label{eqn: a_3 backward}.
\end{align}
Applying the coefficients (\ref{eqn: a_1 backward})-(\ref{eqn: a_3 backward}) to (\ref{eqn: discret Lf}) we get
\begin{equation}
    \hat{\mathcal{L}}f_{i+1} = \frac{-2(f_{i-1}-2f_i + f_{i+1}) + (f_{i-1}- f_{i+1})h\cot x_i}{2h^2}.
    \label{eqn: approx Lf backward}
\end{equation}

The formula (\ref{eqn: approx Lf backward}) gives a backward difference approximation of $\mathcal{L}f$ for a function, $f_i$, evaluated on a uniform grid, $x_i$.

\clearpage

\section{Integral Identity and Green's Function}
\label{Appendix: Fundamental integral identity and a Green's function}
In this appendix, an integral identity and a Green's function for the operator
\begin{equation}
    \mathcal{L}(\cdot) = -\frac{1}{\sin\theta}\parD{}{\theta} \big( \sin\theta \parD{}{\theta} (\cdot) \big)     + \beta (\cdot).
    \label{eqn: operator with beta}
\end{equation}
are derived.

\subsection{Integral Identity}
We seek an integral identity on the form
\begin{equation}
    \int_D dV \: \{u\mathcal{L} v - v\mathcal{L}u  \} = \int_{\partial D} dS \: \{u\mathcal{G}v - v\mathcal{G}u\},
    \label{eqn: integral identity}
\end{equation}
where $u$ and $v$ are some functions defined on the domain, $D$, and $\mathcal{G}$ is some operator that ensures that (\ref{eqn: integral identity}) holds. The stationary solution to (\ref{eqn: EBM on symmetric water sphere}) varies only on the line $\theta \in [0, \pi]$,
\begin{equation*}
    T = T(\theta).
\end{equation*}
For two functions defined on this domain
\begin{align*}
    v&=v(\theta)\\
    u&=u(\theta),
\end{align*}
we must have
\begin{align*}
    \int_{\theta_1}^{\theta_2}dA\: v\mathcal{L}u &= \int_{\theta_1}^{\theta_2}  d\theta\: R^2                                                                             \sin \theta \:v\mathcal{L}u\\
                                                                             &= \int_{\theta_1}^{\theta_2} d\theta\: \sin \theta\: v \left\{-\frac{1}{\sin\theta}\parD{}{\theta} \big( \sin\theta \parD{u}{\theta} \big)     + \beta u \right\}\\
                                                                             &= -\int_{\theta_1}^{\theta_2} d\theta\: v \parD{}{\theta} \big( \sin\theta \parD{u}{\theta} \big) + \beta\int_{\theta_1}^{\theta_2} d\theta\: \sin \theta\: v u  \\
                                                                             &= -v \sin\theta \parD{u}{\theta}\bigg|_{\theta_1}^{\theta_2}  + \int_{\theta_1}^{\theta_2} d\theta\: \parD{v}{\theta}  \sin\theta \parD{u}{\theta} + \beta\int_{\theta_1}^{\theta_2} d\theta\: \sin \theta\: v u  \\
                                                                             &= -v \sin\theta \parD{u}{\theta}\bigg|_{\theta_1}^{\theta_2}  + \parD{v}{\theta}  \sin\theta \:u\bigg|_{\theta_1}^{\theta_2} - \int_{\theta_1}^{\theta_2} d\theta\: u \parD{}{\theta} \big( \sin\theta \parD{v}{\theta} \big)\\
                                                                             & \hspace{3cm}+ \beta\int_{\theta_1}^{\theta_2} d\theta\: \sin \theta\: v u  \\
                                                                             &= \left\{-v \sin\theta \parD{u}{\theta}  + \parD{v}{\theta}  \sin\theta \:u\right\}\bigg|_{\theta_1}^{\theta_2}\\ 
                                                                             & \hspace{3cm}+\int_{\theta_1}^{\theta_2} d\theta\: \sin \theta\: u \left\{-\frac{1}{\sin\theta}\parD{}{\theta} \big( \sin\theta \parD{v}{\theta} \big)     + \beta v \right\}  .\\
\end{align*}
Here the appropriate volume form for a longitude line on the surface of a sphere, $dA= R^2 \sin \theta \: d\theta$, was applied. From the above, we evidently have
\begin{equation*}
    \int_{\theta_1}^{\theta_2} d\theta\: \sin \theta \big\{v\mathcal{L}u - u \mathcal{L}v \big\} = \left\{ u \sin\theta \parD{v}{\theta} -v \sin\theta \parD{u}{\theta}\right\}\bigg|_{\theta_1}^{\theta_2}.
\end{equation*}
This is the desired integral identity for the operator (\ref{eqn: operator with beta}).

\subsection{Find a Green's Function}
\label{section: find a greens function}
A Green's function, $K(\theta, \xi)$, for the operator (\ref{eqn: operator with beta not appendix}), is any solution to the equation (\ref{eqn: equation for greens function not appendix}). The two defining properties of the Dirac-delta function is;
\begin{enumerate}[1)]
    \item for any surface of interest, $S$, we must have
    \begin{equation}
        \int_S dA\: \delta_{\bm{\xi}} = 1,
        \label{eqn: defining propery for dirac-delta 1}
    \end{equation}
    \item for any function, $f(\bm{x})$, defined on $S$ we must have
        \begin{equation}
        \int_S dA\: \delta_{\bm{\xi}} f = f(\bm{\xi}).
        \label{eqn: defining propery for dirac-delta 2}
    \end{equation}
\end{enumerate}
For the line $\theta \in [0, \pi]$ along the surface of a sphere with radius $R$, it can be shown that
\begin{equation}
    \delta_{\xi}(\theta) = \frac{\delta(\theta-\xi)}{2\pi R^2\sin\theta},
    \label{eqn: dirac-delta function}
\end{equation}
where $\delta(\theta-\xi)$ is the usual delayed Dirac-delta function on the line, will ensure that (\ref{eqn: defining propery for dirac-delta 1}) and (\ref{eqn: defining propery for dirac-delta 2}) are satisfied. It is convenient to scale the Green's function we are seeking by a factor such that the right-hand side of (\ref{eqn: dirac-delta function}) becomes unity and (\ref{eqn: equation for greens function not appendix}) becomes
\begin{equation}
    \mathcal{L}K = \frac{\delta(\theta-\xi)}{\sin\theta}.
    \label{eqn: dimensionless equation for greens function}
\end{equation}
We will demand that the Green's function is continuous across $\theta=\xi$, therefore we must have
\begin{equation}
    \lim_{\theta \rightarrow \xi^+} K(\theta, \xi) = \lim_{\theta \rightarrow \xi^-} K(\theta, \xi).
    \label{eqn: condition 1 on K}
\end{equation}
Integrating (\ref{eqn: dimensionless equation for greens function}) over a small interval centered on $\theta=\xi$ we get
\begin{align*}
    \int_{\xi-\varepsilon}^{\xi+\varepsilon} d\theta\:\sin\theta\:\left\{ -\frac{1}{\sin\theta}\parD{}{\theta} \big( \sin\theta \parD{K}{\theta} \big) + \beta K \right\} &= \int_{\xi-\varepsilon}^{\xi+\varepsilon} d\theta\:\sin\theta\:\frac{\delta(\theta-\xi)}{\sin\theta}\\
    -\int_{\xi-\varepsilon}^{\xi+\varepsilon}d\theta\: \parD{}{\theta} \big( \sin\theta \parD{}{\theta} K(\theta, \xi)\big) + \beta \int_{\xi-\varepsilon}^{\xi+\varepsilon} d\theta\:\sin\theta\: K(\theta, \xi) &=1\\
    -\sin\theta \parD{}{\theta} K(\theta, \xi)\bigg|_{\xi-\varepsilon}^{\xi+\varepsilon}  + \beta \int_{\xi-\varepsilon}^{\xi+\varepsilon} d\theta\:\sin\theta\: K(\theta, \xi) &=1.
\end{align*}
Letting $\varepsilon \rightarrow 0$ we must have
\begin{equation}
        \lim_{\theta \rightarrow \xi^+}\parD{}{\theta}K(\theta, \xi) - \lim_{\theta \rightarrow \xi^-}\parD{}{\theta}K(\theta, \xi) = -\frac{1}{\sin\xi}.
    \label{eqn: condition 2 on K}
\end{equation}

At $\theta \neq \xi$ we evidently have
\begin{equation}
    \mathcal{L}K(\theta, \xi) = 0.
    \label{eqn: condition 3 on K}
\end{equation}
We can therefore conclude that $K(\theta, \xi)$ must satisfy the necessary conditions (\ref{eqn: condition 1 on K}), (\ref{eqn: condition 2 on K}) and (\ref{eqn: condition 3 on K}). In order to find a Green's function that solves (\ref{eqn: dimensionless equation for greens function}) we need to find a basis of solutions for an equation on the form 
\begin{equation}
    -\frac{1}{\sin\theta}\parD{}{\theta} \big( \sin\theta \parD{}{\theta}y(\theta) \big) + \beta y(\theta) = 0.
    \label{eqn: General ODE for greens function}
\end{equation}
Introducing a change of variables, $x=\cos\theta$,  and a function, $u$, such that $u(\cos\theta) = y(\theta)$ it can be shown that (\ref{eqn: General ODE for greens function}) can be written on the form
\begin{equation}
        (1-x^2)\frac{\partial^2}{\partial x^2}u(x) + 2x\parD{}{x}u(x) - \beta u(x) =0.
    \label{eqn: ODE for u}
\end{equation}
Let $\lambda$ be a number such that $-\beta = \lambda (\lambda + 1)$. We may now write (\ref{eqn: ODE for u}) as a Legendre equation,
\begin{equation}
    (1-x^2)\frac{\partial^2}{\partial x^2}u(x) + 2x\parD{}{x}u(x) +\lambda (\lambda + 1) u(x) =0,
    \label{eqn: Legendre equation for u}
\end{equation}
which for some arbitrary real or complex value $\lambda$, will have the known basis of solutions $\big\{ P_\lambda,\: Q_\lambda \big\}$. We can therefore use the basis $\big\{ P_\lambda(\cos\theta),\: Q_\lambda(\cos\theta) \big\}$, where $\lambda = \frac{1}{2}\left(\sqrt{1-4\beta}-1 \right)$, to construct the general solution to (\ref{eqn: condition 3 on K}),
\begin{equation}
    K(\theta, \xi) = \begin{cases}
    a(\xi)P_\lambda(\cos\theta) + b(\xi)Q_\lambda(\cos\theta), & \theta > \xi\\[4pt]
    c(\xi)P_\lambda(\cos\theta) + d(\xi)Q_\lambda(\cos\theta), & \theta < \xi
    \end{cases}.
    \label{eqn: Green's function general form}
\end{equation}
The coefficients $a(\xi), b(\xi), c(\xi)$ and $d(\xi)$ can be determined through the conditions (\ref{eqn: condition 1 on K}) and (\ref{eqn: condition 2 on K}). Any choice of these coefficients satisfying (\ref{eqn: condition 1 on K}) and (\ref{eqn: condition 2 on K}) will give a Green's function for the operator (\ref{eqn: operator with beta not appendix}). However, it makes sense for us to seek a Green's function that is non-singular in the domain $\theta \in [0,\pi]$: We want to to develop a set of conditions on the coefficients $a(\xi), b(\xi), c(\xi)$ and $d(\xi)$ to ensure that the Green's function (\ref{eqn: Green's function general form}) is non-singular when $\theta \rightarrow 0$ and $\theta \rightarrow \pi$. Let
\begin{align*}
    &\begin{aligned}
        K_+ (\theta, \xi) &= a(\xi)P_\lambda(\cos\theta) + b(\xi)Q_\lambda(\cos\theta)\\
        K_- (\theta, \xi) &= c(\xi)P_\lambda(\cos\theta) + d(\xi)Q_\lambda(\cos\theta)
    \end{aligned}
    \label{eqn: notation of K_+ and K_- appendix}
    \intertext{such that}
        &K(\theta, \xi) = \begin{cases}
    K_+ (\theta, \xi), & \theta > \xi\\[4pt]
    K_- (\theta, \xi), & \theta < \xi
    \end{cases}.
    \nonumber
\end{align*}
Using a computer algebra system, we find a series expansion of $K_+$ around $\theta = \pi$, and recognize that there is a term in this expansion containing $\log(\pi- \theta)$ with a coefficient $c_{+0}(a(\xi), b(\xi))$. We want to ensure that $K_+$ is non-singular at $\theta = \pi$ and therefore demand that 
\begin{equation}
    c_{+0}(a(\xi), b(\xi))=0 \hsone \forall \: \xi \in [0, \pi].
    \label{eqn: 1st condition on coefficients in the Greens function}
\end{equation}
Similarly, we find a series expansion of $K_-$ around $\theta = 0$. In this expansion there is a term containing $\log(\theta)$ with a coefficient $c_{-0}(c(\xi), d(\xi))$, and we demand that 
\begin{equation}
    c_{-0}(c(\xi), d(\xi))=0 \hsone \forall \: \xi \in [0, \pi].
    \label{eqn: 2nd condition on coefficients in the Greens function}
\end{equation}
Solving the system of equations (\ref{eqn: condition 1 on K}), (\ref{eqn: condition 2 on K}), (\ref{eqn: 1st condition on coefficients in the Greens function}) and (\ref{eqn: 2nd condition on coefficients in the Greens function}) for $a(\xi), b(\xi), c(\xi)$ and $d(\xi)$, we find the following Green's function;
\begin{equation}
           K(\theta, \xi) = \begin{cases}
    \frac{P_\lambda (\cos \xi)(\pi \cot (\pi \lambda)P_\lambda(\cos \theta) - 2 Q_\lambda(\cos \theta) )}{ 2(1+\lambda)(P_\lambda(\cos\xi)Q_{\lambda+1}(\cos\xi) - P_{\lambda+1}(\cos\xi)Q_\lambda(\cos\xi))}        , & \theta > \xi\\[10pt]
    \frac{P_\lambda (\cos \theta)(\pi \cot (\pi \lambda)P_\lambda(\cos \xi) - 2 Q_\lambda(\cos \xi) )}{ 2(1+\lambda)(P_\lambda(\cos\xi)Q_{\lambda+1}(\cos\xi) - P_{\lambda+1}(\cos\xi)Q_\lambda(\cos\xi))}        , & \theta < \xi
    \end{cases}.
    \label{eqn: Green's function appendix}
\end{equation}
The Green's function (\ref{eqn: Green's function appendix}) is non-singular and bounded at $\theta=0$ and $\theta=\pi$. The derivative of the Green's function (\ref{eqn: Green's function appendix}) is also bounded at the boundary and tends to zero $\forall \:\xi \in [0, \pi]$.

\clearpage
\newpage
\clearpage
\section{Boundary Integral Method with a Continent}
\label{appendix: BIEs with a continent}
In this appendix, the BIM is applied to solve the stationary form of the energy balance equation (\ref{eqn: time dependent EBM w continent}) with a zonally symmetric continent, 
\begin{equation*}
    \begin{cases}
    \mathcal{L}T + \beta T = \eta s(\theta)(1 - a_{\text{water}}(T)) - \alpha, &\hsone \theta\in [0, \theta_{l_1}) \cup (\theta_{l_2}, \pi]\\[10pt]
    \mathcal{L}T + \beta T = \eta s(\theta)(1 - a_{\text{land}}(T)) - \alpha, &\hsone  \theta\in [\theta_{l_1}, \theta_{l_2}].
    \end{cases}
\end{equation*}
The first part of this appendix will consist of laying the necessary groundwork for extending the application of the BIM to include a continent. The latter part of this appendix will include the boundary integral equations (BIEs) associated with the various cases, rendered more or less in full for the sake of completeness. 

\vsone

Some conditions must be imposed on the solution on order to include a continent. The solution has to be continuous and must therefore satisfy the following conditions across the boundary between water and land, 
\begin{align*}
    \lim_{\theta \rightarrow \theta_{l_1}^-} T(\theta) &= \lim_{\theta \rightarrow \theta_{l_1}^+} T(\theta) = T(\theta_{l_1})\\
    \lim_{\theta \rightarrow \theta_{l_2}^-} T(\theta) &= \lim_{\theta \rightarrow \theta_{l_2}^+} T(\theta) = T(\theta_{l_2}).
\end{align*}
Furthermore, we can have no build up of heat at the boundary between water and land, nor at the critical latitudes. We must therefore demand that the solution is differentiable across these latitudes,
\begin{align*}
    \lim_{\theta \rightarrow \theta_{l_1}^-} \parD{T}{\theta}(\theta) &= \lim_{\theta \rightarrow \theta_{l_1}^+} \parD{T}{\theta}(\theta) = T(\theta_{l_1})\\
    \lim_{\theta \rightarrow \theta_{l_2}^-} \parD{T}{\theta}(\theta) &= \lim_{\theta \rightarrow \theta_{l_2}^+} \parD{T}{\theta}(\theta) = T(\theta_{l_2})\\
    \lim_{\theta \rightarrow \theta_{c_i}^-} \parD{T}{\theta}(\theta) &= \lim_{\theta \rightarrow \theta_{c_i}^+} \parD{T}{\theta}(\theta) = T(\theta_{c_i}).
\end{align*}
At the poles, the solution must satisfy the boundary conditions 
\begin{align*}
    \lim_{\theta \rightarrow 0} \sin \theta \parD{T}{\theta}(\theta)  &=0
    \intertext{and}
    \lim_{\theta \rightarrow \pi} \sin \theta \parD{T}{\theta}(\theta) &=0.
\end{align*}

\vsone

The albedo value in ice-free conditions is adjusted on the continent. This evidently introduces a dependence on ice conditions in the source term, $h_j$, in the compact energy balance equation
\begin{equation*}
    \mathcal{L} T = h_j.
\end{equation*}
Let
\begin{align}
    h_1 &= \eta s(\theta)(1 - a_1) - \alpha
    \label{eqn: h_1}
    \intertext{be defined for water with no snow/ice cover and let}
    h_2 &= \eta s(\theta)(1 - a_2) - \alpha
    \intertext{be defined for water with snow/ice cover. Similarly, let}
    h_3 &= \eta s(\theta)(1 - a_{1,\tx{land}}) - \alpha
    \intertext{be defined for the continent with no snow/ice cover and let}
    h_4 &= \eta s(\theta)(1 - a_{2,\tx{land}}) - \alpha
    \label{eqn: h_4}
\end{align}
be defined for the continent with snow/ice-cover. The general relation developed in Section \ref{section: Boundary Integral Method},
\begin{equation}
    T(\xi) = \int_{\theta_1}^{\theta_2} d\theta\: \sin \theta \: K(\theta, \xi)h_j(T, \theta) + \left\{ K \sin\theta \parD{T}{\theta} - T \sin\theta \parD{K}{\theta}\right\}\bigg|_{\theta_1}^{\theta_2},
    \label{eqn: general integral relation}
\end{equation}
will be applied in the following. The function $h_j$, where $j=1,2,3,4$, takes the form (\ref{eqn: h_1})-(\ref{eqn: h_4}) depending on the ice conditions on the surface. The application of the BIM depends heavily on the number of critical latitudes, $\theta_{c_i}$, present on the surface, as well as the location of the critical latitudes. In the following sections, cases with an equal number of critical latitudes will be collectively examined, and effects of the critical latitude location on the application of the method will be discussed. The associated BIEs are solved using the numerical technique outlined in Section \ref{section: Boundary integral equation for the case of two ice edges}.

\vsone
\subsection{No Critical Latitudes}
\label{sec: No critical latitudes}
\setcounter{region}{0} 
\setcounter{LowerCaseRoman}{0} 

Consider first the extreme climate states where the planet is either completely covered by ice/snow or completely devoid of any ice/snow cover. The domain is partitioned into region, 
\begin{region}
    \theta \in (0, \theta_{l_1})
    \label{region 1}
\end{region}
\begin{region}
    \theta \in (\theta_{l_1},\theta_{l_2})
    \label{region 2}
\end{region}
\begin{region}
    \theta \in (\theta_{l_2}, \pi).
    \label{region 3}
\end{region}
Focusing first on the ice free planet, the governing equation in region \ref{region 1} is therefore
\begin{align*}
    \mathcal{L} T &= h_1.
    \intertext{In region \ref{region 2} we have}
    \mathcal{L} T &= h_3
    \intertext{and in region \ref{region 3} we have}
    \mathcal{L} T &= h_1.
\end{align*}
The relation (\ref{eqn: general integral relation}) is applied in region \ref{region 1} letting $\theta_1 \rightarrow 0^+$ and $\theta_2 \rightarrow \theta_{l_1}^-$,
\begin{align*}
    T(\xi)  &= \int_0^{\theta_{l_1}}d\theta \: \sin\theta \: K(\theta, \xi) h_1(\theta)\\
            &+ \lim_{\theta_2 \rightarrow \theta_{l_1}^-} K(\theta_{2}, \xi) \sin(\theta_{2}) \parD{T}{\theta}(\theta_{2})\\
            &- \lim_{\theta_2 \rightarrow \theta_{l_1}^-} T(\theta_{2}) \sin(\theta_{2}) \parD{K}{\theta}(\theta_{2},\xi)\\
            &-\lim_{\theta_1  \rightarrow 0^+}K(\theta_1, \xi) \sin (\theta_1) \parD{T}{\theta}(\theta_1)\\
            &+\lim_{\theta_1  \rightarrow 0^+}T(\theta_1) \sin (\theta_1) \parD{K}{\theta}(\theta_1, \xi).
\end{align*}
Applying the boundary condition, as well as the fact that the Green's function is bounded at the boundary, we get the following relation for the solution in region \ref{region 1};
\begin{equation}
    \begin{aligned}
        T(\xi)  &= \int_0^{\theta_{l_1}}d\theta \: \sin\theta \: K(\theta, \xi) h_1(\theta)
            + K(\theta_{l_1}, \xi) \sin(\theta_{l_1}) \parD{T}{\theta}(\theta_{l_1})
            \\&- T(\theta_{l_1}) \sin(\theta_{l_1}) \lim_{\theta  \rightarrow \theta_{l_1}^-} \parD{K}{\theta}(\theta,\xi).
    \end{aligned}
    \label{eqn: no ice general integral relation applied region 1}
\end{equation}
Apply (\ref{eqn: general integral relation}) in region \ref{region 2} letting $\theta_1 \rightarrow \theta_{l_1}^+$ and $\theta_2 \rightarrow \theta_{l_2}^-$ provides a relation for the solution in region \ref{region 2},
\begin{equation}
    \begin{aligned}
        T(\xi)  &= \int_{\theta_{l_1}}^{\theta_{l_2}}d\theta \: \sin\theta \: K(\theta, \xi) h_3(\theta)
            + K(\theta_{l_2}, \xi) \sin(\theta_{l_2}) \parD{T}{\theta}(\theta_{l_2})\\
            &- T(\theta_{l_2}) \sin(\theta_{l_2}) \lim_{\theta  \rightarrow \theta_{l_2}^-}\parD{K}{\theta}(\theta,\xi)
            -K(\theta_{l_1}, \xi) \sin (\theta_{l_1}) \parD{T}{\theta}(\theta_{l_1})\\
            &+T(\theta_{l_1}) \sin (\theta_{l_1}) \lim_{\theta  \rightarrow \theta_{l_1}^+}\parD{K}{\theta}(\theta, \xi).
    \end{aligned}
    \label{eqn: no ice general integral relation applied region 2}
\end{equation}
Finally, a relation for the solution in region \ref{region 3} is obtained by applying (\ref{eqn: general integral relation}) in region \ref{region 3} letting $\theta_1 \rightarrow \theta_{l_2}^+$ and $\theta_2 \rightarrow \pi^-$,
\begin{equation}
    \begin{aligned}
        T(\xi)  &= \int_{\theta_{l_2}}^{\pi}d\theta \: \sin\theta \: K(\theta, \xi) h_1(\theta)
            -K(\theta_{l_2}, \xi) \sin (\theta_{l_2}) \parD{T}{\theta}(\theta_{l_2})\\
            &+T(\theta_{l_2}) \sin (\theta_{l_2}) \lim_{\theta  \rightarrow \theta_{l_2}^+}\parD{K}{\theta}(\theta, \xi).
    \end{aligned}
    \label{eqn: no ice general integral relation applied region 3}
\end{equation}
Here the boundary conditions were applied, as well as the fact that the Green's function is bounded at the boundary. Now we let $\xi$ approach the boundaries of the sub-domains \ref{region 1}, \ref{region 2} and \ref{region 3} in the same fashion as in section \ref{section: Boundary integral equation for the case of two ice edges}. Starting with (\ref{eqn: no ice general integral relation applied region 1}) we let $\xi \rightarrow 0^+$,
\begin{equation}
    \begin{aligned}
        T(0)  &= \int_0^{\theta_{l_1}}d\theta \: \sin\theta \: K(\theta, 0) h_1(\theta)
            + K(\theta_{l_1}, 0) \sin(\theta_{l_1}) \parD{T}{\theta}(\theta_{l_1})
            \\&- T(\theta_{l_1}) \sin(\theta_{l_1}) \lim_{\xi  \rightarrow 0^+}\lim_{\theta  \rightarrow \theta_{l_1}^-} \parD{K}{\theta}(\theta,\xi),
    \end{aligned}
    \label{eqn: no ice BIE 1}
\end{equation}
and $\xi \rightarrow \theta_{l_1}^-$,
\begin{equation}
    \begin{aligned}
        T(\theta_{l_1})  &= \int_0^{\theta_{l_1}}d\theta \: \sin\theta \: K(\theta, \theta_{l_1}) h_1(\theta)
            + \lim_{\xi  \rightarrow \theta_{l_1}^-} K(\theta_{l_1}, \xi) \sin(\theta_{l_1}) \parD{T}{\theta}(\theta_{l_1})
            \\&- T(\theta_{l_1}) \sin(\theta_{l_1}) \lim_{\xi  \rightarrow \theta_{l_1}^-}\lim_{\theta  \rightarrow \theta_{l_1}^-} \parD{K}{\theta}(\theta,\xi).
    \end{aligned}
\end{equation}
Similarly, in (\ref{eqn: no ice general integral relation applied region 2}) we let $\xi \rightarrow \theta_{l_1}^+$,
\begin{equation}
        \begin{aligned}
        T(\theta_{l_1})  &= \int_{\theta_{l_1}}^{\theta_{l_2}}d\theta \: \sin\theta \: K(\theta, \theta_{l_1}) h_3(\theta)
            + K(\theta_{l_2}, \theta_{l_1}) \sin(\theta_{l_2}) \parD{T}{\theta}(\theta_{l_2})\\
            &- T(\theta_{l_2}) \sin(\theta_{l_2}) \lim_{\xi  \rightarrow \theta_{l_1}^+}\lim_{\theta  \rightarrow \theta_{l_2}^-}\parD{K}{\theta}(\theta,\xi)
            -\lim_{\xi  \rightarrow \theta_{l_1}^+}K(\theta_{l_1}, \xi) \sin (\theta_{l_1}) \parD{T}{\theta}(\theta_{l_1})\\
            &+T(\theta_{l_1}) \sin (\theta_{l_1}) \lim_{\xi  \rightarrow \theta_{l_1}^+} \lim_{\theta \rightarrow \theta_{l_1}^+}\parD{K}{\theta}(\theta, \xi),
    \end{aligned}
    \label{eqn: no ice BIE 3}
\end{equation}
and $\xi \rightarrow \theta_{l_2}^-$,
\begin{equation}
        \begin{aligned}
        T(\theta_{l_2})  &= \int_{\theta_{l_1}}^{\theta_{l_2}}d\theta \: \sin\theta \: K(\theta, \theta_{l_2}) h_3(\theta)
            + \lim_{\xi  \rightarrow \theta_{l_2}^-} K(\theta_{l_2}, \xi) \sin(\theta_{l_2}) \parD{T}{\theta}(\theta_{l_2})\\
            &- T(\theta_{l_2}) \sin(\theta_{l_2}) \lim_{\xi  \rightarrow \theta_{l_2}^-}\lim_{\theta  \rightarrow \theta_{l_2}^-}\parD{K}{\theta}(\theta,\xi)
            -K(\theta_{l_1}, \theta_{l_2}) \sin (\theta_{l_1}) \parD{T}{\theta}(\theta_{l_1})\\
            &+T(\theta_{l_1}) \sin (\theta_{l_1}) \lim_{\xi  \rightarrow \theta_{l_2}^-} \lim_{\theta \rightarrow \theta_{l_1}^+}\parD{K}{\theta}(\theta, \xi).
    \end{aligned}
    \label{eqn: no ice BIE 4}
\end{equation}

Finally, in (\ref{eqn: no ice general integral relation applied region 3}) we let $\xi \rightarrow \theta_{l_2}^+$,
\begin{equation}
    \begin{aligned}
        T(\theta_{l_2})  &= \int_{\theta_{l_2}}^{\pi}d\theta \: \sin\theta \: K(\theta, \theta_{l_2}) h_1(\theta)
            - \lim_{\xi  \rightarrow \theta_{l_2}^+} K(\theta_{l_2}, \xi) \sin (\theta_{l_2}) \parD{T}{\theta}(\theta_{l_2})\\
            &+T(\theta_{l_2}) \sin (\theta_{l_2}) \lim_{\xi  \rightarrow \theta_{l_2}^+} \lim_{\theta  \rightarrow \theta_{l_2}^+}\parD{K}{\theta}(\theta, \xi),
    \end{aligned}
    \label{eqn: no ice BIE 5}
\end{equation}
and $\xi \rightarrow \pi^-$,
\begin{equation}
    \begin{aligned}
        T(\pi)  &= \int_{\theta_{l_2}}^{\pi}d\theta \: \sin\theta \: K(\theta, \pi) h_1(\theta)
            -K(\theta_{l_2}, \pi) \sin (\theta_{l_2}) \parD{T}{\theta}(\theta_{l_2})\\
            &+T(\theta_{l_2}) \sin (\theta_{l_2}) \lim_{\xi  \rightarrow \pi^-} \lim_{\theta  \rightarrow \theta_{l_2}^+}\parD{K}{\theta}(\theta, \xi).
    \end{aligned}
    \label{eqn: no ice BIE 6}
\end{equation}
The system of BIEs (\ref{eqn: no ice BIE 1})-(\ref{eqn: no ice BIE 6}) is solved for $T(0)$, $T(\theta_{l_1})$, $T(\theta_{l_2})$, $\parD{T}{\theta}(\theta_{l_1})$, $\parD{T}{\theta}(\theta_{l_2})$ and $T(\pi)$. The solution inside the regions \ref{region 1}, \ref{region 2} and \ref{region 3} is given by (\ref{eqn: no ice general integral relation applied region 1}), (\ref{eqn: no ice general integral relation applied region 2}) and (\ref{eqn: no ice general integral relation applied region 3}), respectively.

\vsone

The other state is the one where the planet is completely ice/snow-covered. The approach here will be analogous to the above process, changing the function $h_j$ in the relevant regions. Hence, the solution inside region \ref{region 1}, \ref{region 2} and \ref{region 3} is given by (in respective order):
\begin{align}
    &\begin{aligned}
        T(\xi)  &= \int_0^{\theta_{l_1}}d\theta \: \sin\theta \: K(\theta, \xi) h_2(\theta)
            + K(\theta_{l_1}, \xi) \sin(\theta_{l_1}) \parD{T}{\theta}(\theta_{l_1})
            \\&- T(\theta_{l_1}) \sin(\theta_{l_1}) \lim_{\theta  \rightarrow \theta_{l_1}^-} \parD{K}{\theta}(\theta,\xi),
    \end{aligned}
    \label{eqn: only ice general integral relation applied region 1}\\
    &\begin{aligned}
        T(\xi)  &= \int_{\theta_{l_1}}^{\theta_{l_2}}d\theta \: \sin\theta \: K(\theta, \xi) h_4(\theta)
            + K(\theta_{l_2}, \xi) \sin(\theta_{l_2}) \parD{T}{\theta}(\theta_{l_2})\\
            &- T(\theta_{l_2}) \sin(\theta_{l_2}) \lim_{\theta  \rightarrow \theta_{l_2}^-}\parD{K}{\theta}(\theta,\xi)
            -K(\theta_{l_1}, \xi) \sin (\theta_{l_1}) \parD{T}{\theta}(\theta_{l_1})\\
            &+T(\theta_{l_1}) \sin (\theta_{l_1}) \lim_{\theta  \rightarrow \theta_{l_1}^+}\parD{K}{\theta}(\theta, \xi)
    \end{aligned}
    \label{eqn: only ice general integral relation applied region 2}\\
    &\begin{aligned}
        T(\xi)  &= \int_{\theta_{l_2}}^{\pi}d\theta \: \sin\theta \: K(\theta, \xi) h_2(\theta)
            -K(\theta_{l_2}, \xi) \sin (\theta_{l_2}) \parD{T}{\theta}(\theta_{l_2})\\
            &+T(\theta_{l_2}) \sin (\theta_{l_2}) \lim_{\theta  \rightarrow \theta_{l_2}^+}\parD{K}{\theta}(\theta, \xi),
    \end{aligned}
    \label{eqn: only ice general integral relation applied region 3}
\end{align}
The BIEs are:
\begin{align}
    &\begin{aligned}
        T(0)  &= \int_0^{\theta_{l_1}}d\theta \: \sin\theta \: K(\theta, 0) h_2(\theta)
            + K(\theta_{l_1}, 0) \sin(\theta_{l_1}) \parD{T}{\theta}(\theta_{l_1})
            \\&- T(\theta_{l_1}) \sin(\theta_{l_1}) \lim_{\xi  \rightarrow 0^+}\lim_{\theta  \rightarrow \theta_{l_1}^-} \parD{K}{\theta}(\theta,\xi)
    \end{aligned}
    \label{eqn: only ice BIE 1}\\
    &\begin{aligned}
        T(\theta_{l_1})  &= \int_0^{\theta_{l_1}}d\theta \: \sin\theta \: K(\theta, \theta_{l_1}) h_2(\theta)
            + \lim_{\xi  \rightarrow \theta_{l_1}^-} K(\theta_{l_1}, \xi) \sin(\theta_{l_1}) \parD{T}{\theta}(\theta_{l_1})
            \\&- T(\theta_{l_1}) \sin(\theta_{l_1}) \lim_{\xi  \rightarrow \theta_{l_1}^-}\lim_{\theta  \rightarrow \theta_{l_1}^-} \parD{K}{\theta}(\theta,\xi)
    \end{aligned}\\
    &\begin{aligned}
        T(\theta_{l_1})  &= \int_{\theta_{l_1}}^{\theta_{l_2}}d\theta \: \sin\theta \: K(\theta, \theta_{l_1}) h_4(\theta)
            + K(\theta_{l_2}, \theta_{l_1}) \sin(\theta_{l_2}) \parD{T}{\theta}(\theta_{l_2})\\
            &- T(\theta_{l_2}) \sin(\theta_{l_2}) \lim_{\xi  \rightarrow \theta_{l_1}^+}\lim_{\theta  \rightarrow \theta_{l_2}^-}\parD{K}{\theta}(\theta,\xi)
            -\lim_{\xi  \rightarrow \theta_{l_1}^+}K(\theta_{l_1}, \xi) \sin (\theta_{l_1}) \parD{T}{\theta}(\theta_{l_1})\\
            &+T(\theta_{l_1}) \sin (\theta_{l_1}) \lim_{\xi  \rightarrow \theta_{l_1}^+} \lim_{\theta \rightarrow \theta_{l_1}^+}\parD{K}{\theta}(\theta, \xi)
    \end{aligned}\\
    &\begin{aligned}
        T(\theta_{l_2})  &= \int_{\theta_{l_1}}^{\theta_{l_2}}d\theta \: \sin\theta \: K(\theta, \theta_{l_2}) h_4(\theta)
            + \lim_{\xi  \rightarrow \theta_{l_2}^-} K(\theta_{l_2}, \xi) \sin(\theta_{l_2}) \parD{T}{\theta}(\theta_{l_2})\\
            &- T(\theta_{l_2}) \sin(\theta_{l_2}) \lim_{\xi  \rightarrow \theta_{l_2}^-}\lim_{\theta  \rightarrow \theta_{l_2}^-}\parD{K}{\theta}(\theta,\xi)
            -K(\theta_{l_1}, \theta_{l_2}) \sin (\theta_{l_1}) \parD{T}{\theta}(\theta_{l_1})\\
            &+T(\theta_{l_1}) \sin (\theta_{l_1}) \lim_{\xi  \rightarrow \theta_{l_2}^-} \lim_{\theta \rightarrow \theta_{l_1}^+}\parD{K}{\theta}(\theta, \xi)
    \end{aligned}
\end{align}
\begin{align}
    &\begin{aligned}
        T(\theta_{l_2})  &= \int_{\theta_{l_2}}^{\pi}d\theta \: \sin\theta \: K(\theta, \theta_{l_2}) h_2(\theta)
            - \lim_{\xi  \rightarrow \theta_{l_2}^+} K(\theta_{l_2}, \xi) \sin (\theta_{l_2}) \parD{T}{\theta}(\theta_{l_2})\\
            &+T(\theta_{l_2}) \sin (\theta_{l_2}) \lim_{\xi  \rightarrow \theta_{l_2}^+} \lim_{\theta  \rightarrow \theta_{l_2}^+}\parD{K}{\theta}(\theta, \xi)
    \end{aligned}\label{eqn: only ice BIE 5}\\
    &\begin{aligned}
        T(\pi)  &= \int_{\theta_{l_2}}^{\pi}d\theta \: \sin\theta \: K(\theta, \pi) h_2(\theta)
            -K(\theta_{l_2}, \pi) \sin (\theta_{l_2}) \parD{T}{\theta}(\theta_{l_2})\\
            &+T(\theta_{l_2}) \sin (\theta_{l_2}) \lim_{\xi  \rightarrow \pi^-} \lim_{\theta  \rightarrow \theta_{l_2}^+}\parD{K}{\theta}(\theta, \xi)
    \end{aligned}
    \label{eqn: only ice BIE 6}
\end{align}

\vsone
Another state with no critical latitudes is possible, with the parameter used in this study: a scenario wherein the continent is entirely covered in ice while the adjacent oceans are ice free. In this state, the continental edges are "ice lines" or "ice edges", but (\ref{eqn: critical latitude condition on continent}) is not satisfied and there are therefore no critical latitudes. Consequently, the same equations from above may be applied to find a solution by by changing $h_3 \to h_4$ in equation (\ref{eqn: no ice general integral relation applied region 1}), (\ref{eqn: no ice BIE 3}) and (\ref{eqn: no ice BIE 3}).

\vsone
Finally, a state where both the northern ocean and the continent is ice covered is possible for asymetrical continent configurations. Consequently, the solution whitin the aforementioned regions are (\ref{eqn: only ice general integral relation applied region 1})-(\ref{eqn: only ice general integral relation applied region 3}), but interchanging $h_2 \to h_1$ in (\ref{eqn: only ice general integral relation applied region 3}) since the southern ocean will be ice-free. Similarly, the BIEs are (\ref{eqn: only ice BIE 1})-(\ref{eqn: only ice BIE 6}), but interchanging $h_2 \to h_1$ in (\ref{eqn: only ice BIE 5}) and (\ref{eqn: only ice BIE 6}).

\vsone
\subsection{One Critical Latitude}
\setcounter{region}{0} 
For asymmetrical continent configurations, climate states with one critical latitude exists. The critical latitude may be either on the northern ocean, on the continent or on the southern ocean. Let us first consider the scenario where the critical latitude is on the northern ocean. The domain is partitioned into the four regions,
\begin{region}
    \theta \in (0, \theta_{c_1}),
    \label{region 1}
\end{region}
\begin{region}
    \theta \in (\theta_{c_1},\theta_{l_1}),
    \label{region 2}
\end{region}
\begin{region}
    \theta \in (\theta_{l_1}, \theta_{l_2}),
    \label{region 3}
\end{region}
\begin{region}
    \theta \in (\theta_{l_2},\pi).
    \label{region 4}
\end{region}
Applying the relations (\ref{eqn: general integral relation}) gives the solution within the respective regions: 
\begin{align}
    &\begin{aligned}
        T(\xi)  &= \int_0^{\theta_{c_1}}d\theta \: \sin\theta \: K(\theta, \xi) h_2(\theta)
            + K(\theta_{c_1}, \xi) \sin(\theta_{c_1}) \parD{T}{\theta}(\theta_{c_1})
            \\&+ \sin(\theta_{c_1}) \lim_{\theta  \rightarrow \theta_{c_1}^-} \parD{K}{\theta}(\theta,\xi)
    \end{aligned}
    \label{eqn: 1 critical latitude on northern ocean general integral relation applied region 1}
\end{align}
\begin{align}
    &\begin{aligned}
        T(\xi)  &= \int_{\theta_{c_1}}^{\theta_{l_1}}d\theta \: \sin\theta \: K(\theta, \xi) h_1(\theta)
            + K(\theta_{l_1}, \xi) \sin(\theta_{l_1}) \parD{T}{\theta}(\theta_{l_1})\\
            &- T(\theta_{l_1}) \sin(\theta_{l_1}) \lim_{\theta  \rightarrow \theta_{l_1}^-}\parD{K}{\theta}(\theta,\xi)
            -K(\theta_{c_1}, \xi) \sin (\theta_{c_1}) \parD{T}{\theta}(\theta_{c_1})\\
            &- \sin (\theta_{c_1}) \lim_{\theta  \rightarrow \theta_{c_1}^+}\parD{K}{\theta}(\theta, \xi)
    \end{aligned}
    \label{eqn: 1 critical latitude on northern ocean general integral relation applied region 2}
    \\
    &\begin{aligned}
        T(\xi)  &= \int_{\theta_{l_1}}^{\theta_{l_2}}d\theta \: \sin\theta \: K(\theta, \xi) h_3(\theta)
            + K(\theta_{l_2}, \xi) \sin(\theta_{l_2}) \parD{T}{\theta}(\theta_{l_2})\\
            &- T(\theta_{l_2}) \sin(\theta_{l_2}) \lim_{\theta  \rightarrow \theta_{l_2}^-}\parD{K}{\theta}(\theta,\xi)
            -K(\theta_{l_1}, \xi) \sin (\theta_{l_1}) \parD{T}{\theta}(\theta_{l_1})\\
            &+T(\theta_{l_1}) \sin (\theta_{l_1}) \lim_{\theta  \rightarrow \theta_{l_1}^+}\parD{K}{\theta}(\theta, \xi)
    \end{aligned}
    \label{eqn: 1 critical latitude on northern ocean general integral relation applied region 3}
    \\
    &\begin{aligned}
        T(\xi)  &= \int_{\theta_{l_2}}^{\pi}d\theta \: \sin\theta \: K(\theta, \xi) h_1(\theta)
            -K(\theta_{l_2}, \xi) \sin (\theta_{l_2}) \parD{T}{\theta}(\theta_{l_2})\\
            &+T(\theta_{l_2}) \sin (\theta_{l_2}) \lim_{\theta  \rightarrow \theta_{l_2}^+}\parD{K}{\theta}(\theta, \xi)
    \end{aligned}
    \label{eqn: 1 critical latitude on northern ocean general integral relation applied region 4}
\end{align}
The associated BIEs are:
\begin{align}
    &\begin{aligned}
        T(0)  &= \int_0^{\theta_{c_1}}d\theta \: \sin\theta \: K(\theta, 0) h_2(\theta)
            + K(\theta_{c_1}, 0) \sin(\theta_{c_1}) \parD{T}{\theta}(\theta_{c_1})
            \\&+ \sin(\theta_{c_1}) \lim_{\xi  \rightarrow 0^+}\lim_{\theta  \rightarrow \theta_{c_1}^-} \parD{K}{\theta}(\theta,\xi)
    \end{aligned}
    \label{eqn: 1 critical latitude critical latitude on northern ocean BIE 1}\\
    &\begin{aligned}
        -1  &= \int_0^{\theta_{c_1}}d\theta \: \sin\theta \: K(\theta, \theta_{c_1}) h_2(\theta)
            + \lim_{\xi  \rightarrow \theta_{c_1}^-} K(\theta_{c_1}, \xi) \sin(\theta_{c_1}) \parD{T}{\theta}(\theta_{c_1})
            \\&+ \sin(\theta_{c_1}) \lim_{\xi  \rightarrow \theta_{c_1}^-}\lim_{\theta  \rightarrow \theta_{c_1}^-} \parD{K}{\theta}(\theta,\xi)
    \end{aligned}\\
    &\begin{aligned}
        -1  &= \int_{\theta_{c_1}}^{\theta_{l_1}}d\theta \: \sin\theta \: K(\theta, \theta_{c_1}) h_1(\theta)
            + K(\theta_{l_1}, \theta_{c_1}) \sin(\theta_{l_1}) \parD{T}{\theta}(\theta_{l_1})\\
            &- T(\theta_{l_1}) \sin(\theta_{l_1}) \lim_{\xi  \rightarrow \theta_{c_1}^+} \lim_{\theta  \rightarrow \theta_{l_1}^-}\parD{K}{\theta}(\theta,\xi)
            -\lim_{\xi  \rightarrow \theta_{c_1}^+}K(\theta_{c_1}, \xi) \sin (\theta_{c_1}) \parD{T}{\theta}(\theta_{c_1})\\
            &- \sin (\theta_{c_1})\lim_{\xi  \rightarrow \theta_{c_1}^+} \lim_{\theta  \rightarrow \theta_{c_1}^+}\parD{K}{\theta}(\theta, \xi)
    \end{aligned}\\
    &\begin{aligned}
        T(\theta_{l_1})  &= \int_{\theta_{c_1}}^{\theta_{l_1}}d\theta \: \sin\theta \: K(\theta, \theta_{l_1}) h_1(\theta)
            + \lim_{\xi  \rightarrow \theta_{l_1}^-}K(\theta_{l_1}, \xi) \sin(\theta_{l_1}) \parD{T}{\theta}(\theta_{l_1})\\
            &- T(\theta_{l_1}) \sin(\theta_{l_1}) \lim_{\xi  \rightarrow \theta_{l_1}^-} \lim_{\theta  \rightarrow \theta_{l_1}^-}\parD{K}{\theta}(\theta,\xi)
            -K(\theta_{c_1}, \theta_{l_1}) \sin (\theta_{c_1}) \parD{T}{\theta}(\theta_{c_1})\\
            &- \sin (\theta_{c_1}) \lim_{\xi  \rightarrow \theta_{l_1}^-} \lim_{\theta  \rightarrow \theta_{c_1}^+}\parD{K}{\theta}(\theta, \xi)
    \end{aligned}
\end{align}
\begin{align}
    &\begin{aligned}
        T(\theta_{l_1})  &= \int_{\theta_{l_1}}^{\theta_{l_2}}d\theta \: \sin\theta \: K(\theta, \theta_{l_1}) h_4(\theta)
            + K(\theta_{l_2}, \theta_{l_1}) \sin(\theta_{l_2}) \parD{T}{\theta}(\theta_{l_2})\\
            &- T(\theta_{l_2}) \sin(\theta_{l_2}) \lim_{\xi  \rightarrow \theta_{l_1}^+}\lim_{\theta  \rightarrow \theta_{l_2}^-}\parD{K}{\theta}(\theta,\xi)
            -\lim_{\xi  \rightarrow \theta_{l_1}^+}K(\theta_{l_1}, \xi) \sin (\theta_{l_1}) \parD{T}{\theta}(\theta_{l_1})\\
            &+T(\theta_{l_1}) \sin (\theta_{l_1}) \lim_{\xi  \rightarrow \theta_{l_1}^+} \lim_{\theta \rightarrow \theta_{l_1}^+}\parD{K}{\theta}(\theta, \xi)
    \end{aligned}\\
    &\begin{aligned}
        T(\theta_{l_2})  &= \int_{\theta_{l_1}}^{\theta_{l_2}}d\theta \: \sin\theta \: K(\theta, \theta_{l_2}) h_4(\theta)
            + \lim_{\xi  \rightarrow \theta_{l_2}^-} K(\theta_{l_2}, \xi) \sin(\theta_{l_2}) \parD{T}{\theta}(\theta_{l_2})\\
            &- T(\theta_{l_2}) \sin(\theta_{l_2}) \lim_{\xi  \rightarrow \theta_{l_2}^-}\lim_{\theta  \rightarrow \theta_{l_2}^-}\parD{K}{\theta}(\theta,\xi)
            -K(\theta_{l_1}, \theta_{l_2}) \sin (\theta_{l_1}) \parD{T}{\theta}(\theta_{l_1})\\
            &+T(\theta_{l_1}) \sin (\theta_{l_1}) \lim_{\xi  \rightarrow \theta_{l_2}^-} \lim_{\theta \rightarrow \theta_{l_1}^+}\parD{K}{\theta}(\theta, \xi).
    \end{aligned}\\
    &\begin{aligned}
        T(\theta_{l_2})  &= \int_{\theta_{l_2}}^{\pi}d\theta \: \sin\theta \: K(\theta, \theta_{l_2}) h_2(\theta)
            - \lim_{\xi  \rightarrow \theta_{l_2}^+} K(\theta_{l_2}, \xi) \sin (\theta_{l_2}) \parD{T}{\theta}(\theta_{l_2})\\
            &+T(\theta_{l_2}) \sin (\theta_{l_2}) \lim_{\xi  \rightarrow \theta_{l_2}^+} \lim_{\theta  \rightarrow \theta_{l_2}^+}\parD{K}{\theta}(\theta, \xi)
    \end{aligned}\\
    &\begin{aligned}
        T(\pi)  &= \int_{\theta_{l_2}}^{\pi}d\theta \: \sin\theta \: K(\theta, \pi) h_2(\theta)
            -K(\theta_{l_2}, \pi) \sin (\theta_{l_2}) \parD{T}{\theta}(\theta_{l_2})\\
            &+T(\theta_{l_2}) \sin (\theta_{l_2}) \lim_{\xi  \rightarrow \pi^-} \lim_{\theta  \rightarrow \theta_{l_2}^+}\parD{K}{\theta}(\theta, \xi)
    \end{aligned}
    \label{eqn: 1 critical latitude critical latitude on northern ocean BIE 8}
\end{align}

\vsone
If the critical latitude is on the continent, the domain must be partitioned into the following regions:
\setcounter{region}{0} 
\begin{region}
    \theta \in (0, \theta_{l_1})
    \label{region 1}
\end{region}
\begin{region}
    \theta \in (\theta_{l_1},\theta_{c_1})
    \label{region 2}
\end{region}
\begin{region}
    \theta \in (\theta_{c_1}, \theta_{l_2})
    \label{region 3}
\end{region}
\begin{region}
    \theta \in (\theta_{l_2},\pi)
    \label{region 4}
\end{region}
The relations (\ref{eqn: general integral relation}) is applied in these regions giving the following solution: 
\begin{align}
    &\begin{aligned}
        T(\xi)  &= \int_0^{\theta_{l_1}}d\theta \: \sin\theta \: K(\theta, \xi) h_2(\theta)
            + K(\theta_{l_1}, \xi) \sin(\theta_{l_1}) \parD{T}{\theta}(\theta_{l_1})
            \\&- T(\theta_{l_1}) \sin(\theta_{l_1}) \lim_{\theta  \rightarrow \theta_{l_1}^-} \parD{K}{\theta}(\theta,\xi)
    \end{aligned}
    \label{eqn: 1 critical latitude on continent general integral relation applied region 1}
\end{align}
\begin{align}
    &\begin{aligned}
        T(\xi)  &= \int_{\theta_{l_1}}^{\theta_{c_1}}d\theta \: \sin\theta \: K(\theta, \xi) h_4(\theta)
            + K(\theta_{c_1}, \xi) \sin(\theta_{c_1}) \parD{T}{\theta}(\theta_{c_1})\\
            &+ T_c \sin(\theta_{c_1}) \lim_{\theta  \rightarrow \theta_{c_1}^-}\parD{K}{\theta}(\theta,\xi)
            -K(\theta_{l_1}, \xi) \sin (\theta_{l_1}) \parD{T}{\theta}(\theta_{l_1})\\
            &+T(\theta_{l_1}) \sin (\theta_{l_1}) \lim_{\theta  \rightarrow \theta_{l_1}^+}\parD{K}{\theta}(\theta, \xi)
    \end{aligned}
    \label{eqn: 1 critical latitude on continent general integral relation applied region 2}
    \\
    &\begin{aligned}
                T(\xi)  &= \int_{\theta_{c_1}}^{\theta_{l_2}}d\theta \: \sin\theta \: K(\theta, \xi) h_3(\theta)
            + K(\theta_{l_2}, \xi) \sin(\theta_{l_2}) \parD{T}{\theta}(\theta_{l_2})\\
            &- T(\theta_{l_2}) \sin(\theta_{l_2}) \lim_{\theta  \rightarrow \theta_{l_2}^-}\parD{K}{\theta}(\theta,\xi)
            -K(\theta_{c_1}, \xi) \sin (\theta_{c_1}) \parD{T}{\theta}(\theta_{c_1})\\
            &-T_c \sin (\theta_{c_1}) \lim_{\theta  \rightarrow \theta_{c_1}^+}\parD{K}{\theta}(\theta, \xi)
    \end{aligned}
    \label{eqn: 1 critical latitude on continent general integral relation applied region 3}
    \\
    &\begin{aligned}
            T(\xi)  &= \int_{\theta_{l_2}}^{\pi}d\theta \: \sin\theta \: K(\theta, \xi) h_1(\theta)
            -K(\theta_{l_2}, \xi) \sin (\theta_{l_2}) \parD{T}{\theta}(\theta_{l_2})\\
            &+T(\theta_{l_2}) \sin (\theta_{l_2}) \lim_{\theta  \rightarrow \theta_{l_2}^+}\parD{K}{\theta}(\theta, \xi)
    \end{aligned}
    \label{eqn: 1 critical latitude on continent general integral relation applied region 4}
\end{align}
In addition to the boundary conditions and the aforementioned properties of the Green's functions, the demand that $T(\theta_{c_1})=-T_c$ on the continent was applied above. The associated BIEs are:
\begin{align}
    &\begin{aligned}
        T(0)  &= \int_0^{\theta_{l_1}}d\theta \: \sin\theta \: K(\theta, 0) h_2(\theta)
            + K(\theta_{l_1}, 0) \sin(\theta_{l_1}) \parD{T}{\theta}(\theta_{l_1})
            \\&- T(\theta_{l_1}) \sin(\theta_{l_1}) \lim_{\xi  \rightarrow 0^+}\lim_{\theta  \rightarrow \theta_{l_1}^-} \parD{K}{\theta}(\theta,\xi)
    \end{aligned}
    \label{eqn: 1 critical latitude on continent BIE 1}
        \end{align}
\begin{align}
    &\begin{aligned}
        T(\theta_{l_1})  &= \int_0^{\theta_{l_1}}d\theta \: \sin\theta \: K(\theta, \theta_{l_1}) h_2(\theta)
            + \lim_{\xi  \rightarrow \theta_{l_1}^-} K(\theta_{l_1}, \xi) \sin(\theta_{l_1}) \parD{T}{\theta}(\theta_{l_1})
            \\&- T(\theta_{l_1}) \sin(\theta_{l_1}) \lim_{\xi  \rightarrow \theta_{l_1}^-}\lim_{\theta  \rightarrow \theta_{l_1}^-} \parD{K}{\theta}(\theta,\xi)
    \end{aligned}\\
    &\begin{aligned}
        T(\theta_{l_1})  &= \int_{\theta_{l_1}}^{\theta_{c_1}}d\theta \: \sin\theta \: K(\theta, \theta_{l_1}) h_4(\theta)
            + K(\theta_{c_1}, \theta_{l_1}) \sin(\theta_{c_1}) \parD{T}{\theta}(\theta_{c_1})\\
            &+ T_c \sin(\theta_{c_1}) \lim_{\xi  \rightarrow \theta_{l_1}^+}\lim_{\theta  \rightarrow \theta_{l_2}^-}\parD{K}{\theta}(\theta,\xi)
            -\lim_{\xi  \rightarrow \theta_{l_1}^+}K(\theta_{l_1}, \xi) \sin (\theta_{l_1}) \parD{T}{\theta}(\theta_{l_1})\\
            &+T(\theta_{l_1}) \sin (\theta_{l_1}) \lim_{\xi  \rightarrow \theta_{l_1}^+} \lim_{\theta \rightarrow \theta_{l_1}^+}\parD{K}{\theta}(\theta, \xi)
    \end{aligned}\\
    &\begin{aligned}
        -T_c  &= \int_{\theta_{l_1}}^{\theta_{c_1}}d\theta \: \sin\theta \: K(\theta, \theta_{c_1}) h_4(\theta)
            + \lim_{\xi  \rightarrow \theta_{c_1}^-} K(\theta_{c_1}, \xi) \sin(\theta_{c_1}) \parD{T}{\theta}(\theta_{c_1})\\
            &+ T_c \sin(\theta_{c_1}) \lim_{\xi  \rightarrow \theta_{c_1}^-} \lim_{\theta  \rightarrow \theta_{c_1}^-}\parD{K}{\theta}(\theta,\xi)
            -K(\theta_{l_1}, \theta_{c_1}) \sin (\theta_{l_1}) \parD{T}{\theta}(\theta_{l_1})\\
            &+T(\theta_{l_1}) \sin (\theta_{l_1}) \lim_{\xi  \rightarrow \theta_{c_1}^-} \lim_{\theta  \rightarrow \theta_{l_1}^+}\parD{K}{\theta}(\theta, \xi)
    \end{aligned}
\end{align}
\begin{align}
    &\begin{aligned}
        -T_c  &= \int_{\theta_{c_1}}^{\theta_{l_2}}d\theta \: \sin\theta \: K(\theta, \theta_{c_1}) h_3(\theta)
            + K(\theta_{l_2}, \theta_{c_1}) \sin(\theta_{l_2}) \parD{T}{\theta}(\theta_{l_2})\\
            &- T(\theta_{l_2}) \sin(\theta_{l_2})\lim_{\xi  \rightarrow \theta_{c_1}^+} \lim_{\theta  \rightarrow \theta_{l_2}^-}\parD{K}{\theta}(\theta,\xi)
            -\lim_{\xi  \rightarrow \theta_{c_1}^+}K(\theta_{c_1}, \xi) \sin (\theta_{c_1}) \parD{T}{\theta}(\theta_{c_1})\\
            &-T_c \sin (\theta_{c_1})\lim_{\xi  \rightarrow \theta_{c_1}^+} \lim_{\theta  \rightarrow \theta_{c_1}^+}\parD{K}{\theta}(\theta, \xi)
    \end{aligned}\\
    &\begin{aligned}
        T(\theta_{l_2})  &= \int_{\theta_{c_1}}^{\theta_{l_2}}d\theta \: \sin\theta \: K(\theta, \theta_{l_2}) h_3(\theta)
            + \lim_{\xi  \rightarrow \theta_{l_2}^-}K(\theta_{l_2}, \xi) \sin(\theta_{l_2}) \parD{T}{\theta}(\theta_{l_2})\\
            &- T(\theta_{l_2}) \sin(\theta_{l_2}) \lim_{\xi  \rightarrow \theta_{l_2}^-} \lim_{\theta  \rightarrow \theta_{l_2}^-}\parD{K}{\theta}(\theta,\xi)
            -K(\theta_{c_1}, \theta_{l_2}) \sin (\theta_{c_1}) \parD{T}{\theta}(\theta_{c_1})\\
            &-T_c \sin (\theta_{c_1}) \lim_{\xi  \rightarrow \theta_{l_2}^-} \lim_{\theta  \rightarrow \theta_{c_1}^+}\parD{K}{\theta}(\theta, \xi)
    \end{aligned}\\
    &\begin{aligned}
        T(\theta_{l_2})  &= \int_{\theta_{l_2}}^{\pi}d\theta \: \sin\theta \: K(\theta, \theta_{l_2}) h_2(\theta)
            - \lim_{\xi  \rightarrow \theta_{l_2}^+} K(\theta_{l_2}, \xi) \sin (\theta_{l_2}) \parD{T}{\theta}(\theta_{l_2})\\
            &+T(\theta_{l_2}) \sin (\theta_{l_2}) \lim_{\xi  \rightarrow \theta_{l_2}^+} \lim_{\theta  \rightarrow \theta_{l_2}^+}\parD{K}{\theta}(\theta, \xi),
    \end{aligned}\\
    &\begin{aligned}
        T(\pi)  &= \int_{\theta_{l_2}}^{\pi}d\theta \: \sin\theta \: K(\theta, \pi) h_2(\theta)
            -K(\theta_{l_2}, \pi) \sin (\theta_{l_2}) \parD{T}{\theta}(\theta_{l_2})\\
            &+T(\theta_{l_2}) \sin (\theta_{l_2}) \lim_{\xi  \rightarrow \pi^-} \lim_{\theta  \rightarrow \theta_{l_2}^+}\parD{K}{\theta}(\theta, \xi)
    \end{aligned}
    \label{eqn: 1 critical latitude on continent BIE 8}
\end{align}

\vsone
Finally, if the critical latitude is on the southern ocean, the domain must be partitioned accordingly. This state is only realized if both the continent and the northern ocean is completely ice covered. However, the continental edge, $\theta_{l_2}$, is not considered a critical latitude since the condition (\ref{eqn: critical latitude condition on continent}) does not hold and this state therefore has only one critical latitude. Consequently, the domain is to be partitioned as follows:
\setcounter{region}{0} 
\begin{region}
    \theta \in (0, \theta_{l_1})
    \label{region 1}
\end{region}
\begin{region}
    \theta \in (\theta_{l_1},\theta_{l_2})
    \label{region 2}
\end{region}
\begin{region}
    \theta \in (\theta_{l_2}, \theta_{c_1})
    \label{region 3}
\end{region}
\begin{region}
    \theta \in (\theta_{c_1},\pi)
    \label{region 4}
\end{region}
The solution within the regions (\ref{region 1})-(\ref{region 4}) is given by:
\begin{align}
    &\begin{aligned}
        T(\xi)  &= \int_0^{\theta_{l_1}}d\theta \: \sin\theta \: K(\theta, \xi) h_2(\theta)
            + K(\theta_{l_1}, \xi) \sin(\theta_{l_1}) \parD{T}{\theta}(\theta_{l_1})
            \\&- T(\theta_{l_1}) \sin(\theta_{l_1}) \lim_{\theta  \rightarrow \theta_{l_1}^-} \parD{K}{\theta}(\theta,\xi),
    \end{aligned}
\end{align}
\begin{align}
    &\begin{aligned}
        T(\xi)  &= \int_{\theta_{l_1}}^{\theta_{l_2}}d\theta \: \sin\theta \: K(\theta, \xi) h_4(\theta)
            + K(\theta_{l_2}, \xi) \sin(\theta_{l_2}) \parD{T}{\theta}(\theta_{l_2})\\
            &- T(\theta_{l_2}) \sin(\theta_{l_2}) \lim_{\theta  \rightarrow \theta_{l_2}^-}\parD{K}{\theta}(\theta,\xi)
            -K(\theta_{l_1}, \xi) \sin (\theta_{l_1}) \parD{T}{\theta}(\theta_{l_1})\\
            &+T(\theta_{l_1}) \sin (\theta_{l_1}) \lim_{\theta  \rightarrow \theta_{l_1}^+}\parD{K}{\theta}(\theta, \xi)
    \end{aligned}\\
    &\begin{aligned}
        T(\xi)  &= \int_{\theta_{l_2}}^{\theta_{c_1}}d\theta \: \sin\theta \: K(\theta, \xi) h_1(\theta)
            + K(\theta_{c_1}, \xi) \sin(\theta_{c_1}) \parD{T}{\theta}(\theta_{c_1})\\
            &+ \sin(\theta_{c_1}) \lim_{\theta  \rightarrow \theta_{c_1}^-}\parD{K}{\theta}(\theta,\xi)
            -K(\theta_{l_2}, \xi) \sin (\theta_{l_1}) \parD{T}{\theta}(\theta_{l_2})\\
            &+T(\theta_{l_2}) \sin (\theta_{l_2}) \lim_{\theta  \rightarrow \theta_{l_2}^+}\parD{K}{\theta}(\theta, \xi)
    \end{aligned}
    \\
    &\begin{aligned}
        T(\xi)  &= \int_{\theta_{c_1}}^{\pi}d\theta \: \sin\theta \: K(\theta, \xi) h_2(\theta)
            -K(\theta_{c_1}, \xi) \sin (\theta_{c_1}) \parD{T}{\theta}(\theta_{c_1})\\
            &- \sin (\theta_{c_1}) \lim_{\theta  \rightarrow \theta_{c_1}^+}\parD{K}{\theta}(\theta, \xi).
    \end{aligned}
\end{align}
The BIEs are:
\begin{align}
    &\begin{aligned}
        T(0)  &= \int_0^{\theta_{l_1}}d\theta \: \sin\theta \: K(\theta, 0) h_2(\theta)
            + K(\theta_{l_1}, 0) \sin(\theta_{l_1}) \parD{T}{\theta}(\theta_{l_1})
            \\&- T(\theta_{l_1}) \sin(\theta_{l_1}) \lim_{\xi  \rightarrow 0^+}\lim_{\theta  \rightarrow \theta_{l_1}^-} \parD{K}{\theta}(\theta,\xi)
    \end{aligned}\\
    &\begin{aligned}
        T(\theta_{l_1})  &= \int_0^{\theta_{l_1}}d\theta \: \sin\theta \: K(\theta, \theta_{l_1}) h_2(\theta)
            + \lim_{\xi  \rightarrow \theta_{l_1}^-} K(\theta_{l_1}, \xi) \sin(\theta_{l_1}) \parD{T}{\theta}(\theta_{l_1})
            \\&- T(\theta_{l_1}) \sin(\theta_{l_1}) \lim_{\xi  \rightarrow \theta_{l_1}^-}\lim_{\theta  \rightarrow \theta_{l_1}^-} \parD{K}{\theta}(\theta,\xi)
    \end{aligned}\\
    &\begin{aligned}
        T(\theta_{l_1})  &= \int_{\theta_{l_1}}^{\theta_{l_2}}d\theta \: \sin\theta \: K(\theta, \theta_{l_1}) h_4(\theta)
            + K(\theta_{l_2}, \theta_{l_1}) \sin(\theta_{l_2}) \parD{T}{\theta}(\theta_{l_2})\\
            &- T(\theta_{l_2}) \sin(\theta_{l_2}) \lim_{\xi  \rightarrow \theta_{l_1}^+}\lim_{\theta  \rightarrow \theta_{l_2}^-}\parD{K}{\theta}(\theta,\xi)
            -\lim_{\xi  \rightarrow \theta_{l_1}^+}K(\theta_{l_1}, \xi) \sin (\theta_{l_1}) \parD{T}{\theta}(\theta_{l_1})\\
            &+T(\theta_{l_1}) \sin (\theta_{l_1}) \lim_{\xi  \rightarrow \theta_{l_1}^+} \lim_{\theta \rightarrow \theta_{l_1}^+}\parD{K}{\theta}(\theta, \xi)
    \end{aligned}\\
    &\begin{aligned}
        T(\theta_{l_2})  &= \int_{\theta_{l_1}}^{\theta_{l_2}}d\theta \: \sin\theta \: K(\theta, \theta_{l_2}) h_4(\theta)
            + \lim_{\xi  \rightarrow \theta_{l_2}^-} K(\theta_{l_2}, \xi) \sin(\theta_{l_2}) \parD{T}{\theta}(\theta_{l_2})\\
            &- T(\theta_{l_2}) \sin(\theta_{l_2}) \lim_{\xi  \rightarrow \theta_{l_2}^-}\lim_{\theta  \rightarrow \theta_{l_2}^-}\parD{K}{\theta}(\theta,\xi)
            -K(\theta_{l_1}, \theta_{l_2}) \sin (\theta_{l_1}) \parD{T}{\theta}(\theta_{l_1})\\
            &+T(\theta_{l_1}) \sin (\theta_{l_1}) \lim_{\xi  \rightarrow \theta_{l_2}^-} \lim_{\theta \rightarrow \theta_{l_1}^+}\parD{K}{\theta}(\theta, \xi)
    \end{aligned}
\end{align}
\begin{align}
    &\begin{aligned}
        T(\theta_{l_2})  &= \int_{\theta_{l_2}}^{\theta_{c_1}}d\theta \: \sin\theta \: K(\theta, \theta_{l_2}) h_1(\theta)
            + K(\theta_{c_1}, \theta_{l_2}) \sin(\theta_{c_1}) \parD{T}{\theta}(\theta_{c_1})\\
            &+ \sin(\theta_{c_1}) \lim_{\xi  \rightarrow \theta_{l_2}^+} \lim_{\theta  \rightarrow \theta_{c_1}^-}\parD{K}{\theta}(\theta,\xi)
            -\lim_{\xi  \rightarrow \theta_{l_2}^+} K(\theta_{l_2}, \xi) \sin (\theta_{l_1}) \parD{T}{\theta}(\theta_{l_2})\\
            &+T(\theta_{l_2}) \sin (\theta_{l_2}) \lim_{\xi  \rightarrow \theta_{l_2}^+} \lim_{\theta  \rightarrow \theta_{l_2}^+}\parD{K}{\theta}(\theta, \xi)
    \end{aligned}\\   
    &\begin{aligned}
        -1  &= \int_{\theta_{l_2}}^{\theta_{c_1}}d\theta \: \sin\theta \: K(\theta, \theta_{c_1}) h_1(\theta)
            + \lim_{\xi  \rightarrow \theta_{c_1}^-} K(\theta_{c_1}, \xi) \sin(\theta_{c_1}) \parD{T}{\theta}(\theta_{c_1})\\
            &+ \sin(\theta_{c_1}) \lim_{\xi  \rightarrow \theta_{c_1}^-} \lim_{\theta  \rightarrow \theta_{c_1}^-}\parD{K}{\theta}(\theta,\xi)
            -K(\theta_{l_2}, \theta_{c_1}) \sin (\theta_{l_1}) \parD{T}{\theta}(\theta_{l_2})\\
            &+T(\theta_{l_2}) \sin (\theta_{l_2}) \lim_{\xi  \rightarrow \theta_{c_1}^-} \lim_{\theta  \rightarrow \theta_{l_2}^+}\parD{K}{\theta}(\theta, \xi)
    \end{aligned}\\
    &\begin{aligned}
        -1  &= \int_{\theta_{c_1}}^{\pi}d\theta \: \sin\theta \: K(\theta, \theta_{c_1}) h_2(\theta)
            - \lim_{\xi  \rightarrow \theta_{c_1}^+} K(\theta_{c_1}, \xi) \sin (\theta_{c_1}) \parD{T}{\theta}(\theta_{c_1})\\
            &- \sin (\theta_{c_1}) \lim_{\xi  \rightarrow \theta_{c_1}^+} \lim_{\theta  \rightarrow \theta_{c_1}^+}\parD{K}{\theta}(\theta, \xi)
    \end{aligned}\\
    &\begin{aligned}
        T(\pi)  &= \int_{\theta_{c_1}}^{\pi}d\theta \: \sin\theta \: K(\theta, \pi) h_2(\theta)
            -K(\theta_{c_1}, \pi) \sin (\theta_{c_1}) \parD{T}{\theta}(\theta_{c_1})\\
            &- \sin (\theta_{c_1}) \lim_{\xi  \rightarrow \pi^-} \lim_{\theta  \rightarrow \theta_{c_1}^+}\parD{K}{\theta}(\theta, \xi)
    \end{aligned}
\end{align}
\vsone
\subsection{Two Critical Latitudes}
\setcounter{region}{0} 
Two critical latitudes are realized by an ice cover centered at both poles, thus having two critical latitudes at both the northern and southern ocean. The domain is partitioned into:
\begin{region}
    \theta \in (0, \theta_{c_1}),
    \label{region 1}
\end{region}
\begin{region}
    \theta \in (\theta_{c_1},\theta_{l_1}),
    \label{region 2}
\end{region}
\begin{region}
    \theta \in (\theta_{l_1}, \theta_{l_2}),
    \label{region 3}
\end{region}
\begin{region}
    \theta \in (\theta_{l_2},\theta_{c_2}),
    \label{region 4}
\end{region}
\begin{region}
    \theta \in (\theta_{c_2}, \pi)
    \label{region 5}
\end{region}
The solutions within these region are found by applying (\ref{eqn: general integral relation}):
\begin{align}
    &\begin{aligned}
        T(\xi)  &= \int_0^{\theta_{c_1}}d\theta \: \sin\theta \: K(\theta, \xi) h_2(\theta)
            + K(\theta_{c_1}, \xi) \sin(\theta_{c_1}) \parD{T}{\theta}(\theta_{c_1})
            \\&+ \sin(\theta_{c_1}) \lim_{\theta  \rightarrow \theta_{c_1}^-} \parD{K}{\theta}(\theta,\xi)
    \end{aligned}
    \label{eqn: 2 critical latitudes ice on continent general integral relation applied region 1}
\end{align}
\begin{align}
    &\begin{aligned}
        T(\xi)  &= \int_{\theta_{c_1}}^{\theta_{l_1}}d\theta \: \sin\theta \: K(\theta, \xi) h_1(\theta)
            + K(\theta_{l_1}, \xi) \sin(\theta_{l_1}) \parD{T}{\theta}(\theta_{l_1})\\
            &- T(\theta_{l_1}) \sin(\theta_{l_1}) \lim_{\theta  \rightarrow \theta_{l_1}^-}\parD{K}{\theta}(\theta,\xi)
            -K(\theta_{c_1}, \xi) \sin (\theta_{c_1}) \parD{T}{\theta}(\theta_{c_1})\\
            &- \sin (\theta_{c_1}) \lim_{\theta  \rightarrow \theta_{c_1}^+}\parD{K}{\theta}(\theta, \xi)
    \end{aligned}
    \label{eqn: 2 critical latitudes ice on continent general integral relation applied region 2}
    \\
    &\begin{aligned}
        T(\xi)  &= \int_{\theta_{l_1}}^{\theta_{l_2}}d\theta \: \sin\theta \: K(\theta, \xi) h_3(\theta)
            + K(\theta_{l_2}, \xi) \sin(\theta_{l_2}) \parD{T}{\theta}(\theta_{l_2})\\
            &- T(\theta_{l_2}) \sin(\theta_{l_2}) \lim_{\theta  \rightarrow \theta_{l_2}^-}\parD{K}{\theta}(\theta,\xi)
            -K(\theta_{l_1}, \xi) \sin (\theta_{l_1}) \parD{T}{\theta}(\theta_{l_1})\\
            &+T(\theta_{l_1}) \sin (\theta_{l_1}) \lim_{\theta  \rightarrow \theta_{l_1}^+}\parD{K}{\theta}(\theta, \xi)
    \end{aligned}
    \label{eqn: 2 critical latitudes ice on continent general integral relation applied region 3}
    \\
    &\begin{aligned}
        T(\xi)  &= \int_{\theta_{l_2}}^{\theta_{c_2}}d\theta \: \sin\theta \: K(\theta, \xi) h_1(\theta)
            + K(\theta_{c_2}, \xi) \sin(\theta_{c_2}) \parD{T}{\theta}(\theta_{c_2})\\
            &+ \sin(\theta_{c_2}) \lim_{\theta  \rightarrow \theta_{c_2}^-}\parD{K}{\theta}(\theta,\xi)
            -K(\theta_{l_2}, \xi) \sin (\theta_{l_1}) \parD{T}{\theta}(\theta_{l_2})\\
            &+T(\theta_{l_2}) \sin (\theta_{l_2}) \lim_{\theta  \rightarrow \theta_{l_2}^+}\parD{K}{\theta}(\theta, \xi)
    \end{aligned}
    \label{eqn: 2 critical latitudes ice on continent general integral relation applied region 4}
    \\
    &\begin{aligned}
        T(\xi)  &= \int_{\theta_{c_2}}^{\pi}d\theta \: \sin\theta \: K(\theta, \xi) h_2(\theta)
            -K(\theta_{c_2}, \xi) \sin (\theta_{c_2}) \parD{T}{\theta}(\theta_{c_2})\\
            &- \sin (\theta_{c_2}) \lim_{\theta  \rightarrow \theta_{c_2}^+}\parD{K}{\theta}(\theta, \xi).
    \end{aligned}
    \label{eqn: 2 critical latitudes ice on continent general integral relation applied region 5}
\end{align}
The associated BIEs are: 
\begin{align}
    &\begin{aligned}
        T(0)  &= \int_0^{\theta_{c_1}}d\theta \: \sin\theta \: K(\theta, 0) h_2(\theta)
            + K(\theta_{c_1}, 0) \sin(\theta_{c_1}) \parD{T}{\theta}(\theta_{c_1})
            \\&+ \sin(\theta_{c_1}) \lim_{\xi  \rightarrow 0^+}\lim_{\theta  \rightarrow \theta_{c_1}^-} \parD{K}{\theta}(\theta,\xi)
    \end{aligned}
    \label{eqn: 2 critical latitudes ice on continent BIE 1}\\
    &\begin{aligned}
        -1  &= \int_0^{\theta_{c_1}}d\theta \: \sin\theta \: K(\theta, \theta_{c_1}) h_2(\theta)
            + \lim_{\xi  \rightarrow \theta_{c_1}^-} K(\theta_{c_1}, \xi) \sin(\theta_{c_1}) \parD{T}{\theta}(\theta_{c_1})
            \\&+ \sin(\theta_{c_1}) \lim_{\xi  \rightarrow \theta_{c_1}^-}\lim_{\theta  \rightarrow \theta_{c_1}^-} \parD{K}{\theta}(\theta,\xi)
    \end{aligned}\\
    &\begin{aligned}
        -1  &= \int_{\theta_{c_1}}^{\theta_{l_1}}d\theta \: \sin\theta \: K(\theta, \theta_{c_1}) h_1(\theta)
            + K(\theta_{l_1}, \theta_{c_1}) \sin(\theta_{l_1}) \parD{T}{\theta}(\theta_{l_1})\\
            &- T(\theta_{l_1}) \sin(\theta_{l_1}) \lim_{\xi  \rightarrow \theta_{c_1}^+} \lim_{\theta  \rightarrow \theta_{l_1}^-}\parD{K}{\theta}(\theta,\xi)
            -\lim_{\xi  \rightarrow \theta_{c_1}^+}K(\theta_{c_1}, \xi) \sin (\theta_{c_1}) \parD{T}{\theta}(\theta_{c_1})\\
            &- \sin (\theta_{c_1})\lim_{\xi  \rightarrow \theta_{c_1}^+} \lim_{\theta  \rightarrow \theta_{c_1}^+}\parD{K}{\theta}(\theta, \xi)
    \end{aligned}
\end{align}
\begin{align}
    &\begin{aligned}
        T(\theta_{l_1})  &= \int_{\theta_{c_1}}^{\theta_{l_1}}d\theta \: \sin\theta \: K(\theta, \theta_{l_1}) h_1(\theta)
            + \lim_{\xi  \rightarrow \theta_{l_1}^-}K(\theta_{l_1}, \xi) \sin(\theta_{l_1}) \parD{T}{\theta}(\theta_{l_1})\\
            &- T(\theta_{l_1}) \sin(\theta_{l_1}) \lim_{\xi  \rightarrow \theta_{l_1}^-} \lim_{\theta  \rightarrow \theta_{l_1}^-}\parD{K}{\theta}(\theta,\xi)
            -K(\theta_{c_1}, \theta_{l_1}) \sin (\theta_{c_1}) \parD{T}{\theta}(\theta_{c_1})\\
            &- \sin (\theta_{c_1}) \lim_{\xi  \rightarrow \theta_{l_1}^-} \lim_{\theta  \rightarrow \theta_{c_1}^+}\parD{K}{\theta}(\theta, \xi)
    \end{aligned}\\
    &\begin{aligned}
        T(\theta_{l_1})  &= \int_{\theta_{l_1}}^{\theta_{l_2}}d\theta \: \sin\theta \: K(\theta, \theta_{l_1}) h_3(\theta)
            + K(\theta_{l_2}, \theta_{l_1}) \sin(\theta_{l_2}) \parD{T}{\theta}(\theta_{l_2})\\
            &- T(\theta_{l_2}) \sin(\theta_{l_2}) \lim_{\xi  \rightarrow \theta_{l_1}^+}\lim_{\theta  \rightarrow \theta_{l_2}^-}\parD{K}{\theta}(\theta,\xi)
            -\lim_{\xi  \rightarrow \theta_{l_1}^+}K(\theta_{l_1}, \xi) \sin (\theta_{l_1}) \parD{T}{\theta}(\theta_{l_1})\\
            &+T(\theta_{l_1}) \sin (\theta_{l_1}) \lim_{\xi  \rightarrow \theta_{l_1}^+} \lim_{\theta \rightarrow \theta_{l_1}^+}\parD{K}{\theta}(\theta, \xi)
            \label{eqn: 2 critical latitudes ice on continent BIE 5}
    \end{aligned}
\end{align}
\begin{align}
    &\begin{aligned}
        T(\theta_{l_2})  &= \int_{\theta_{l_1}}^{\theta_{l_2}}d\theta \: \sin\theta \: K(\theta, \theta_{l_2}) h_3(\theta)
            + \lim_{\xi  \rightarrow \theta_{l_2}^-} K(\theta_{l_2}, \xi) \sin(\theta_{l_2}) \parD{T}{\theta}(\theta_{l_2})\\
            &- T(\theta_{l_2}) \sin(\theta_{l_2}) \lim_{\xi  \rightarrow \theta_{l_2}^-}\lim_{\theta  \rightarrow \theta_{l_2}^-}\parD{K}{\theta}(\theta,\xi)
            -K(\theta_{l_1}, \theta_{l_2}) \sin (\theta_{l_1}) \parD{T}{\theta}(\theta_{l_1})\\
            &+T(\theta_{l_1}) \sin (\theta_{l_1}) \lim_{\xi  \rightarrow \theta_{l_2}^-} \lim_{\theta \rightarrow \theta_{l_1}^+}\parD{K}{\theta}(\theta, \xi)
            \label{eqn: 2 critical latitudes ice on continent BIE 6}
    \end{aligned}\\
    &\begin{aligned}
        T(\theta_{l_2})  &= \int_{\theta_{l_2}}^{\theta_{c_2}}d\theta \: \sin\theta \: K(\theta, \theta_{l_2}) h_1(\theta)
            + K(\theta_{c_2}, \theta_{l_2}) \sin(\theta_{c_2}) \parD{T}{\theta}(\theta_{c_2})\\
            &+ \sin(\theta_{c_2}) \lim_{\xi  \rightarrow \theta_{l_2}^+} \lim_{\theta  \rightarrow \theta_{c_2}^-}\parD{K}{\theta}(\theta,\xi)
            -\lim_{\xi  \rightarrow \theta_{l_2}^+} K(\theta_{l_2}, \xi) \sin (\theta_{l_1}) \parD{T}{\theta}(\theta_{l_2})\\
            &+T(\theta_{l_2}) \sin (\theta_{l_2}) \lim_{\xi  \rightarrow \theta_{l_2}^+} \lim_{\theta  \rightarrow \theta_{l_2}^+}\parD{K}{\theta}(\theta, \xi)
    \end{aligned}   
\end{align}
\begin{align}    
    &\begin{aligned}
        -1  &= \int_{\theta_{l_2}}^{\theta_{c_2}}d\theta \: \sin\theta \: K(\theta, \theta_{c_2}) h_1(\theta)
            + \lim_{\xi  \rightarrow \theta_{c_2}^-} K(\theta_{c_2}, \xi) \sin(\theta_{c_2}) \parD{T}{\theta}(\theta_{c_2})\\
            &+ \sin(\theta_{c_2}) \lim_{\xi  \rightarrow \theta_{c_2}^-} \lim_{\theta  \rightarrow \theta_{c_2}^-}\parD{K}{\theta}(\theta,\xi)
            -K(\theta_{l_2}, \theta_{c_2}) \sin (\theta_{l_1}) \parD{T}{\theta}(\theta_{l_2})\\
            &+T(\theta_{l_2}) \sin (\theta_{l_2}) \lim_{\xi  \rightarrow \theta_{c_2}^-} \lim_{\theta  \rightarrow \theta_{l_2}^+}\parD{K}{\theta}(\theta, \xi)
    \end{aligned}\\
    &\begin{aligned}
        -1  &= \int_{\theta_{c_2}}^{\pi}d\theta \: \sin\theta \: K(\theta, \theta_{c_2}) h_2(\theta)
            - \lim_{\xi  \rightarrow \theta_{c_2}^+} K(\theta_{c_2}, \xi) \sin (\theta_{c_2}) \parD{T}{\theta}(\theta_{c_2})\\
            &- \sin (\theta_{c_2}) \lim_{\xi  \rightarrow \theta_{c_2}^+} \lim_{\theta  \rightarrow \theta_{c_2}^+}\parD{K}{\theta}(\theta, \xi)
    \end{aligned}\\
    &\begin{aligned}
        T(\pi)  &= \int_{\theta_{c_2}}^{\pi}d\theta \: \sin\theta \: K(\theta, \pi) h_2(\theta)
            -K(\theta_{c_2}, \pi) \sin (\theta_{c_2}) \parD{T}{\theta}(\theta_{c_2})\\
            &- \sin (\theta_{c_2}) \lim_{\xi  \rightarrow \pi^-} \lim_{\theta  \rightarrow \theta_{c_2}^+}\parD{K}{\theta}(\theta, \xi)
    \end{aligned}
    \label{eqn: 2 critical latitudes ice on continent BIE 10}
\end{align}

\vsone
States with two critical latitudes on either ocean are also possible with a completely ice covered continent. Solutions are found using the equations above, but the function $h_3$ must be changed to $h_4$ in region \ref{region 3}. Consequently, the solution is given by (\ref{eqn: 2 critical latitudes ice on continent general integral relation applied region 1})-(\ref{eqn: 2 critical latitudes ice on continent general integral relation applied region 5}), but changing $h_3 \to h_4$ in (\ref{eqn: 2 critical latitudes ice on continent general integral relation applied region 3}). BIE are (\ref{eqn: 2 critical latitudes ice on continent BIE 1})-(\ref{eqn: 2 critical latitudes ice on continent BIE 10}), changing $h_3 \to h_4$ in (\ref{eqn: 2 critical latitudes ice on continent BIE 5}) and (\ref{eqn: 2 critical latitudes ice on continent BIE 6}).

\vsone
\setcounter{region}{0} 
For the asymmetric continent configurations studied, there may be states where there are two critical latitudes on the southern ocean. The domain is partitioned into the following sub-domains:
\begin{region}
    \theta \in (0, \theta_{l_1})
    \label{region 1}
\end{region}
\begin{region}
    \theta \in (\theta_{l_1},\theta_{l_2})
    \label{region 2}
\end{region}
\begin{region}
    \theta \in (\theta_{l_2}, \theta_{c_1}),
    \label{region 3}
\end{region}
\begin{region}
    \theta \in (\theta_{c_1},\theta_{c_2}),
    \label{region 4}
\end{region}
\begin{region}
    \theta \in (\theta_{c_2}, \pi)
    \label{region 5}
\end{region}
The solution within these regions are (in order):
\begin{align}
    &\begin{aligned}
        T(\xi)  &= \int_0^{\theta_{l_1}}d\theta \: \sin\theta \: K(\theta, \xi) h_2(\theta)
            + K(\theta_{l_1}, \xi) \sin(\theta_{l_1}) \parD{T}{\theta}(\theta_{l_1})
            \\&- T(\theta_{l_1}) \sin(\theta_{l_1}) \lim_{\theta  \rightarrow \theta_{l_1}^-} \parD{K}{\theta}(\theta,\xi),
    \end{aligned}\\
    &\begin{aligned}
        T(\xi)  &= \int_{\theta_{l_1}}^{\theta_{l_2}}d\theta \: \sin\theta \: K(\theta, \xi) h_4(\theta)
            + K(\theta_{l_2}, \xi) \sin(\theta_{l_2}) \parD{T}{\theta}(\theta_{l_2})\\
            &- T(\theta_{l_2}) \sin(\theta_{l_2}) \lim_{\theta  \rightarrow \theta_{l_2}^-}\parD{K}{\theta}(\theta,\xi)
            -K(\theta_{l_1}, \xi) \sin (\theta_{l_1}) \parD{T}{\theta}(\theta_{l_1})\\
            &+T(\theta_{l_1}) \sin (\theta_{l_1}) \lim_{\theta  \rightarrow \theta_{l_1}^+}\parD{K}{\theta}(\theta, \xi)
    \end{aligned}\\
    &\begin{aligned}
        T(\xi)  &= \int_{\theta_{l_2}}^{\theta_{c_1}}d\theta \: \sin\theta \: K(\theta, \xi) h_2(\theta)
            + K(\theta_{c_1}, \xi) \sin(\theta_{c_1}) \parD{T}{\theta}(\theta_{c_1})\\
            &+ \sin(\theta_{c_1}) \lim_{\theta  \rightarrow \theta_{c_1}^-}\parD{K}{\theta}(\theta,\xi)
            -K(\theta_{l_2}, \xi) \sin (\theta_{l_1}) \parD{T}{\theta}(\theta_{l_2})\\
            &+T(\theta_{l_2}) \sin (\theta_{l_2}) \lim_{\theta  \rightarrow \theta_{l_2}^+}\parD{K}{\theta}(\theta, \xi)
    \end{aligned}\\
    &\begin{aligned}
        T(\xi) &= \int_{\theta_{c_1}}^{\theta_{c_2}} d\theta\: \sin\theta \:K(\theta, \xi) h_1(\theta)
                +  K(\theta_{c_2}, \xi) \sin(\theta_{c_2}) \lim_{\theta_2 \rightarrow \theta_{c_2}^-}\parD{T}{\theta}(\theta_{2})\\
                &+ \sin(\theta_{c_2}) \lim_{\theta_2 \rightarrow \theta_{c_2}^-}\parD{K}{\theta}(\theta_{2},\xi)\\
                &-K(\theta_{c_1}, \xi) \sin (\theta_{c_1}) \lim_{\theta_1  \rightarrow \theta_{c_1}^+}\parD{T}{\theta}(\theta_1)\\
                &- \sin (\theta_{c_1}) \lim_{\theta_1  \rightarrow \theta_{c_1}^+} \parD{K}{\theta}(\theta_1, \xi).
    \end{aligned}
\end{align}
\begin{align}
    &\begin{aligned}
        T(\xi) &= \int_{\theta_{c_2}}^{\pi} d\theta\: \sin\theta \:K(\theta, \xi) h_2(\theta)
                -K(\theta_{c_2}, \xi) \sin (\theta_{c_2}) \lim_{\theta_1  \rightarrow \theta_{c_2}^+}\parD{T}{\theta}(\theta_1)
                \\&- \sin (\theta_{c_2}) \lim_{\theta_1  \rightarrow \theta_{c_2}^+}\parD{K}{\theta}(\theta_1, \xi).
    \end{aligned}
\end{align}
The associated BIEs are:
\begin{align}
    &\begin{aligned}
        T(0)  &= \int_0^{\theta_{l_1}}d\theta \: \sin\theta \: K(\theta, 0) h_2(\theta)
            + K(\theta_{l_1}, 0) \sin(\theta_{l_1}) \parD{T}{\theta}(\theta_{l_1})
            \\&- T(\theta_{l_1}) \sin(\theta_{l_1}) \lim_{\xi  \rightarrow 0^+}\lim_{\theta  \rightarrow \theta_{l_1}^-} \parD{K}{\theta}(\theta,\xi)
    \end{aligned}\\
    &\begin{aligned}
        T(\theta_{l_1})  &= \int_0^{\theta_{l_1}}d\theta \: \sin\theta \: K(\theta, \theta_{l_1}) h_2(\theta)
            + \lim_{\xi  \rightarrow \theta_{l_1}^-} K(\theta_{l_1}, \xi) \sin(\theta_{l_1}) \parD{T}{\theta}(\theta_{l_1})
            \\&- T(\theta_{l_1}) \sin(\theta_{l_1}) \lim_{\xi  \rightarrow \theta_{l_1}^-}\lim_{\theta  \rightarrow \theta_{l_1}^-} \parD{K}{\theta}(\theta,\xi)
    \end{aligned}\\
    &\begin{aligned}
        T(\theta_{l_1})  &= \int_{\theta_{l_1}}^{\theta_{l_2}}d\theta \: \sin\theta \: K(\theta, \theta_{l_1}) h_4(\theta)
            + K(\theta_{l_2}, \theta_{l_1}) \sin(\theta_{l_2}) \parD{T}{\theta}(\theta_{l_2})\\
            &- T(\theta_{l_2}) \sin(\theta_{l_2}) \lim_{\xi  \rightarrow \theta_{l_1}^+}\lim_{\theta  \rightarrow \theta_{l_2}^-}\parD{K}{\theta}(\theta,\xi)
            -\lim_{\xi  \rightarrow \theta_{l_1}^+}K(\theta_{l_1}, \xi) \sin (\theta_{l_1}) \parD{T}{\theta}(\theta_{l_1})\\
            &+T(\theta_{l_1}) \sin (\theta_{l_1}) \lim_{\xi  \rightarrow \theta_{l_1}^+} \lim_{\theta \rightarrow \theta_{l_1}^+}\parD{K}{\theta}(\theta, \xi)
    \end{aligned}\\
    &\begin{aligned}
        T(\theta_{l_2})  &= \int_{\theta_{l_1}}^{\theta_{l_2}}d\theta \: \sin\theta \: K(\theta, \theta_{l_2}) h_4(\theta)
            + \lim_{\xi  \rightarrow \theta_{l_2}^-} K(\theta_{l_2}, \xi) \sin(\theta_{l_2}) \parD{T}{\theta}(\theta_{l_2})\\
            &- T(\theta_{l_2}) \sin(\theta_{l_2}) \lim_{\xi  \rightarrow \theta_{l_2}^-}\lim_{\theta  \rightarrow \theta_{l_2}^-}\parD{K}{\theta}(\theta,\xi)
            -K(\theta_{l_1}, \theta_{l_2}) \sin (\theta_{l_1}) \parD{T}{\theta}(\theta_{l_1})\\
            &+T(\theta_{l_1}) \sin (\theta_{l_1}) \lim_{\xi  \rightarrow \theta_{l_2}^-} \lim_{\theta \rightarrow \theta_{l_1}^+}\parD{K}{\theta}(\theta, \xi)
    \end{aligned}\\
    &\begin{aligned}
        T(\theta_{l_2})  &= \int_{\theta_{l_2}}^{\theta_{c_1}}d\theta \: \sin\theta \: K(\theta, \theta_{l_2}) h_2(\theta)
            + K(\theta_{c_1}, \theta_{l_2}) \sin(\theta_{c_1}) \parD{T}{\theta}(\theta_{c_1})\\
            &+ \sin(\theta_{c_1}) \lim_{\xi  \rightarrow \theta_{l_2}^+} \lim_{\theta  \rightarrow \theta_{c_1}^-}\parD{K}{\theta}(\theta,\xi)
            -\lim_{\xi  \rightarrow \theta_{l_2}^+} K(\theta_{l_2}, \xi) \sin (\theta_{l_1}) \parD{T}{\theta}(\theta_{l_2})\\
            &+T(\theta_{l_2}) \sin (\theta_{l_2}) \lim_{\xi  \rightarrow \theta_{l_2}^+} \lim_{\theta  \rightarrow \theta_{l_2}^+}\parD{K}{\theta}(\theta, \xi)
    \end{aligned}\\    
    &\begin{aligned}
        -1  &= \int_{\theta_{l_2}}^{\theta_{c_1}}d\theta \: \sin\theta \: K(\theta, \theta_{c_1}) h_2(\theta)
            + \lim_{\xi  \rightarrow \theta_{c_1}^-} K(\theta_{c_1}, \xi) \sin(\theta_{c_1}) \parD{T}{\theta}(\theta_{c_1})\\
            &+ \sin(\theta_{c_1}) \lim_{\xi  \rightarrow \theta_{c_1}^-} \lim_{\theta  \rightarrow \theta_{c_1}^-}\parD{K}{\theta}(\theta,\xi)
            -K(\theta_{l_2}, \theta_{c_1}) \sin (\theta_{l_1}) \parD{T}{\theta}(\theta_{l_2})\\
            &+T(\theta_{l_2}) \sin (\theta_{l_2}) \lim_{\xi  \rightarrow \theta_{c_1}^-} \lim_{\theta  \rightarrow \theta_{l_2}^+}\parD{K}{\theta}(\theta, \xi)
    \end{aligned}
\end{align}
\begin{align}
    &\begin{aligned}
        -1 &=\int_{\theta_{c_1}}^{\theta_{c_2}}d\theta\:\sin\theta \:K(\theta, \theta_{c_1}) h_1(\theta)
        + K(\theta_{c_2}, \theta_{c_1})\sin(\theta_{c_2})\parD{T}{\theta}(\theta_{c_2}) 
        \\&+ \sin(\theta_{c_2})\lim_{\xi \rightarrow \theta_{c_1}^+}\lim_{\theta \rightarrow \theta_{c_2}^-}\parD{K}{\theta}(\theta,\xi)
        -K(\theta_{c_1}, \theta_{c_1})\sin(\theta_{c_1})\parD{T}{\theta}(\theta_{c_1})
        \\&-\sin(\theta_{c_1})\lim_{\xi \rightarrow \theta_{c_1}^+}\lim_{\theta \rightarrow \theta_{c_1}^+}\parD{K}{\theta}(\theta, \xi).
    \end{aligned}\\ 
    &\begin{aligned}
        -1 &=\int_{\theta_{c_1}}^{\theta_{c_2}}d\theta\:\sin\theta \:K(\theta, \theta_{c_2}) h_1(\theta)
        + K(\theta_{c_2}, \theta_{c_2})\sin(\theta_{c_2})\parD{T}{\theta}(\theta_{c_2}) 
        \\&+ \sin(\theta_{c_2})\lim_{\xi \rightarrow \theta_{c_2}^-}\lim_{\theta \rightarrow \theta_{c_2}^-}\parD{K}{\theta}(\theta,\xi)
        -K(\theta_{c_1}, \theta_{c_2})\sin(\theta_{c_1})\parD{T}{\theta}(\theta_{c_1})
        \\&-\sin(\theta_{c_1})\lim_{\xi \rightarrow \theta_{c_2}^-}\lim_{\theta \rightarrow \theta_{c_1}^+}\parD{K}{\theta}(\theta, \xi).
    \end{aligned}\\
    &\begin{aligned}
        -1 &= \int_{\theta_{c_2}}^\pi d\theta \: \sin\theta \: K(\theta, \theta_{c_2}) h_2(\theta) 
        -K(\theta_{c_2}, \theta_{c_2}) \sin (\theta_{c_2}) \parD{T}{\theta}(\theta_{c_2})
        \\&- \sin (\theta_{c_2})\lim_{\xi \rightarrow \theta_{c_2}^+}\lim_{\theta \rightarrow \theta_{c_2}^+}\parD{K}{\theta}(\theta, \xi).
\end{aligned}\\
&\begin{aligned}
        T(\pi) &= \int_{\theta_{c_2}}^\pi d\theta \: \sin\theta \: K(\theta, \pi) h_2(\theta) 
        -K(\theta_{c_2}, \pi) \sin (\theta_{c_2}) \parD{T}{\theta}(\theta_{c_2})
        \\&- \sin (\theta_{c_2})\lim_{\xi \rightarrow \pi^-}\lim_{\theta \rightarrow \theta_{c_2}^+}\parD{K}{\theta}(\theta, \xi).
\end{aligned}
\end{align}

\vsone
\setcounter{region}{0} 
For the asymmetrical continent configurations, two critical latitudes are also realized through one critical latitude on the continent and one critical latitude the northern ocean. The domain is partitioned into the following sub-domains:
\begin{region}
    \theta \in (0, \theta_{c_1}),
    \label{region 1}
\end{region}
\begin{region}
    \theta \in (\theta_{c_1},\theta_{l_1}),
    \label{region 2}
\end{region}
\begin{region}
    \theta \in (\theta_{l_1}, \theta_{c_2}),
    \label{region 3}
\end{region}
\begin{region}
    \theta \in (\theta_{c_2},\theta_{l_2}),
    \label{region 4}
\end{region}
\begin{region}
    \theta \in (\theta_{l_2}, \pi).
    \label{region 5}
\end{region}
The solution in each of these regions are given by:
\begin{align}
    &\begin{aligned}
        T(\xi)  &= \int_0^{\theta_{c_1}}d\theta \: \sin\theta \: K(\theta, \xi) h_2(\theta)
            + K(\theta_{c_1}, \xi) \sin(\theta_{c_1}) \parD{T}{\theta}(\theta_{c_1})
            \\&+ \sin(\theta_{c_1}) \lim_{\theta  \rightarrow \theta_{c_1}^-} \parD{K}{\theta}(\theta,\xi)
    \end{aligned}
    \label{eqn: 2 critical latitudes asymmetrical general integral relation applied region 1}
\end{align}
\begin{align}
    &\begin{aligned}
        T(\xi)  &= \int_{\theta_{c_1}}^{\theta_{l_1}}d\theta \: \sin\theta \: K(\theta, \xi) h_1(\theta)
            + K(\theta_{l_1}, \xi) \sin(\theta_{l_1}) \parD{T}{\theta}(\theta_{l_1})\\
            &- T(\theta_{l_1}) \sin(\theta_{l_1}) \lim_{\theta  \rightarrow \theta_{l_1}^-}\parD{K}{\theta}(\theta,\xi)
            -K(\theta_{c_1}, \xi) \sin (\theta_{c_1}) \parD{T}{\theta}(\theta_{c_1})\\
            &- \sin (\theta_{c_1}) \lim_{\theta  \rightarrow \theta_{c_1}^+}\parD{K}{\theta}(\theta, \xi)
    \end{aligned}\\
    &\begin{aligned}
        T(\xi)  &= \int_{\theta_{l_1}}^{\theta_{c_2}}d\theta \: \sin\theta \: K(\theta, \xi) h_4(\theta)
            + K(\theta_{c_2}, \xi) \sin(\theta_{c_2}) \parD{T}{\theta}(\theta_{c_2})\\
            &+ T_c \sin(\theta_{c_2}) \lim_{\theta  \rightarrow \theta_{c_2}^-}\parD{K}{\theta}(\theta,\xi)
            -K(\theta_{l_1}, \xi) \sin (\theta_{l_1}) \parD{T}{\theta}(\theta_{l_1})\\
            &+T(\theta_{l_1}) \sin (\theta_{l_1}) \lim_{\theta  \rightarrow \theta_{l_1}^+}\parD{K}{\theta}(\theta, \xi)
    \end{aligned}\\
    &\begin{aligned}
        T(\xi)  &= \int_{\theta_{c_2}}^{\theta_{l_2}}d\theta \: \sin\theta \: K(\theta, \xi) h_3(\theta)
            + K(\theta_{l_2}, \xi) \sin(\theta_{l_2}) \parD{T}{\theta}(\theta_{l_2})\\
            &- T(\theta_{l_2}) \sin(\theta_{l_2}) \lim_{\theta  \rightarrow \theta_{l_2}^-}\parD{K}{\theta}(\theta,\xi)
            -K(\theta_{c_2}, \xi) \sin (\theta_{c_2}) \parD{T}{\theta}(\theta_{c_2})\\
            &-T_c \sin (\theta_{c_2}) \lim_{\theta  \rightarrow \theta_{c_2}^+}\parD{K}{\theta}(\theta, \xi)
    \end{aligned}\\
    &\begin{aligned}
        T(\xi)  &= \int_{\theta_{l_2}}^{\pi}d\theta \: \sin\theta \: K(\theta, \xi) h_1(\theta)
            -K(\theta_{l_2}, \xi) \sin (\theta_{l_2}) \parD{T}{\theta}(\theta_{l_2})\\
            &+T(\theta_{l_2}) \sin (\theta_{l_2}) \lim_{\theta  \rightarrow \theta_{l_2}^+}\parD{K}{\theta}(\theta, \xi)
    \end{aligned}
    \label{eqn: 2 critical latitudes asymmetrical general integral relation applied region 5}
\end{align}
The associated BIEs are:
\begin{align}
    &\begin{aligned}
        T(0)  &= \int_0^{\theta_{c_1}}d\theta \: \sin\theta \: K(\theta, 0) h_2(\theta)
            + K(\theta_{c_1}, 0) \sin(\theta_{c_1}) \parD{T}{\theta}(\theta_{c_1})
            \\&+ \sin(\theta_{c_1}) \lim_{\xi  \rightarrow 0^+}\lim_{\theta  \rightarrow \theta_{c_1}^-} \parD{K}{\theta}(\theta,\xi)
    \end{aligned}
    \label{eqn: 2 critical latitudes asymmetrical BIE 1}\\
    &\begin{aligned}
        -1  &= \int_0^{\theta_{c_1}}d\theta \: \sin\theta \: K(\theta, \theta_{c_1}) h_2(\theta)
            + \lim_{\xi  \rightarrow \theta_{c_1}^-} K(\theta_{c_1}, \xi) \sin(\theta_{c_1}) \parD{T}{\theta}(\theta_{c_1})
            \\&+ \sin(\theta_{c_1}) \lim_{\xi  \rightarrow \theta_{c_1}^-}\lim_{\theta  \rightarrow \theta_{c_1}^-} \parD{K}{\theta}(\theta,\xi)
    \end{aligned}\\
    &\begin{aligned}
        -1  &= \int_{\theta_{c_1}}^{\theta_{l_1}}d\theta \: \sin\theta \: K(\theta, \theta_{c_1}) h_1(\theta)
            + K(\theta_{l_1}, \theta_{c_1}) \sin(\theta_{l_1}) \parD{T}{\theta}(\theta_{l_1})\\
            &- T(\theta_{l_1}) \sin(\theta_{l_1}) \lim_{\xi  \rightarrow \theta_{c_1}^+} \lim_{\theta  \rightarrow \theta_{l_1}^-}\parD{K}{\theta}(\theta,\xi)
            -\lim_{\xi  \rightarrow \theta_{c_1}^+}K(\theta_{c_1}, \xi) \sin (\theta_{c_1}) \parD{T}{\theta}(\theta_{c_1})\\
            &- \sin (\theta_{c_1})\lim_{\xi  \rightarrow \theta_{c_1}^+} \lim_{\theta  \rightarrow \theta_{c_1}^+}\parD{K}{\theta}(\theta, \xi)
    \end{aligned}
\end{align}
\begin{align}
    &\begin{aligned}
        T(\theta_{l_1})  &= \int_{\theta_{c_1}}^{\theta_{l_1}}d\theta \: \sin\theta \: K(\theta, \theta_{l_1}) h_1(\theta)
            + \lim_{\xi  \rightarrow \theta_{l_1}^-}K(\theta_{l_1}, \xi) \sin(\theta_{l_1}) \parD{T}{\theta}(\theta_{l_1})\\
            &- T(\theta_{l_1}) \sin(\theta_{l_1}) \lim_{\xi  \rightarrow \theta_{l_1}^-} \lim_{\theta  \rightarrow \theta_{l_1}^-}\parD{K}{\theta}(\theta,\xi)
            -K(\theta_{c_1}, \theta_{l_1}) \sin (\theta_{c_1}) \parD{T}{\theta}(\theta_{c_1})\\
            &- \sin (\theta_{c_1}) \lim_{\xi  \rightarrow \theta_{l_1}^-} \lim_{\theta  \rightarrow \theta_{c_1}^+}\parD{K}{\theta}(\theta, \xi)
    \end{aligned}\\
    &\begin{aligned}
        T(\theta_{l_1})  &= \int_{\theta_{l_1}}^{\theta_{c_2}}d\theta \: \sin\theta \: K(\theta, \theta_{l_1}) h_4(\theta)
            + K(\theta_{c_2}, \theta_{l_1}) \sin(\theta_{c_2}) \parD{T}{\theta}(\theta_{c_2})\\
            &+ T_c \sin(\theta_{c_2}) \lim_{\xi  \rightarrow \theta_{l_1}^+}\lim_{\theta  \rightarrow \theta_{l_2}^-}\parD{K}{\theta}(\theta,\xi)
            -\lim_{\xi  \rightarrow \theta_{l_1}^+}K(\theta_{l_1}, \xi) \sin (\theta_{l_1}) \parD{T}{\theta}(\theta_{l_1})\\
            &+T(\theta_{l_1}) \sin (\theta_{l_1}) \lim_{\xi  \rightarrow \theta_{l_1}^+} \lim_{\theta \rightarrow \theta_{l_1}^+}\parD{K}{\theta}(\theta, \xi)
    \end{aligned}\\
    &\begin{aligned}
        -T_c  &= \int_{\theta_{l_1}}^{\theta_{c_2}}d\theta \: \sin\theta \: K(\theta, \theta_{c_2}) h_4(\theta)
            + \lim_{\xi  \rightarrow \theta_{c_2}^-} K(\theta_{c_2}, \xi) \sin(\theta_{c_2}) \parD{T}{\theta}(\theta_{c_2})\\
            &+ T_c \sin(\theta_{c_2}) \lim_{\xi  \rightarrow \theta_{c_2}^-} \lim_{\theta  \rightarrow \theta_{c_2}^-}\parD{K}{\theta}(\theta,\xi)
            -K(\theta_{l_1}, \theta_{c_2}) \sin (\theta_{l_1}) \parD{T}{\theta}(\theta_{l_1})\\
            &+T(\theta_{l_1}) \sin (\theta_{l_1}) \lim_{\xi  \rightarrow \theta_{c_2}^-} \lim_{\theta  \rightarrow \theta_{l_1}^+}\parD{K}{\theta}(\theta, \xi)
    \end{aligned}\\ 
    &\begin{aligned}
        -T_c  &= \int_{\theta_{c_2}}^{\theta_{l_2}}d\theta \: \sin\theta \: K(\theta, \theta_{c_2}) h_3(\theta)
            + K(\theta_{l_2}, \theta_{c_2}) \sin(\theta_{l_2}) \parD{T}{\theta}(\theta_{l_2})\\
            &- T(\theta_{l_2}) \sin(\theta_{l_2})\lim_{\xi  \rightarrow \theta_{c_2}^+} \lim_{\theta  \rightarrow \theta_{l_2}^-}\parD{K}{\theta}(\theta,\xi)
            -\lim_{\xi  \rightarrow \theta_{c_2}^+}K(\theta_{c_2}, \xi) \sin (\theta_{c_2}) \parD{T}{\theta}(\theta_{c_2})\\
            &-T_c \sin (\theta_{c_2})\lim_{\xi  \rightarrow \theta_{c_2}^+} \lim_{\theta  \rightarrow \theta_{c_2}^+}\parD{K}{\theta}(\theta, \xi)
    \end{aligned}\\
    &\begin{aligned}
        T(\theta_{l_2})  &= \int_{\theta_{c_2}}^{\theta_{l_2}}d\theta \: \sin\theta \: K(\theta, \theta_{l_2}) h_3(\theta)
            + \lim_{\xi  \rightarrow \theta_{l_2}^-}K(\theta_{l_2}, \xi) \sin(\theta_{l_2}) \parD{T}{\theta}(\theta_{l_2})\\
            &- T(\theta_{l_2}) \sin(\theta_{l_2}) \lim_{\xi  \rightarrow \theta_{l_2}^-} \lim_{\theta  \rightarrow \theta_{l_2}^-}\parD{K}{\theta}(\theta,\xi)
            -K(\theta_{c_2}, \theta_{l_2}) \sin (\theta_{c_2}) \parD{T}{\theta}(\theta_{c_2})\\
            &-T_c \sin (\theta_{c_2}) \lim_{\xi  \rightarrow \theta_{l_2}^-} \lim_{\theta  \rightarrow \theta_{c_2}^+}\parD{K}{\theta}(\theta, \xi)
    \end{aligned}\\
    &\begin{aligned}
        T(\theta_{l_2})  &= \int_{\theta_{l_2}}^{\pi}d\theta \: \sin\theta \: K(\theta, \theta_{l_2}) h_1(\theta)
            - \lim_{\xi  \rightarrow \theta_{l_2}^+} K(\theta_{l_2}, \xi) \sin (\theta_{l_2}) \parD{T}{\theta}(\theta_{l_2})\\
            &+T(\theta_{l_2}) \sin (\theta_{l_2}) \lim_{\xi  \rightarrow \theta_{l_2}^+} \lim_{\theta  \rightarrow \theta_{l_2}^+}\parD{K}{\theta}(\theta, \xi)
    \end{aligned}
    \\
    &\begin{aligned}
        T(\pi)  &= \int_{\theta_{l_2}}^{\pi}d\theta \: \sin\theta \: K(\theta, \pi) h_1(\theta)
            -K(\theta_{l_2}, \pi) \sin (\theta_{l_2}) \parD{T}{\theta}(\theta_{l_2})\\
            &+T(\theta_{l_2}) \sin (\theta_{l_2}) \lim_{\xi  \rightarrow \pi^-} \lim_{\theta  \rightarrow \theta_{l_2}^+}\parD{K}{\theta}(\theta, \xi)
    \end{aligned}
    \label{eqn: 2 critical latitudes asymmetrical BIE 10}
\end{align}

\vsone
\subsection{Four Critical Latitudes}
\label{appendix: The case of 4 critical latitudes}
\setcounter{region}{0} 
States with four critical latitudes are possible for the meridionally symmetric continent configuration. As states with four critical latitudes only appear given the meridional symmetry, it is advisable to apply an assumption of symmetry to save on computational cost. The model's inherent symmetry ensures a symmetric temperature profiles around the equator. Consequently, the BIM is applied on the truncated domain $\theta \in [0, \frac{\pi}{2}]$. A symmetric solution implies that $T(\theta)=T(\pi -\theta)\implies T'(\theta)=-T'(\pi -\theta)$. At $\theta=\frac{\pi}{2}$, it is evident that $T'(\frac{\pi}{2}) = 0$. Four critical latitudes are realized by a partial ice cover on the ocean and the continent in both hemispheres, therefore the truncated domain is partitioned into the four sub-domains:
\begin{region}
    \theta \in (0, \theta_{c_1})
\end{region}
\begin{region}
    \theta \in (\theta_{c_1},\theta_{l_1})
\end{region}
\begin{region}
    \theta \in (\theta_{l_1}, \theta_{c_2})
\end{region}
\begin{region}
    \theta \in (\theta_{c_2}, \frac{\pi}{2})
\end{region}
Here the ice cover may be in region (\ref{region 3}) and region (\ref{region 4}) is ice free, or vice versa. Consider the former fist: The solution within the regions are (in respective order):
\begin{align}
    &\begin{aligned}
        T(\xi)  &= \int_0^{\theta_{c_1}}d\theta \: \sin\theta \: K(\theta, \xi) h_2(\theta)
            + K(\theta_{c_1}, \xi) \sin(\theta_{c_1}) \parD{T}{\theta}(\theta_{c_1})
            \\&+ \sin(\theta_{c_1}) \lim_{\theta  \rightarrow \theta_{c_1}^-} \parD{K}{\theta}(\theta,\xi)
    \end{aligned}
    \label{eqn: 4 critical latitudes general integral relation applied region 1 symmetry}\\
    &\begin{aligned}
        T(\xi)  &= \int_{\theta_{c_1}}^{\theta_{l_1}}d\theta \: \sin\theta \: K(\theta, \xi) h_1(\theta)
            + K(\theta_{l_1}, \xi) \sin(\theta_{l_1}) \parD{T}{\theta}(\theta_{l_1})\\
            &- T(\theta_{l_1}) \sin(\theta_{l_1}) \lim_{\theta  \rightarrow \theta_{l_1}^-}\parD{K}{\theta}(\theta,\xi)
            -K(\theta_{c_1}, \xi) \sin (\theta_{c_1}) \parD{T}{\theta}(\theta_{c_1})\\
            &- \sin (\theta_{c_1}) \lim_{\theta  \rightarrow \theta_{c_1}^+}\parD{K}{\theta}(\theta, \xi)
    \end{aligned}
    \label{eqn: 4 critical latitudes general integral relation applied region 2 symmetry}\\
    &\begin{aligned}
        T(\xi)  &= \int_{\theta_{l_1}}^{\theta_{c_2}}d\theta \: \sin\theta \: K(\theta, \xi) h_4(\theta)
            + K(\theta_{c_2}, \xi) \sin(\theta_{c_2}) \parD{T}{\theta}(\theta_{c_2})\\
            &+ T_c \sin(\theta_{c_2}) \lim_{\theta  \rightarrow \theta_{c_2}^-}\parD{K}{\theta}(\theta,\xi)
            -K(\theta_{l_1}, \xi) \sin (\theta_{l_1}) \parD{T}{\theta}(\theta_{l_1})\\
            &+T(\theta_{l_1}) \sin (\theta_{l_1}) \lim_{\theta  \rightarrow \theta_{l_1}^+}\parD{K}{\theta}(\theta, \xi)
    \end{aligned}
    \label{eqn: 4 critical latitudes general integral relation applied region 3 symmetry}\\
    &\begin{aligned}
        T(\xi)  &= \int_{\theta_{c_2}}^{\frac{\pi}{2}}d\theta \: \sin\theta \: K(\theta, \xi) h_3(\theta)
            - T(\frac{\pi}{2}) \lim_{\theta  \rightarrow \frac{\pi}{2}^-}\parD{K}{\theta}(\theta,\xi)
            \\&-K(\theta_{c_2}, \xi) \sin (\theta_{c_2}) \parD{T}{\theta}(\theta_{c_2})
            -T_c \sin (\theta_{c_2}) \lim_{\theta  \rightarrow \theta_{c_2}^+}\parD{K}{\theta}(\theta, \xi)
    \end{aligned}
    \label{eqn: 4 critical latitudes general integral relation applied region 4 symmetry}
\end{align}
The associated BIEs are:
\begin{align}
    &\begin{aligned}
        T(0)  &= \int_0^{\theta_{c_1}}d\theta \: \sin\theta \: K(\theta, 0) h_2(\theta)
            + K(\theta_{c_1}, 0) \sin(\theta_{c_1}) \parD{T}{\theta}(\theta_{c_1})
            \\&+ \sin(\theta_{c_1}) \lim_{\xi  \rightarrow 0^+}\lim_{\theta  \rightarrow \theta_{c_1}^-} \parD{K}{\theta}(\theta,\xi)
    \end{aligned}
    \label{eqn: 4 critical latitudes BIE 1 symmetry}\\
    &\begin{aligned}
        -1  &= \int_0^{\theta_{c_1}}d\theta \: \sin\theta \: K(\theta, \theta_{c_1}) h_2(\theta)
            + \lim_{\xi  \rightarrow \theta_{c_1}^-} K(\theta_{c_1}, \xi) \sin(\theta_{c_1}) \parD{T}{\theta}(\theta_{c_1})
            \\&+ \sin(\theta_{c_1}) \lim_{\xi  \rightarrow \theta_{c_1}^-}\lim_{\theta  \rightarrow \theta_{c_1}^-} \parD{K}{\theta}(\theta,\xi)
    \end{aligned}
    \label{eqn: 4 critical latitudes BIE 2 symmetry}\\
    &\begin{aligned}
        -1  &= \int_{\theta_{c_1}}^{\theta_{l_1}}d\theta \: \sin\theta \: K(\theta, \theta_{c_1}) h_1(\theta)
            + K(\theta_{l_1}, \theta_{c_1}) \sin(\theta_{l_1}) \parD{T}{\theta}(\theta_{l_1})\\
            &- T(\theta_{l_1}) \sin(\theta_{l_1}) \lim_{\xi  \rightarrow \theta_{c_1}^+} \lim_{\theta  \rightarrow \theta_{l_1}^-}\parD{K}{\theta}(\theta,\xi)
            -\lim_{\xi  \rightarrow \theta_{c_1}^+}K(\theta_{c_1}, \xi) \sin (\theta_{c_1}) \parD{T}{\theta}(\theta_{c_1})\\
            &- \sin (\theta_{c_1})\lim_{\xi  \rightarrow \theta_{c_1}^+} \lim_{\theta  \rightarrow \theta_{c_1}^+}\parD{K}{\theta}(\theta, \xi)
    \end{aligned}
    \label{eqn: 4 critical latitudes BIE 3 symmetry}\\
    &\begin{aligned}
        T(\theta_{l_1})  &= \int_{\theta_{c_1}}^{\theta_{l_1}}d\theta \: \sin\theta \: K(\theta, \theta_{l_1}) h_1(\theta)
            + \lim_{\xi  \rightarrow \theta_{l_1}^-}K(\theta_{l_1}, \xi) \sin(\theta_{l_1}) \parD{T}{\theta}(\theta_{l_1})\\
            &- T(\theta_{l_1}) \sin(\theta_{l_1}) \lim_{\xi  \rightarrow \theta_{l_1}^-} \lim_{\theta  \rightarrow \theta_{l_1}^-}\parD{K}{\theta}(\theta,\xi)
            -K(\theta_{c_1}, \theta_{l_1}) \sin (\theta_{c_1}) \parD{T}{\theta}(\theta_{c_1})\\
            &- \sin (\theta_{c_1}) \lim_{\xi  \rightarrow \theta_{l_1}^-} \lim_{\theta  \rightarrow \theta_{c_1}^+}\parD{K}{\theta}(\theta, \xi)
    \end{aligned}
    \label{eqn: 4 critical latitudes BIE 4 symmetry}\\
    &\begin{aligned}
        T(\theta_{l_1})  &= \int_{\theta_{l_1}}^{\theta_{c_2}}d\theta \: \sin\theta \: K(\theta, \theta_{l_1}) h_4(\theta)
            + K(\theta_{c_2}, \theta_{l_1}) \sin(\theta_{c_2}) \parD{T}{\theta}(\theta_{c_2})\\
            &+ T_c \sin(\theta_{c_2}) \lim_{\xi  \rightarrow \theta_{l_1}^+}\lim_{\theta  \rightarrow \theta_{l_2}^-}\parD{K}{\theta}(\theta,\xi)
            -\lim_{\xi  \rightarrow \theta_{l_1}^+}K(\theta_{l_1}, \xi) \sin (\theta_{l_1}) \parD{T}{\theta}(\theta_{l_1})\\
            &+T(\theta_{l_1}) \sin (\theta_{l_1}) \lim_{\xi  \rightarrow \theta_{l_1}^+} \lim_{\theta \rightarrow \theta_{l_1}^+}\parD{K}{\theta}(\theta, \xi)
    \end{aligned}
    \label{eqn: 4 critical latitudes BIE 5 symmetry}\\
    &\begin{aligned}
        -T_c  &= \int_{\theta_{l_1}}^{\theta_{c_2}}d\theta \: \sin\theta \: K(\theta, \theta_{c_2}) h_4(\theta)
            + \lim_{\xi  \rightarrow \theta_{c_2}^-} K(\theta_{c_2}, \xi) \sin(\theta_{c_2}) \parD{T}{\theta}(\theta_{c_2})\\
            &+ T_c \sin(\theta_{c_2}) \lim_{\xi  \rightarrow \theta_{c_2}^-} \lim_{\theta  \rightarrow \theta_{c_2}^-}\parD{K}{\theta}(\theta,\xi)
            -K(\theta_{l_1}, \theta_{c_2}) \sin (\theta_{l_1}) \parD{T}{\theta}(\theta_{l_1})\\
            &+T(\theta_{l_1}) \sin (\theta_{l_1}) \lim_{\xi  \rightarrow \theta_{c_2}^-} \lim_{\theta  \rightarrow \theta_{l_1}^+}\parD{K}{\theta}(\theta, \xi)
    \end{aligned}
    \label{eqn: 4 critical latitudes BIE 6 symmetry}\\
    &\begin{aligned}
        -T_c  &= \int_{\theta_{c_2}}^{\frac{\pi}{2}}d\theta \: \sin\theta \: K(\theta, \theta_{c_2}) h_3(\theta)
            - T(\frac{\pi}{2})\lim_{\xi  \rightarrow \theta_{c_2}^+} \lim_{\theta  \rightarrow \frac{\pi}{2}^-}\parD{K}{\theta}(\theta,\xi)
            \\&-\lim_{\xi  \rightarrow \theta_{c_2}^+}K(\theta_{c_2}, \xi) \sin (\theta_{c_2}) \parD{T}{\theta}(\theta_{c_2})
            -T_c \sin (\theta_{c_2}) \lim_{\xi  \rightarrow \theta_{c_2}^+}\lim_{\theta  \rightarrow \theta_{c_2}^+}\parD{K}{\theta}(\theta, \xi)
    \end{aligned}
    \label{eqn: 4 critical latitudes BIE 7 symmetry}
\end{align}
\begin{align}
    &\begin{aligned}
          T(\frac{\pi}{2})  &= \int_{\theta_{c_2}}^{\frac{\pi}{2}}d\theta \: \sin\theta \: K(\theta, \frac{\pi}{2}) h_3(\theta)
            - T(\frac{\pi}{2})\lim_{\xi  \rightarrow \frac{\pi}{2}^-} \lim_{\theta  \rightarrow \frac{\pi}{2}^-}\parD{K}{\theta}(\theta,\xi)
            \\&-K(\theta_{c_2}, \frac{\pi}{2}) \sin (\theta_{c_2}) \parD{T}{\theta}(\theta_{c_2})
            -T_c \sin (\theta_{c_2}) \lim_{\xi  \rightarrow \frac{\pi}{2}^-}\lim_{\theta  \rightarrow \theta_{c_2}^+}\parD{K}{\theta}(\theta, \xi).
    \end{aligned}
    \label{eqn: 4 critical latitudes BIE 8 symmetry}
\end{align}
\vsone

If region (\ref{region 3}) is ice free and region (\ref{region 4}) is ice covered, the equations from above must be changed accordingly. The solution is given by (\ref{eqn: 4 critical latitudes general integral relation applied region 1 symmetry})-(\ref{eqn: 4 critical latitudes general integral relation applied region 4 symmetry}), but changing $h_4 \to h_3$ in (\ref{eqn: 4 critical latitudes general integral relation applied region 3 symmetry}) and $h_3 \to h_4$ in (\ref{eqn: 4 critical latitudes general integral relation applied region 4 symmetry}). BIEs are (\ref{eqn: 4 critical latitudes BIE 1 symmetry})-(\ref{eqn: 4 critical latitudes BIE 8 symmetry}), interchanging $h_4 \to h_3$ in (\ref{eqn: 4 critical latitudes BIE 5 symmetry}) and (\ref{eqn: 4 critical latitudes BIE 6 symmetry}), as well as $h_3 \to h_4$ in (\ref{eqn: 4 critical latitudes BIE 7 symmetry}) and (\ref{eqn: 4 critical latitudes BIE 8 symmetry}).

\vsone
\subsection{Six Critical Latitudes}
\label{appendix: The case of 6 critical latitudes}
\setcounter{region}{0} 
States with six critical latitudes are possible with the meridionally symmetric continent configuration studied. These states are realized by two zonal strips of ice cover on the continent, forming four critical latitudes on the continent and six in total. Again, aa assumption of symmetry is applied and the truncated domain is partitioned into the five sub-domains:
\begin{region}
    \theta \in (0, \theta_{c_1})
\end{region}
\begin{region}
    \theta \in (\theta_{c_1},\theta_{l_1})
\end{region}
\begin{region}
    \theta \in (\theta_{l_1}, \theta_{c_2})
\end{region}
\begin{region}
    \theta \in (\theta_{c_2}, \theta_{c_3})
\end{region}
\begin{region}
    \theta \in (\theta_{c_3}, \frac{\pi}{2})
\end{region}
The ice distribution within the continent may take on two forms: Ice in region (\ref{region 4}) or ice in region (\ref{region 3}) and (\ref{region 5}). For the former the solution within the five regions are given by (in order):
\begin{align}
    &\begin{aligned}
        T(\xi)  &= \int_0^{\theta_{c_1}}d\theta \: \sin\theta \: K(\theta, \xi) h_2(\theta)
            + K(\theta_{c_1}, \xi) \sin(\theta_{c_1}) \parD{T}{\theta}(\theta_{c_1})
            \\&+ \sin(\theta_{c_1}) \lim_{\theta  \rightarrow \theta_{c_1}^-} \parD{K}{\theta}(\theta,\xi)
    \end{aligned}
    \label{eqn: 6 critical latitudes general integral relation applied region 1 symmetry}\\
    &\begin{aligned}
        T(\xi)  &= \int_{\theta_{c_1}}^{\theta_{l_1}}d\theta \: \sin\theta \: K(\theta, \xi) h_1(\theta)
            + K(\theta_{l_1}, \xi) \sin(\theta_{l_1}) \parD{T}{\theta}(\theta_{l_1})\\
            &- T(\theta_{l_1}) \sin(\theta_{l_1}) \lim_{\theta  \rightarrow \theta_{l_1}^-}\parD{K}{\theta}(\theta,\xi)
            -K(\theta_{c_1}, \xi) \sin (\theta_{c_1}) \parD{T}{\theta}(\theta_{c_1})\\
            &- \sin (\theta_{c_1}) \lim_{\theta  \rightarrow \theta_{c_1}^+}\parD{K}{\theta}(\theta, \xi)
    \end{aligned}
    \label{eqn: 6 critical latitudes general integral relation applied region 2 symmetry}
\end{align}
\begin{align}
    &\begin{aligned}
        T(\xi)  &= \int_{\theta_{l_1}}^{\theta_{c_2}}d\theta \: \sin\theta \: K(\theta, \xi) h_3(\theta)
            + K(\theta_{c_2}, \xi) \sin(\theta_{c_2}) \parD{T}{\theta}(\theta_{c_2})\\
            &+ T_c \sin(\theta_{c_2}) \lim_{\theta  \rightarrow \theta_{c_2}^-}\parD{K}{\theta}(\theta,\xi)
            -K(\theta_{l_1}, \xi) \sin (\theta_{l_1}) \parD{T}{\theta}(\theta_{l_1})\\
            &+T(\theta_{l_1}) \sin (\theta_{l_1}) \lim_{\theta  \rightarrow \theta_{l_1}^+}\parD{K}{\theta}(\theta, \xi)
    \end{aligned}
    \label{eqn: 6 critical latitudes general integral relation applied region 3 symmetry}\\
        &\begin{aligned}
        T(\xi)  &= \int_{\theta_{c_2}}^{\theta_{c_3}}d\theta \: \sin\theta \: K(\theta, \xi) h_4(\theta)
            + K(\theta_{c_3}, \xi) \sin(\theta_{c_3}) \parD{T}{\theta}(\theta_{c_3})\\
            &+ T_c \sin(\theta_{c_3}) \lim_{\theta  \rightarrow \theta_{c_3}^-}\parD{K}{\theta}(\theta,\xi)
            -K(\theta_{c_2}, \xi) \sin (\theta_{c_2}) \parD{T}{\theta}(\theta_{c_2})\\
            &-T_c \sin (\theta_{c_2}) \lim_{\theta  \rightarrow \theta_{c_2}^+}\parD{K}{\theta}(\theta, \xi)
    \end{aligned}
    \label{eqn: 6 critical latitudes general integral relation applied region 4 symmetry}\\
    &\begin{aligned}
        T(\xi)  &= \int_{\theta_{c_3}}^{\frac{\pi}{2}}d\theta \: \sin\theta \: K(\theta, \xi) h_3(\theta)
            - T(\frac{\pi}{2}) \lim_{\theta  \rightarrow \frac{\pi}{2}^-}\parD{K}{\theta}(\theta,\xi)
            \\&-K(\theta_{c_3}, \xi) \sin (\theta_{c_3}) \parD{T}{\theta}(\theta_{c_3})
            -T_c \sin (\theta_{c_3}) \lim_{\theta  \rightarrow \theta_{c_3}^+}\parD{K}{\theta}(\theta, \xi)
    \end{aligned}
    \label{eqn: 6 critical latitudes general integral relation applied region 5 symmetry}
\end{align}
The associated BIEs are:
\begin{align}
    &\begin{aligned}
        T(0)  &= \int_0^{\theta_{c_1}}d\theta \: \sin\theta \: K(\theta, 0) h_2(\theta)
            + K(\theta_{c_1}, 0) \sin(\theta_{c_1}) \parD{T}{\theta}(\theta_{c_1})
            \\&+ \sin(\theta_{c_1}) \lim_{\xi  \rightarrow 0^+}\lim_{\theta  \rightarrow \theta_{c_1}^-} \parD{K}{\theta}(\theta,\xi)
    \end{aligned}
    \label{eqn: 6 critical latitudes BIE 1 symmetry}\\
    &\begin{aligned}
        -1  &= \int_0^{\theta_{c_1}}d\theta \: \sin\theta \: K(\theta, \theta_{c_1}) h_2(\theta)
            + \lim_{\xi  \rightarrow \theta_{c_1}^-} K(\theta_{c_1}, \xi) \sin(\theta_{c_1}) \parD{T}{\theta}(\theta_{c_1})
            \\&+ \sin(\theta_{c_1}) \lim_{\xi  \rightarrow \theta_{c_1}^-}\lim_{\theta  \rightarrow \theta_{c_1}^-} \parD{K}{\theta}(\theta,\xi)
    \end{aligned}
    \label{eqn: 6 critical latitudes BIE 2 symmetry}\\
    &\begin{aligned}
        -1  &= \int_{\theta_{c_1}}^{\theta_{l_1}}d\theta \: \sin\theta \: K(\theta, \theta_{c_1}) h_1(\theta)
            + K(\theta_{l_1}, \theta_{c_1}) \sin(\theta_{l_1}) \parD{T}{\theta}(\theta_{l_1})\\
            &- T(\theta_{l_1}) \sin(\theta_{l_1}) \lim_{\xi  \rightarrow \theta_{c_1}^+} \lim_{\theta  \rightarrow \theta_{l_1}^-}\parD{K}{\theta}(\theta,\xi)
            -\lim_{\xi  \rightarrow \theta_{c_1}^+}K(\theta_{c_1}, \xi) \sin (\theta_{c_1}) \parD{T}{\theta}(\theta_{c_1})\\
            &- \sin (\theta_{c_1})\lim_{\xi  \rightarrow \theta_{c_1}^+} \lim_{\theta  \rightarrow \theta_{c_1}^+}\parD{K}{\theta}(\theta, \xi)
    \end{aligned}
    \label{eqn: 6 critical latitudes BIE 3 symmetry}\\
    &\begin{aligned}
        T(\theta_{l_1})  &= \int_{\theta_{c_1}}^{\theta_{l_1}}d\theta \: \sin\theta \: K(\theta, \theta_{l_1}) h_1(\theta)
            + \lim_{\xi  \rightarrow \theta_{l_1}^-}K(\theta_{l_1}, \xi) \sin(\theta_{l_1}) \parD{T}{\theta}(\theta_{l_1})\\
            &- T(\theta_{l_1}) \sin(\theta_{l_1}) \lim_{\xi  \rightarrow \theta_{l_1}^-} \lim_{\theta  \rightarrow \theta_{l_1}^-}\parD{K}{\theta}(\theta,\xi)
            -K(\theta_{c_1}, \theta_{l_1}) \sin (\theta_{c_1}) \parD{T}{\theta}(\theta_{c_1})\\
            &- \sin (\theta_{c_1}) \lim_{\xi  \rightarrow \theta_{l_1}^-} \lim_{\theta  \rightarrow \theta_{c_1}^+}\parD{K}{\theta}(\theta, \xi)
    \end{aligned}
    \label{eqn: 6 critical latitudes BIE 4 symmetry}
\end{align}
\begin{align}
    &\begin{aligned}
        T(\theta_{l_1})  &= \int_{\theta_{l_1}}^{\theta_{c_2}}d\theta \: \sin\theta \: K(\theta, \theta_{l_1}) h_3(\theta)
            + K(\theta_{c_2}, \theta_{l_1}) \sin(\theta_{c_2}) \parD{T}{\theta}(\theta_{c_2})\\
            &+ T_c \sin(\theta_{c_2}) \lim_{\xi  \rightarrow \theta_{l_1}^+}\lim_{\theta  \rightarrow \theta_{l_2}^-}\parD{K}{\theta}(\theta,\xi)
            -\lim_{\xi  \rightarrow \theta_{l_1}^+}K(\theta_{l_1}, \xi) \sin (\theta_{l_1}) \parD{T}{\theta}(\theta_{l_1})\\
            &+T(\theta_{l_1}) \sin (\theta_{l_1}) \lim_{\xi  \rightarrow \theta_{l_1}^+} \lim_{\theta \rightarrow \theta_{l_1}^+}\parD{K}{\theta}(\theta, \xi)
    \end{aligned}
    \label{eqn: 6 critical latitudes BIE 5 symmetry}\\
    &\begin{aligned}
        -T_c  &= \int_{\theta_{l_1}}^{\theta_{c_2}}d\theta \: \sin\theta \: K(\theta, \theta_{c_2}) h_3(\theta)
            + \lim_{\xi  \rightarrow \theta_{c_2}^-} K(\theta_{c_2}, \xi) \sin(\theta_{c_2}) \parD{T}{\theta}(\theta_{c_2})\\
            &+ T_c \sin(\theta_{c_2}) \lim_{\xi  \rightarrow \theta_{c_2}^-} \lim_{\theta  \rightarrow \theta_{c_2}^-}\parD{K}{\theta}(\theta,\xi)
            -K(\theta_{l_1}, \theta_{c_2}) \sin (\theta_{l_1}) \parD{T}{\theta}(\theta_{l_1})\\
            &+T(\theta_{l_1}) \sin (\theta_{l_1}) \lim_{\xi  \rightarrow \theta_{c_2}^-} \lim_{\theta  \rightarrow \theta_{l_1}^+}\parD{K}{\theta}(\theta, \xi)
    \end{aligned}
    \label{eqn: 6 critical latitudes BIE 6 symmetry}\\
    &\begin{aligned}
        -T_c  &= \int_{\theta_{c_2}}^{\theta_{c_3}}d\theta \: \sin\theta \: K(\theta, \theta_{c_2}) h_4(\theta)
            + K(\theta_{c_3}, \theta_{c_2}) \sin(\theta_{c_3}) \parD{T}{\theta}(\theta_{c_3})\\
            &+ T_c \sin(\theta_{c_3}) \lim_{\xi  \rightarrow \theta_{c_2}^+}\lim_{\theta  \rightarrow \theta_{c_3}^-}\parD{K}{\theta}(\theta,\xi)
            - \lim_{\xi  \rightarrow \theta_{c_2}^+} K(\theta_{c_2}, \xi) \sin (\theta_{c_2}) \parD{T}{\theta}(\theta_{c_2})\\
            &-T_c \sin (\theta_{c_2})\lim_{\xi  \rightarrow \theta_{c_2}^-} \lim_{\theta \rightarrow \theta_{c_2}^+}\parD{K}{\theta}(\theta, \xi)
    \end{aligned}
    \label{eqn: 6 critical latitudes BIE 7 symmetry}\\
    &\begin{aligned}
        -T_c  &= \int_{\theta_{c_2}}^{\theta_{c_3}}d\theta \: \sin\theta \: K(\theta, \theta_{c_3}) h_4(\theta)
            + \lim_{\xi  \rightarrow \theta_{c_3}^-}K(\theta_{c_3}, \xi) \sin(\theta_{c_3}) \parD{T}{\theta}(\theta_{c_3})\\
            &+ T_c \sin(\theta_{c_3}) \lim_{\xi  \rightarrow \theta_{c_2}^-}\lim_{\theta  \rightarrow \theta_{c_3}^-}\parD{K}{\theta}(\theta,\xi)
            -K(\theta_{c_2}, \theta_{c_3}) \sin (\theta_{c_2}) \parD{T}{\theta}(\theta_{c_2})\\
            &-T_c \sin (\theta_{c_2}) \lim_{\xi  \rightarrow \theta_{c_2}^-}\lim_{\theta  \rightarrow \theta_{c_2}^+}\parD{K}{\theta}(\theta, \xi)
    \end{aligned}
    \label{eqn: 6 critical latitudes BIE 8 symmetry}\\
    &\begin{aligned}
        -T_c  &= \int_{\theta_{c_3}}^{\frac{\pi}{2}}d\theta \: \sin\theta \: K(\theta, \theta_{c_3}) h_3(\theta)
            - T(\frac{\pi}{2})\lim_{\xi  \rightarrow \theta_{c_3}^+} \lim_{\theta  \rightarrow \frac{\pi}{2}^-}\parD{K}{\theta}(\theta,\xi)
            \\&-\lim_{\xi  \rightarrow \theta_{c_3}^+}K(\theta_{c_3}, \xi) \sin (\theta_{c_3}) \parD{T}{\theta}(\theta_{c_3})
            -T_c \sin (\theta_{c_3}) \lim_{\xi  \rightarrow \theta_{c_3}^+}\lim_{\theta  \rightarrow \theta_{c_3}^+}\parD{K}{\theta}(\theta, \xi)
    \end{aligned}
    \label{eqn: 6 critical latitudes BIE 9 symmetry}\\
    &\begin{aligned}
          T(\frac{\pi}{2})  &= \int_{\theta_{c_3}}^{\frac{\pi}{2}}d\theta \: \sin\theta \: K(\theta, \frac{\pi}{2}) h_3(\theta)
            - T(\frac{\pi}{2})\lim_{\xi  \rightarrow \frac{\pi}{2}^-} \lim_{\theta  \rightarrow \frac{\pi}{2}^-}\parD{K}{\theta}(\theta,\xi)
            \\&-K(\theta_{c_3}, \frac{\pi}{2}) \sin (\theta_{c_3}) \parD{T}{\theta}(\theta_{c_3})
            -T_c \sin (\theta_{c_3}) \lim_{\xi  \rightarrow \frac{\pi}{2}^-}\lim_{\theta  \rightarrow \theta_{c_3}^+}\parD{K}{\theta}(\theta, \xi)
    \end{aligned}
    \label{eqn: 6 critical latitudes BIE 10 symmetry}
\end{align}

\vsone
Solutions corresponding to an ice distribution with ice present in region (\ref{region 3}) and (\ref{region 5}) are given by (\ref{eqn: 6 critical latitudes general integral relation applied region 1 symmetry})-(\ref{eqn: 6 critical latitudes general integral relation applied region 1 symmetry}), but changing $h_3 \to h_4$ in (\ref{eqn: 6 critical latitudes general integral relation applied region 3 symmetry}) and (\ref{eqn: 6 critical latitudes general integral relation applied region 5 symmetry}), as well as $h_4 \to h_3$ in (\ref{eqn: 6 critical latitudes general integral relation applied region 4 symmetry}). BIEs are (\ref{eqn: 6 critical latitudes BIE 1 symmetry})-(\ref{eqn: 6 critical latitudes BIE 10 symmetry}), changing $h_3 \to h_4$ in (\ref{eqn: 6 critical latitudes BIE 5 symmetry}), (\ref{eqn: 6 critical latitudes BIE 6 symmetry}), (\ref{eqn: 6 critical latitudes BIE 9 symmetry}) and (\ref{eqn: 6 critical latitudes BIE 10 symmetry}), and $h_4 \to h_3$ in (\ref{eqn: 6 critical latitudes BIE 7 symmetry}) and (\ref{eqn: 6 critical latitudes BIE 8 symmetry}).
\end{appendices}

\clearpage
\addcontentsline{toc}{section}{References}
\printbibliography

\end{document}